\documentclass[moor]{informs1}              

\usepackage{algorithm}
\usepackage{algpseudocode}
\usepackage[euler]{textgreek}
\usepackage{enumitem}
\usepackage{relsize}
\usepackage{subfigure}
\usepackage{multirow}
\usepackage{capt-of}
\usepackage{bm}

\newcommand{\Real}{\mathbb{R}}
\newcommand{\Diag}{{\textrm{diag}}}
\newcommand{\st}{\mathrm{s.t.}}

\newcommand{\sfW}{\mathsf{N}}
\newcommand{\sfq}{\mathsf{d}}
\newcommand{\sfC}{\mathsf{C}}
\newcommand{\sfD}{\mathsf{D}}
\newcommand{\sfH}{\mathsf{H}}
\newcommand{\KL}{\mathsf{KL}}
\newcommand{\sfR}{\mathsf{R}}
\newcommand{\bfP}{{\mathbf{P}}}

\newcommand{\sfV}{{\mathsf{V}}}

\newcommand{\sfv}{{\mathsf{v}}}
\newcommand{\rmx}{{\mathrm{x}}}
\newcommand{\rmX}{{\mathrm{X}}}
\newcommand{\rmY}{{\mathrm{Y}}}
\newcommand{\rmW}{{\mathrm{W}}}
\newcommand{\calC}{\mathcal{C}}
\newcommand{\calH}{\mathcal{H}}
\newcommand{\calO}{\mathcal{O}}
\newcommand{\dom}{\mathrm{dom}}
\newcommand{\intr}{\mathrm{int}}
\newcommand{\blk}{{\textrm{blkdiag}}}
\newcommand{\eig}{{\textrm{eig}}}
\newcommand{\xUL}{y}
\newcommand{\yLL}{x}

\newcommand{\tphi}{\widetilde{\phi}}
\newcommand{\tlambda}{\widetilde{\lambda}}
\newcommand{\txteps}{\text{\textepsilon}}

\usepackage{natbib}
 \NatBibNumeric
 \def\bibfont{\small}%
 \bibpunct[, ]{[}{]}{,}{n}{}{,}%

\usepackage[colorlinks=true,breaklinks=true,bookmarks=true,urlcolor=blue,
     citecolor=blue,linkcolor=blue,bookmarksopen=false,draft=false]{hyperref}

\def\EMAIL#1{\href{mailto:#1}{#1}}

\TheoremsNumberedThrough     

\EquationsNumberedThrough    

\begin{document}


\RUNAUTHOR{Goyal and Lamperski}

\RUNTITLE{Bilevel optimization with traffic equilibrium constraints}

\TITLE{An algorithm for bilevel optimization with traffic equilibrium constraints: convergence rate analysis}

\ARTICLEAUTHORS{%
\AUTHOR{Akshit Goyal}
\AFF{Department of Industrial and Systems Engineering, University of Minnesota, Minneapolis, Minnesota 55455, \EMAIL{goyal080@umn.edu}}
\AUTHOR{Andrew Lamperski}
\AFF{Department of Electrical and Computer Engineering, University of Minnesota, Minneapolis, Minnesota 55455, \EMAIL{alampers@umn.edu}}
} 

\ABSTRACT{%
Bilevel optimization with traffic equilibrium constraints plays an important role in transportation planning and management problems such as traffic control, transport network design, and congestion pricing. In this paper, we consider a double-loop gradient-based algorithm to solve such bilevel problems and provide a non-asymptotic convergence guarantee of $\calO(K^{-1})+\calO(\lambda^D)$ where $K$, $D$ are respectively the number of upper- and lower-level iterations, and $0<\lambda<1$ is a constant. Compared to existing literature, which either provides asymptotic convergence or makes strong assumptions and requires a complex design of step sizes, we establish convergence for choice of simple constant step sizes and considering fewer assumptions. The analysis techniques in this paper use concepts from the field of robust control and can potentially serve as a guiding framework for analyzing more general bilevel optimization algorithms.
}%


\MSCCLASS{}
\ORMSCLASS{Primary: ; secondary: }
\HISTORY{}

\maketitle
%
\section{Introduction.}\label{sec:intro} 
Bilevel optimization with traffic equilibrium constraints (BOTEC), also called the Stackelberg congestion game by \citet{li2022differentiable}, is an important problem widely studied in the transportation literature.  It has applications in many areas of transportation planning and management, including area traffic control (\citet{chiou1999optimization,chiou2010efficient}), transport network design  (\citet{suwansirikul1987equilibrium, meng2001equivalent, chiou2005bilevel, josefsson2007sensitivity}), and congestion pricing  (\citet{verhoef2002second}). At the core of this problem is the concept of traffic equilibria, also called Wardrop equilibria, which refers to the Nash equilibrium of a routing or congestion game with a finite number of populations. In this game, each population (i.e., group of travelers) with travel demand between a unique pair of origin-destination nodes in a network competes to minimize their travel time. 
The equilibrium is then defined as the set of network flows that satisfy the travel demand and minimize the travel time of each population, 
such that there is no incentive for any population to change their route to improve their travel time. This equilibrium condition is often used as a constraint in bilevel optimization problems, where the upper-level problem seeks to optimize some system-level objective, subject to the condition that the network is in a traffic equilibrium state. The equilibrium solution in fact can be found by solving an equivalent optimization problem since the routing game is a potential game as argued by \citet{beckmann1956studies, rosenthal1973class} in which all populations together aim to minimize a potential function with each flow vector constrained to a simplex.
As a result, in this paper, we consider BOTEC problems of the form \eqref{eq:BOTEC_Intro} where an extra regularization term is added to the lower-level objective. This regularization makes the analysis more tractable.
\begin{subequations}\label{eq:BOTEC_Intro}
\begin{align}
    \min_{\xUL\in\calC} \ & F(\xUL) = f\left( h^*(\xUL),\xUL\right) \\
    \st \ & h^*(\xUL) \in \begin{aligned}[t]
    \underset{h=(h_i, i\in[\sfW])}{\arg\min}\;\; & g(h,\xUL) + \eta\ \overline{\psi}(h) \quad \st \;\; h_i\geq0,\ \mathbf{1}^\top h_i = 1\;\;\forall i\in[\sfW]  
    \end{aligned}\label{eq:LowerLevel}
\end{align}    
\end{subequations}
where $f(h,\xUL)$ and $g(h,\xUL)$ are respectively the upper- and lower-level objective functions, $\xUL$ is the upper-level decision variable, $h=(h_i, i\in[\sfW])$ is the lower-decision variable (where $h_i$ corresponds to the decision of $i$\textsuperscript{th} population), $h^*(\xUL)$ is the optimal lower-level solution for a given $\xUL$ and the regularizer $\overline\psi(h)$ is a strongly convex function composed of negative Shannon entropies.
Note that $\mathbf{1}$ is the vector of ones.

BOTEC is a non-convex problem and various gradient-based solution approaches have been proposed for different transportation applications in the literature. In context of transport network design problem, \citet{suwansirikul1987equilibrium} proposed a heuristic for solving BOTEC and \citet{meng2001equivalent} apply an augmented Lagrangian method with known local convergence. 
Extending these works, 
\citet{chiou2005bilevel} conduct extensive numerical experiments 
to demonstrate the computational efficacy of four different variants of gradient-based methods. All of the four approaches rely on obtaining the approximate gradient of $F(\xUL)$ using the sensitivity analysis of lower-level equilibrium network flows w.r.t upper-level decision as discussed by \citet{tobin1988sensitivity, patriksson2004sensitivity, josefsson2007sensitivity}. This approximate gradient is used in a gradient projection algorithm to update the upper-level decision, for which \citet{chiou1999optimization} and \citet{chiou2005bilevel} argue convergence to the KKT point. However, these works typically focus on asymptotic convergence without providing convergence rates.

Only recently, \citet{liu2022inducing} provide the convergence rate of a single-loop algorithm for solving BOTEC by extending the work of \citet{hong2023two} to a lower-level problem with simplex constraints. Furthermore, the recent work of \citet{li2023achieving} construct upper and lower bounding problems to bilevel programs with general equilibrium constraints, and bounds the rate of convergence to the optimal objective.
Compared to standard assumptions in bilevel literature, \citet{liu2022inducing} and \citet{li2023achieving} make additional assumptions on $F(\xUL)$ to be strongly convex. Further, \citet{liu2022inducing} assumes its gradient estimator (which uses implicit differentiation of lower-level dynamics) to be Lipshchitz continuous. Since $F(\xUL)$ depends on the optimal solution of lower-level problem \eqref{eq:LowerLevel}, it is unclear if the convexity of $F(\xUL)$ can be checked. Contrary to this, our work only makes smoothness and boundedness assumptions on $f$, $g$ where $f$ can be non-convex and $g$ is convex but makes no assumption on $F(\xUL)$. Moreover, the bound on gradient estimation error of $F(\xUL)$ is derived explicitly (as given in Lemma \ref{lemma:GradEstErrBnd-2}) without assuming Lipschitz continuity. 

\subsection{Our contributions.}
Overall, this paper makes the following contributions.
\begin{enumerate}
    \item Prove $\calO(K^{-1})+\calO(\lambda^D)$ convergence of the proposed algorithm for solving BOTEC to a stationary point where $K$, $D$ are respectively the number of upper- and lower-level iterations, $0<\lambda<1$ is a constant.
    \item As compared to analysis in \citet{liu2022inducing} which require complex design of step-sizes, we prove convergence by considering constant upper-level and lower-level step sizes satisfying appropriate conditions. Also, only empirical study by \citet{li2022differentiable} has been done on the convergence of routing game as well as BOTEC under fixed step-size whereas our work conducts theoretical analysis. 
    \item Introduce analysis techniques which use concepts from theory of robust control (\citet{skelton2013unified}). This provides a novel framework to analyze algorithms for a more general class of bilevel optimization problems in future.
\end{enumerate}

\subsection{Outline.}\label{outline1} 
The paper is outlined as follows. In Section \ref{sec:problem_setup}, we first introduce the problem in more detail followed by assumptions in Section \ref{sec:assumptions}
and relevant notations in \ref{sec:notations}. Section \ref{sec:main-results} highlights main results of the paper. In particular, Section \ref{sec:Algorithm} describes the proposed algorithm to solve problem \eqref{eq:BOTEC_Intro} and Section \ref{sec:ConvgRate} specify its convergence rate guarantee for both unconstrained \& constrained cases. Section \ref{sec:ExampleApps} specifies some of the key applications of BOTEC such as transport network design (Section \ref{sec:App1-NetowrkDesign}) and traffic control (Section \ref{sec:App2-TrafficControl}) which rely on solving traffic equilibirum problems explained in Section \ref{sec:TrafficEqbm}. In Section \ref{sec:NumericExprmnt} we conduct numerical experiments on benchmark Sioux Fall network and report computational results. Lastly, section 
\ref{sec:Proofs} details the proof of convergence rate results in Section \ref{sec:ConvgRate}. 

\section{Problem setup.}\label{sec:problem_setup} 
Consider problem \eqref{eq:BOTEC_Intro} with entropic regularization $\overline{\psi}_\sfH(\cdot)$ as follows with $\eta>0$
\begin{subequations}
\begin{align}
    \min_{\xUL\in\calC} \ & F(\xUL) = f\left( h^*(\xUL),\xUL\right) \\
    \st \ & h^*(\xUL) \in 
    \begin{aligned}[t]
    \underset{h=(h_i, i\in[\sfW])\in\calH}{\arg\min}\;\; & g^\eta(h,\xUL) = g(h,\xUL) + \eta\ \overline{\psi}_\sfH(h) 
    \end{aligned}\label{eq:LowerLevel-EntropyReg}
\end{align}    
\end{subequations}
where $\calC\subseteq \Real^{\sfq_u}$ is a closed convex set, $\calH=\calH_1\times\cdots\times\calH_\sfW$, $\calH_i = \big\{h_i\in\Real^{\sfq_\ell^i}| h_i\geq0,\ \mathbf{1}^\top h_i = 1\big\}$ for each $i=1,\hdots,\sfW$, and the regularizer $\overline{\psi}_\sfH(h)=\sum_{i=1}^\sfW \psi_\sfH(h_i)$ with $\psi_\sfH(h_i)=\sum_{j=1}^{\sfq_\ell^i} \left(h_{i,j} \ln(h_{i,j})-h_{i,j}\right)$ as the negative Shanon entropy. For any given $\xUL$, the function  $g(h,\xUL)$ is convex in $h$ and therefore, the lower-level objective $g^\eta(h,\xUL)$ is $\eta$--strongly convex in $h$ (see Appendix \ref{appendix:Prelim_MD}--(\ref{propty:StrongCnvx})). As a result, $h^*(\xUL)$ is unique for any given $\xUL$.

\subsection{Assumptions.}\label{sec:assumptions}
For vectors, $\|\cdot\|$ means the Euclidean norm and in case of matrices, $\|\cdot\|$ refers to the spectral norm.
\begin{assumption}
Denote $\bm{z}=(h,\xUL)$. For upper-level problem,  we assume
\begin{enumerate}[label=\normalfont(U\arabic*)]
    \item \label{assump:Lipschtz_ULGrad} The gradient $\nabla_{\bm{z}} f(\bm{z})$ is $L_f$-Lipschitz.
    $$\|\nabla_{\bm{z}} f(\bm{z})-\nabla_{\bm{z}} f(\bm{z}')\| \leq L_{f} \|\bm{z}-\bm{z}'\| \;\; \forall \bm{z}, \bm{z}'\in\calH\times\calC$$
    \item \label{assump:Bounded_ULGrad} The norm of gradient is bounded i.e. $\|\nabla_{\bm{z}} f(\bm{z})\|\leq \Omega_f\;\forall \bm{z}\in\calH\times\calC$.
\end{enumerate}
\end{assumption}

\begin{assumption} \label{assump:lower-level}
For lower-level problem, we assume for any given $\xUL$
    \begin{enumerate}[label=\normalfont(L\arabic*)]    
    \item \label{assump:Lipschtz_LLGrad} The gradient $\nabla_h g(h,\xUL)$ is $L_g$-Lipschitz w.r.t. $h$.
    $$\|\nabla_h g(h,\xUL)-\nabla_h g(h',\xUL)\| \leq L_g\|h-h'\| \ \forall h, h'\in\calH$$
            
    \item \label{assump:Bounded_LLGrad} The norm of gradient is bounded i.e. $\|\nabla_h g(h,\xUL)\|\leq \Omega_g\;\forall h\in\calH$.

    \item \label{assump:Lipschtz_LL2ndOrder} The derivatives $\nabla_\xUL\nabla_h g(h,\xUL)$ and $\nabla^2_h g(h,\xUL)$ are $\tau_g^h$- and $\rho_g^h$-Lipschitz w.r.t. $h$ respectively. 
    $$\|\nabla_\xUL\nabla_h g(h,\xUL)-\nabla_\xUL\nabla_h g(h',\xUL)\|\leq\tau_g^h\|h-h'\| \ \forall h, h'\in\calH$$
    $$\|\nabla^2_h g(h,\xUL)-\nabla^2_h g(h',\xUL)\|\leq\rho_g^h\|h-h'\| \ \forall h, h'\in\calH$$

    \item \label{assump:Bounded_LLHessian} The spectral norm of Hessian matrix is bounded  i.e. $\|\nabla^2_{h} g(h,\xUL)\|\leq L_g\ \forall h\in\calH$.

    \item \label{assump:Bounded_LLCrossHessian} The spectral norm of cross Hessian matrix is bounded  i.e. $\|\nabla_\xUL\nabla_{h} g(h,\xUL)\|\leq \Lambda_g^{\xUL h}\ \forall h\in\calH$.    
\end{enumerate}

For any given $h$, 
\begin{enumerate}[label=\normalfont(L\arabic*),resume]
\item \label{assump:Lipschtz_LL2ndOrder-x} The derivatives $\nabla_\xUL\nabla_h g(h,\xUL)$ and $\nabla^2_h g(h,\xUL)$ are $\tau_g^\xUL$- and $\rho_g^\xUL$-Lipschitz w.r.t. $\xUL$ respectively. 
$$\|\nabla_\xUL\nabla_h g(h,\xUL)-\nabla_\xUL\nabla_h g(h,\xUL')\|\leq\tau^\xUL_g\|\xUL-\xUL'\| \ \forall \xUL, \xUL'\in\calC$$
$$\|\nabla^2_h g(h,\xUL)-\nabla^2_h g(h,\xUL')\|\leq\rho_g^\xUL\|\xUL-\xUL'\| \ \forall \xUL, \xUL'\in\calC$$
\end{enumerate}
\end{assumption}

\subsection{Notation.}\label{sec:notations}
In this paper, the notation $[d]$ refers to the set $\{1,\hdots,d\}$ for $d>0$. For a vector $a\in\Real^d$, its square $a^2=({a_1^2},\hdots,{a_d^2})^\top$, square root $\sqrt{a}=(\sqrt{a_1},\hdots,\sqrt{a_d})^\top$ and diagonalization $\Diag(a)=\begin{bmatrix} a_1 & \hdots & 0 \\ 0 & \ddots & 0 \\ 0 & 0 & a_d\end{bmatrix}.$ For vectors $a,b\in\Real^d$, their multiplication $a\circ b=(a_1b_1,\hdots,a_db_d)^\top$ and division $a/b=(a_1/b_1,\hdots,a_d/b_d)$. For matrix $A\in\Real^{d\times d}$, $\eig(A)$ is the vector of eigenvalues of $A$. The notation $\intr(S)$ means the interior of a set $S$.

In algorithm, the upper-level iteration is indexed using $k$ and the lower-level using $t$. For a fixed upper level iterate $\xUL^k$, the lower-level iterate is denoted by $h^{k,t}$ and the optimal lower-level solution $h^*(\xUL^k)$ by $h^{k,*}$. The approximate Jacobian matrix $\partial h^{k,t}/\partial \xUL^k$ is denoted by $\sfR_{k,t}$ and the true Jacobian matrix $\partial h^{k,*}/\partial \xUL^k$ by $\sfR_{k,*}$.

\section{Main results.}\label{sec:main-results}

\subsection{Gradient estimation.} \label{sec:GradientEst} 
The ideal algorithm is gradient descent with exact gradient $\nabla F(\xUL)$ (at $k$\textsuperscript{th} upper-level iteration) given as
\begin{equation}\label{eq:ExactGrad}
    \nabla F(\xUL^k)= \nabla_\xUL f\big( h^{k,*},\xUL^k\big) + \sfR_{k,*}^\top \nabla_h f\big( h^{k,*},\xUL^k\big)
\end{equation}
In practice, it is generally difficult to find $h^{k,*}$, $\sfR_{k,*}^\top$. The solution is to iteratively solve the lower-level problem for a fixed upper-level decision and use the resulting lower-level iterates to estimate the exact gradient \eqref{eq:ExactGrad}. In this paper, we consider the following gradient estimator
\begin{equation}\label{eq:HyperGrad}
    \widehat\nabla F(\xUL^k)= \nabla_\xUL f\big( h^{k,D},\xUL^k\big) + \sfR_{k,D}^\top \nabla_h f\big( h^{k,D},\xUL^k\big)
\end{equation}
where $D$ is the fixed number of lower-level iterations for each $k$.

\subsection{Algorithm.} \label{sec:Algorithm}
In the algorithm considered in this paper, a projected gradient descent (PGD) step updates the upper-level decision and a projected mirror descent (PMD) step w.r.t. Kullback–Leibler (KL) divergence updates the lower-level (also discussed in \citet{krichene2015convergence}) since the following projection step \eqref{eq:PMDStep_Opt}  on simplex has a closed form \eqref{eq:Alg_PMDStep} in Algorithm \ref{alg:Bilevel_ITD}. 
\begin{align}\label{eq:PMDStep_Opt}
    & h^{k,t+1} := \underset{h\in\calH}{\arg\min}\ g^\eta(h^{k,t},\xUL^k) + \left\langle \nabla g^\eta(h^{k,t},\xUL^k), h-h^{k,t} \right\rangle + \frac{1}{\alpha_k} \overline{\sfD}_\KL(h,h^{k,t})
\end{align}
where KL divergence $\overline{\sfD}_\KL(h,h')$ is specified in Appendix \ref{appendix:Prelim_MD}. The Jacobian matrix $\sfR_{k,t}$ is obtained by iteratively differentiating through the PMD dynamics. The overall algorithm is outlined as below.
\begin{algorithm}
\begin{algorithmic}[1]
\State {Input}: $K$, $D$, $\eta$, step sizes: $\left(\beta_k,\ \alpha_{k}\right)_{k=0}^{K-1}$.
\State {Initialize:} $\xUL^0\in\calC$, $h^0_i\in\widetilde \calH_i\ i=1,\hdots,\sfW$. 
\For{$k=0,\hdots,K-1$}
    \State If $k>0$ set $h^{k,0} = h^{k-1,D}$, otherwise $h^{k,0}=h^0$. Set $\sfR_{k,0}={\partial h^{k,0}}/{\partial{\xUL^k}}=0$.
   \For{$t=0,\hdots,D-1$}
    \begin{align}
        & h_i^{k,t+1}=\frac{h_i^{k,t}\circ \exp\left(-\alpha_k \nabla_{h_i}g^\eta(h^{k,t},\xUL^k)\right)} {{h_i^{k,t}}^\top \exp\left(-\alpha_k \nabla_{h_i}g^\eta(h^{k,t},\xUL^k)\right)} \quad \forall i\in[\sfW] \label{eq:Alg_PMDStep} \\ \nonumber\\
        & \sfR_{k,t+1} = \Phi(\sfR_{k,t},h^{k,t},h^{k,t+1},\xUL^k) \label{eq:Alg_JacbStep}
    \end{align}    
    \EndFor
    \State Obtain estimation of gradient $\widehat\nabla F(\xUL^{k})$ using \eqref{eq:HyperGrad},
    \begin{align}
         \xUL^{k+1}=\text{Proj}_{\calC}\big(\xUL^{k}-\beta_k \widehat\nabla F(\xUL^{k})\big)
    \end{align}
\EndFor
\caption{PGD with Iterative Differentiation of PMD on simplex}
\label{alg:Bilevel_ITD}
\end{algorithmic}
\end{algorithm}
\newline 
Consider $\widetilde\calH=\widetilde\calH_1\times\cdots\times\widetilde\calH_\sfW\subset\intr(\calH)$  derived in Lemma \ref{lemma:Iters_SmplxIntr} (under Assumption \ref{assump:Bounded_LLGrad}) where for each $i=1,\hdots,\sfW$ $$\widetilde \calH_i=\big\{h_i\in\Real^{\sfq_\ell^i}: h_i\geq \nu^{\min} \cdot \mathbf{1}_{\sfq_\ell^i},\ \mathbf{1}^\top h_i = 1\big\}\subset\intr(\calH_i),$$
$0 < \nu^{\min}:=e^{-\frac{2\Omega_g}{\eta}}/\sfq_\ell^{\max}<1$. 
Note in Algorithm \ref{alg:Bilevel_ITD}, the initialization $h^{0}=(h^{0}_i,i\in[\sfW])\in \widetilde \calH$. Consequently, from Lemma \ref{lemma:Iters_SmplxIntr} it holds that $h^{k,t}\in\widetilde\calH\ \forall k,t$ and $h^{k,*}\in\widetilde\calH\ \forall k$.

\subsection{Convergence rates.}\label{sec:ConvgRate}
\begin{theorem}\label{theorem:Unconstrained}
Consider $\calC=\Real^{\sfq_u}$.  There exist constants $\overline \alpha$, $\underline{\beta}$, $\overline{\beta}$ such that if $\alpha_k=\alpha$, $\beta_k=\beta$ satisfying $0<\alpha\leq\overline \alpha$ and $\underline{\beta}\leq\beta\leq \overline \beta$, then the output of Algorithm \ref{alg:Bilevel_ITD} satisfies the following for $D\geq D^0:=D^0(\eta,\alpha)$   \begin{equation}
    \begin{aligned}
    \frac{1}{K} \sum_{k=0}^{K-1} \|\nabla F(\xUL^k)\|^2 \leq \calO\left(\frac{1}{K}\right) + \calO(\lambda^D)
    \end{aligned}
\end{equation}
where $0<\lambda:=\lambda(\eta,\alpha)<1$ is a constant.
\end{theorem}

\begin{theorem}\label{theorem:Constrained}
Consider closed convex set $\calC\subset\Real^{\sfq_u}$. There exist constants $\overline \alpha$, $\underline{\beta}$, $\overline{\beta}$ such that if $\alpha_k=\alpha$, $\beta_k=\beta$ satisfying $0<\alpha\leq\overline \alpha$ and $\underline{\beta}\leq\beta\leq \overline \beta$, 
then the output of Algorithm \ref{alg:Bilevel_ITD} satisfies the following for $D\geq D^0:=D^0(\eta,\alpha)$ 
\begin{equation}
    \begin{aligned}
    \frac{1}{K} \sum_{k=1}^{K} \mathrm{dist}^2\left(-\nabla F(\xUL^{k}),N_\calC(\xUL^{k})\right) \leq \calO\left(\frac{1}{K}\right) + \calO(\lambda^D)
    \end{aligned}
\end{equation}
where $0<\lambda:=\lambda(\eta,\alpha)<1$ is a constant, $N_\calC(\xUL)=\left\{p\in\Real^{\sfq_u} |\langle p ,\ \xUL' - \xUL \rangle \leq 0 \;\; \forall \xUL'\in \calC\right\}$ is the normal cone to the set $\calC$ at point $\xUL$ and 
$\mathrm{dist}\left(-\nabla F(\xUL),N_\calC(\xUL)\right)$ is the distance of $-\nabla F(\xUL)$ to $N_\calC(\xUL)$.
\end{theorem}

\section{Example Applications.}\label{sec:ExampleApps} \ 
\subsection{Traffic Equilibrium/Routing Game (Lower-level problem)}\label{sec:TrafficEqbm}
Consider a road network with $\mathsf{A}$ number of links and $\sfW$ number of populations where each population is characterized by a unique origin-destination (OD) pair. For each population or OD pair $i\in[\sfW]$, let $\sfq^i_\ell$ be the number of different path (or route) choices for travel 
and $\xi_i$ be the fixed travel demand. 
Define 
$\delta^{ai}(\xi_i)$
as the vector of link-paths weighted incidence between link $a$ and OD pair $i$ . Its $j$\textsuperscript{th} entry is 
\begin{align*}
    \delta^{ai}_j(\xi_i)= \begin{cases}
        \xi_i & \text{ if path } j \text{ uses link } a \\
        0 & \text{ o.w.}
    \end{cases}
\end{align*}
where $j\in[\sfq^i_\ell]$. For each link $a\in[\mathsf{A}]$, the overall vector of link-paths weighted incidence is given by 
$\delta^a(\xi) = [\delta^{a1}(\xi_1);\hdots;\delta^{a\sfW}(\xi_\sfW)]_{\sfq_\ell\times 1}$ where $\sfq_\ell=\sum_{i\in[\sfW]} \sfq^i_\ell$. Let $h_i=(h_{i,j},\ j\in [\sfq^i_\ell])_{\sfq^i_\ell\times 1}$ denote the vector of path flows for population $i$ where $h_{i,j}$ is the fraction of population $i$ traveling on path $j\in[\sfq^i_\ell]$. The overall vector of path flows is $h=(h_i,\ i\in[\sfW])_{\sfq_\ell\times 1}$. Say $\xUL$ is some fixed upper-level decision such as capacity expansion of road networks or traffic signal setting. 
For a given $\xUL$, if $c_a(\yLL_a,\xUL)$ is the travel time on link $a$ where $\yLL_a$ is the amount of flow through $a$, then the lower level objective function $g(\cdot)$ which is the potential function of the routing game is given by
\begin{equation}\label{eq:LL_Eqbm}
g(h,\xUL) = \sum_{a\in [\mathsf{A}]} \int_{0}^{x_a}c_a(u,\xUL)du 
\end{equation}
where 
$x_a={\delta^a(\xi)}^\top h=\sum_{i\in[\sfW]}\sum_{j\in[\sfq^i_\ell]} \delta^{ai}_j(\xi_i)\cdot h_{i,j}$ is the total flow through link $a$ and $h_i\in \mathcal{H}_i= \big\{h_i\in\Real^{\sfq_\ell^i}: h_i\geq0,\; \mathbf{1}^\top h_i = 1\big\}$.

Define 
$\Delta(\xi)=\big[\delta^1(\xi), \ \hdots, \ \delta^\mathsf{A}(\xi)\big]_{\sfq_\ell\times\mathsf{A}}$
as the full links-paths weighted incidence matrix. Note that $\nabla_h g(h,\xUL) = \Delta(\xi) \begin{bmatrix} c_1(\yLL_1,\xUL); & \hdots; &  c_\mathsf{A}(\yLL_\mathsf{A},\xUL) \end{bmatrix}_{\mathsf{A}\times1}$, $\nabla^2_h g(h,\xUL) = \Delta(\xi) \Diag\left(\frac{\partial c_1(\yLL_1,\xUL)}{\partial \yLL_1},\hdots,\frac{\partial c_\mathsf{A}(\yLL_\mathsf{A},\xUL)}{\partial \yLL_\mathsf{A}}\right) {\Delta(\xi)}^\top$
where $\yLL = {\Delta(\xi)}^\top h$.    
Then  $\nabla^2_h g(h,\xUL)\succeq0$ if $c_a(\yLL_a,\xUL)$ is monotonically increasing in $\yLL_a$ i.e ${\partial c_a(\yLL_a,\xUL)}/{\partial \yLL_a}\geq0\ \forall a$. This monotonicity condition implies convexity of $g$. Therefore, if $c_a(\cdot)$ is differentiable w.r.t. $\yLL_a$ and $\xUL$ then convexity of $g$ as well as assumptions \ref{assump:Lipschtz_LLGrad}-\ref{assump:Lipschtz_LL2ndOrder-x} can be checked.

\subsection{Application-1: Transportation Network Design}\label{sec:App1-NetowrkDesign} Consider a system planner who wants to make a decision on expanding the capacity of roads in a transportation network that minimizes the total time of travel plus the construction cost. But the travel time in the network actually depends on the equilibrium traffic pattern induced by capacity decisions. More formally, let $\xUL_a$ be the decision on capacity expansion of link $a$ upper bounded by $u_a$, $\widetilde c_a(\cdot)$ be the cost of travel on link $a$, $G_a(\cdot)$ be the investment cost for construction on link $a$ and $\theta$ is the conversion factor from investment cost to travel times. With $\xUL=(\xUL_a, {a\in [\mathsf{A}]})$, the objective function and constraints of the upper-level (system planner) problem are
\begin{equation}
\begin{aligned}
& F(\xUL) = f(h^*(\xUL),\xUL) = \textstyle \sum_{a\in [\mathsf{A}]} \left(\widetilde c_a(\yLL_a^*(\xUL),\xUL_a)\cdot\yLL_a^*(\xUL)+\theta G_a(\xUL_a)\right) \\
& \calC = \big\{\xUL|\ 0\leq \xUL_a\leq u_a\;\;\forall a\in [\mathsf{A}]\big\}\\
\end{aligned}
\end{equation}
where 
$\yLL^*(\xUL) = {\Delta(\xi)}^\top h^*(\xUL)$ is the vector of link flows and $h^*(\xUL)$ is the equilibrium path flow of the routing game described in Section \ref{sec:TrafficEqbm}. The assumptions \ref{assump:Lipschtz_ULGrad}, \ref{assump:Bounded_ULGrad} are checkable for well-behaved functions $\widetilde c_a(\cdot)$, $G_a(\cdot)$ such as differentiable once w.r.t. $\yLL_a$ and $\xUL$.

\subsection{Application-2: Traffic Signal Control}\label{sec:App2-TrafficControl} 
Consider the problem of optimization of traffic signal timings for area traffic control in order to minimize the total rate of delay and number of stops in a traffic network where both of the performance indices depend on the resulting traffic pattern from the setting of traffic signals. More formally, the decision needs to be made on the start and duration of green signal groups at various junctions in a road network. The detailed description of parameters and decision variables are given below.\\
\underline{\smash{Parameters}}\\
$\lambda_{\min}$: minimum green\\ 
$\tau_{jlm}$: clearance time between the green for signal group $j$ and $l$ at junction $m$ \\ 
$\Omega_m(j,l)$: 0 if the start of green for signal group $j$ proceeds that of $l$ at junction $m$, 1 otherwise\\
$D_a(\cdot)$: rate of delay on link $a$ \\
$S_a(\cdot)$: number of stops per unit time on link $a$
\\
\underline{\smash{Decision variables}}\\
$\varsigma$: reciprocal of common cycle time\\ 
$\theta=[\theta_{jm}]$: vector of start of green, where $\theta_{jm}$ is the start of green for signal group $j$ at junction $m$ \\
$\phi=[\phi_{jm}]$: vector of duration of green where $\phi_{jm}$ is the duration of green for group $j$ at junction $m$  
\\ \\
Let $\xUL:=(\varsigma,\theta,\phi)$. The constraint set and objective function (which is the weighted combination of rate of delay and number of stops with weights $W_{aD}$ and $W_{aS}$ respectively on link $a$) of the upper-level problem are given below.
\begin{equation}
\begin{aligned}
& F(\xUL)=f(h^*(\xUL),\xUL) = \sum_{a\in [\mathsf{A}]} W_{aD} D_a(\yLL^*(\xUL),\xUL) + W_{aS}S_a(\yLL^*(\xUL),\xUL)\\
& \calC = \begin{aligned}[t] \big\{(\varsigma,\theta,\phi)| \
& \varsigma_{\min}\leq \varsigma\leq \varsigma_{\max},\ \lambda_{\min}\varsigma\leq \phi_{jm}\leq 1\;\;\forall j,m,\\ 
& \theta_{jm}+\phi_{jm}+\tau_{jlm}\varsigma \leq \theta_{lm} + \Omega_m(j,l)\;\;j\neq l,\forall j,l,m\big\}
\end{aligned}
\end{aligned}
\end{equation}
where $\yLL^*(\xUL) = {\Delta(\xi)}^\top h^*(\xUL)$ is the vector of link flows and $h^*(\xUL)$ is the equilibrium path flow of the routing game described in Section \ref{sec:TrafficEqbm}.
Again the assumptions \ref{assump:Lipschtz_ULGrad}, \ref{assump:Bounded_ULGrad} are checkable for functions $D_a(\cdot)$, $S_a(\cdot)$ differentiable once w.r.t. $\yLL$ and $\xUL$.

\subsection{Numerical Experiments}\label{sec:NumericExprmnt}
\begin{figure}[htbp]
   \begin{tabular}{cc}
     \qquad \includegraphics[width=0.4\textwidth]{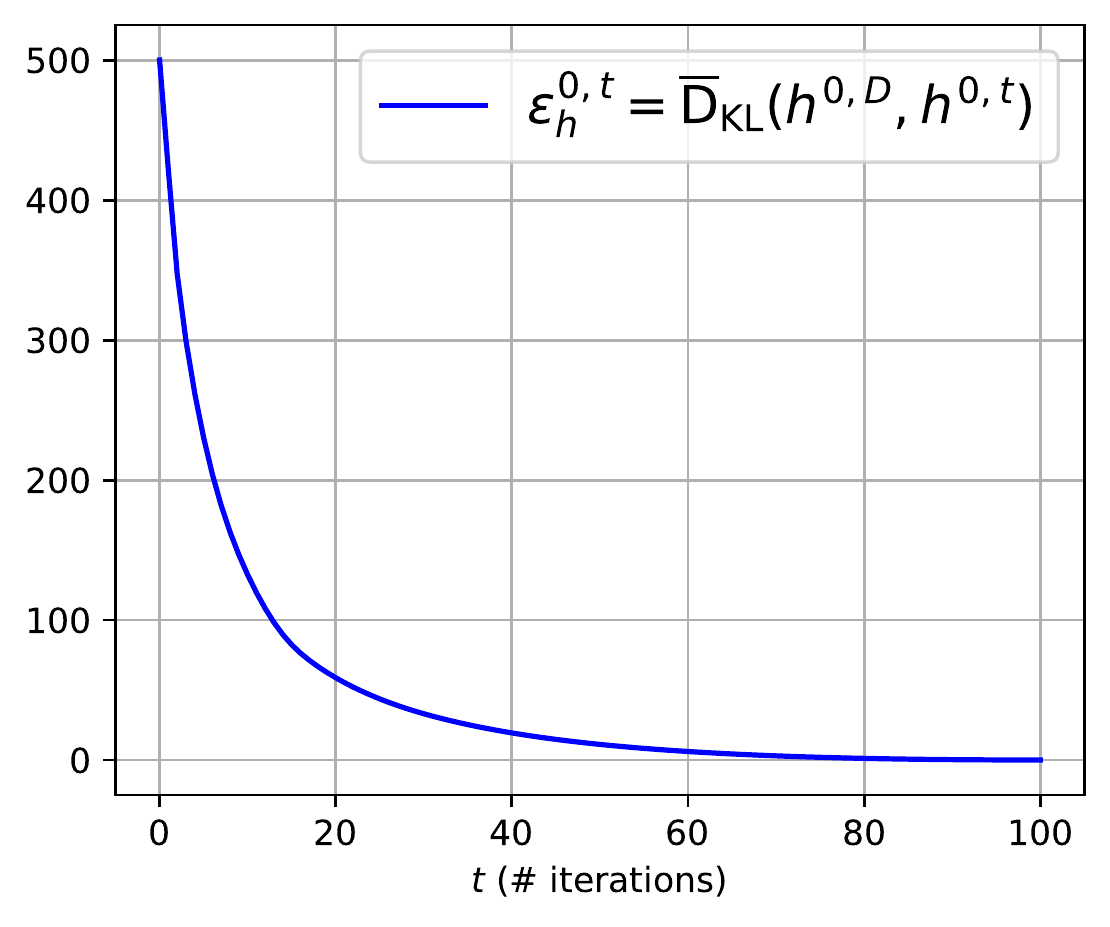} &  
     \qquad 
     \includegraphics[width=0.4\textwidth]{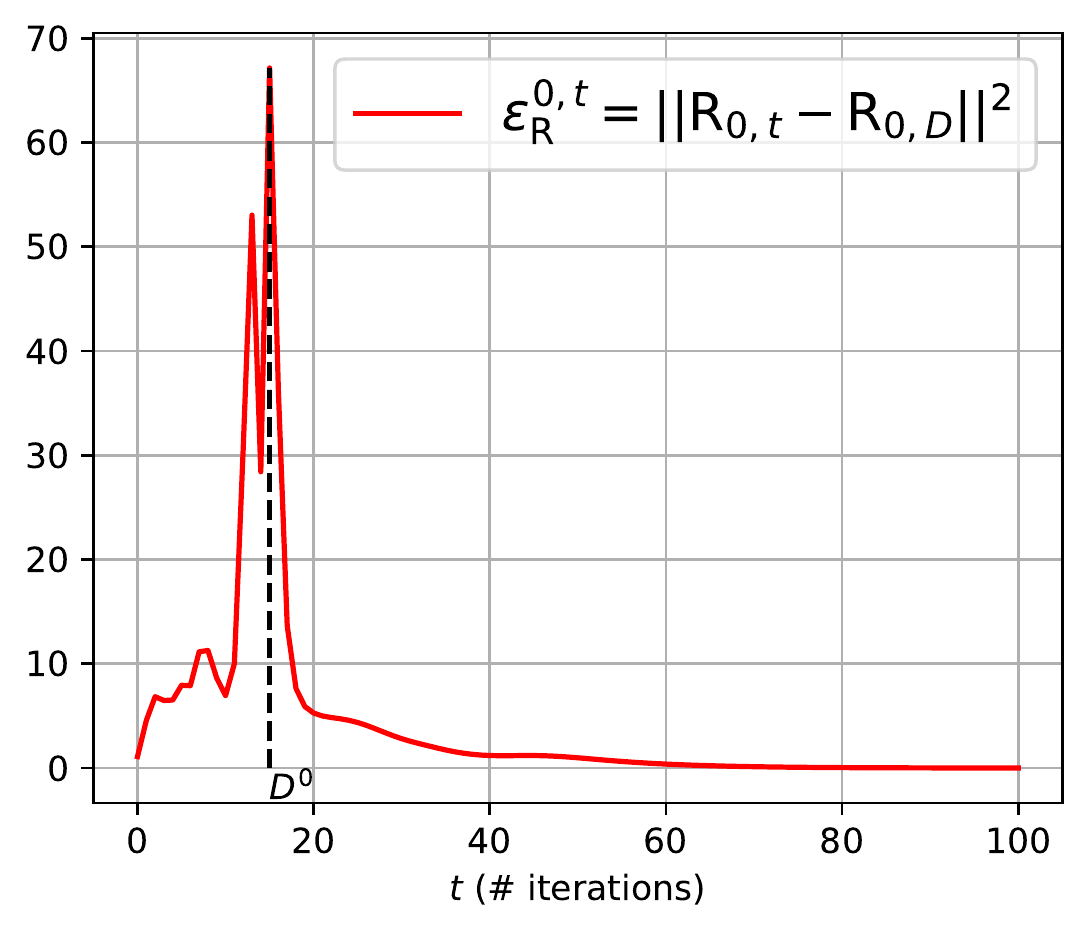}
     \\
      \qquad \scriptsize (a) & \qquad \scriptsize (b)
\end{tabular}
\caption{Convergence behavior of lower-level ($k=0$ and $D=100$).}
\label{fig:SiouxFalls-2}
\end{figure}

We test Algorithm \ref{alg:Bilevel_ITD} on Sioux Falls network which is a standard in traffic equilibrium literature. It has 24 nodes, 76 arcs, 552 OD pairs. For the computations in this section, we consider the application of transport network design (discussed in Section \ref{sec:App1-NetowrkDesign}). All the problem data has been taken from the paper by \citet{suwansirikul1987equilibrium}. The travel time function is as given below
\[
c_a(\yLL_a,\xUL) = A_a + B_a \left(\frac{\yLL_a}{K_a+\xUL_a}\right)^4
\]
where $A_a$ is the free flow travel time, $K_a$ is the existing capacity of link $a$ and $B_a>0$ is a parameter. All these values are specified in \citet{suwansirikul1987equilibrium} for each link $a$. Also, the travel demand matrix between OD pairs is given. In the upper-level objective, the travel cost function $\widetilde c_a(\yLL_a,\xUL_a)$ on link $a$ is set to be the same as $c_a(\yLL_a,\xUL)$, the investment cost function, $G_a(\xUL_a)=d_a \cdot \xUL_a^2$, and the conversion factor $\theta=0.001$. Again, the values of $d_a$ for each link $a$ are given in \citet{suwansirikul1987equilibrium}. Out of $76$ arcs, \citet{suwansirikul1987equilibrium} select 10 arcs (16, 17, 19, 20, 25, 26, 29, 39, 48, 74) for link capacity expansion upper bounded by $25$. Specifically,
\[
u_a = 
\begin{cases}
25 & \text{if } a\in\{16, 17, 19, 20, 25, 26, 29, 39, 48, 74\}\\
0 & \text{o.w.}
\end{cases}
\]
For each OD pair, we extract 5 shortest paths based on free flow travel time $A_a$ and use them along with travel demand data to construct matrix $\Delta(\xi)$ described in Section \ref{sec:TrafficEqbm}. We initialize upper-level iterate $\xUL^0=\mathbf{0}$ and lower-level iterate $h^0_i=\begin{bmatrix}\frac{1}{5} & \frac{1}{5} & \frac{1}{5} & \frac{1}{5} & \frac{1}{5} \end{bmatrix}$ for each OD pair $i$.
We run Algorithm \ref{alg:Bilevel_ITD} for different values of fixed lower-level iterations, $D\in\{40,60,80,100,120\}$. Other algorithm parameters: $\eta = 0.01$, lower-level step size $\alpha=0.50$, upper-level step size $\beta=0.25$, and no. of upper-level iterations $K=100$. All computations run on a Windows 10 machine with 16GB RAM and 1.80GHz Intel Core i7 CPU. We plot the upper-level objective $f(h^{k,D},\xUL^k)$ versus upper-level iteration $k$ in Figure \ref{fig:SiouxFalls-1} and summarize the computational time results in Table \ref{tab:time_results}. In
order to illustrate the significance of $D^0$ in Theorems \ref{theorem:Unconstrained} and \ref{theorem:Constrained}, we also plot the convergence of lower-level iterates $h^{k,t}$, $\sfR^{k,t}$ for $k=0$ in Figure \ref{fig:SiouxFalls-2}. 

\begin{table}[ht]
\hspace{0.25cm}
\begin{minipage}[h]{0.55\linewidth}
    \centering
    \includegraphics[scale=0.5]{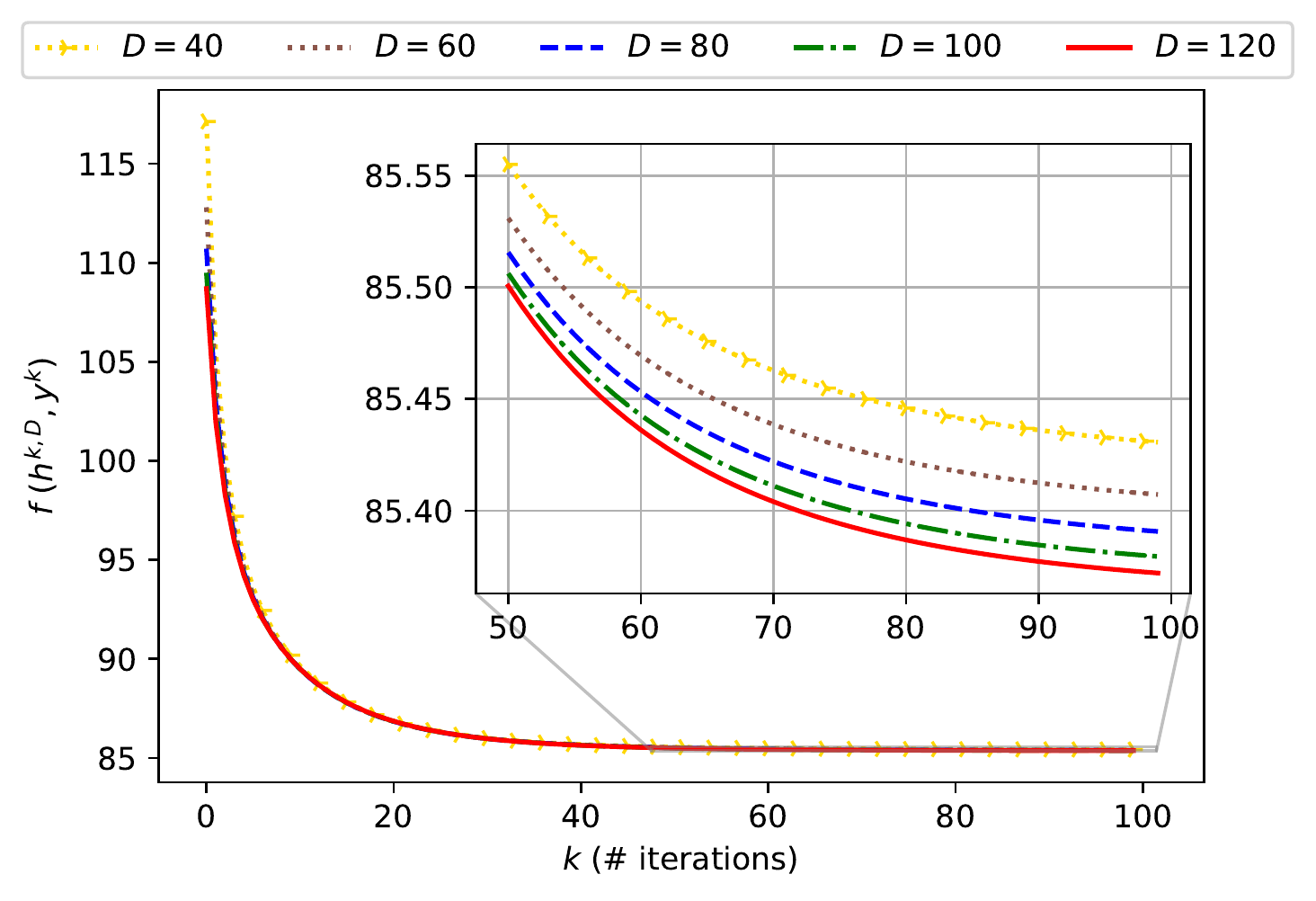}
    \captionof{figure}{$f(h^{k,D},\xUL^k)$ vs $k$.}
    \label{fig:SiouxFalls-1}
\end{minipage}
\hfill
\hspace{0.5cm}
\begin{minipage}[h]{0.35\linewidth}
    \centering
    \begin{tabular}{|c|c|}
        \hline
        \multirow{2}{*}{\quad$D$\quad} &  CPU Time per $k$ \\
                    & (in seconds) \\ \hline
        $40$ &  30.83  \\ 
        $60$ &  50.00  \\ 
        $80$ &  68.11  \\ 
        $100$ & 87.48  \\ 
        $120$ & 112.00   \\ \hline
    \end{tabular}
    \caption{Average time per upper-level iteration.}
    \label{tab:time_results}
\end{minipage}
\end{table}
Following observations can be made:
\begin{enumerate}
    \item From Figure \ref{fig:SiouxFalls-1}, observe for a fixed $k$ the upper-level objective improves with increasing $D$. This improvement comes at a cost of increasing the runtime per upper-level iteration as reported in Table \ref{tab:time_results}.
    There is an increase by roughly $20$ seconds for every 20 iterations increment in $D$.
    \item We can also observe the law of diminishing returns for successive increments in $D$ as it results in smaller improvement in the upper-level objective compared to the computational effort needed.
    \item From Figure \ref{fig:SiouxFalls-2}, observe that the error in $h^{k,t}$ converges monotonically to zero whereas error in $\sfR_{k,t}$ is not necessarily monotonic. It first increases to achieve a peak and then decays to zero. This behavior is an empirical evidence to our convergence rate results in Theorems \ref{theorem:Unconstrained} and \ref{theorem:Constrained} in the sense that a threshold $D^0$ lower-level iterations need to be performed in order to ensure convergence.   
\end{enumerate}

\section{Proofs.}\label{sec:Proofs}
In this section, we provide the proof of theorems stated in Section \ref{sec:ConvgRate}.
In Section \ref{sec:Lemma-Lowerlevel}, we state lemmas corresponding to lower-level problem which will be used in later proofs. In Section \ref{sec:Lemma-Smooth-F_x}, we provide lemmas that prove smoothness of $F(\xUL)$. In Section \ref{sec:GradEstErr}, we  derive the bound on gradient estimation error, and in Section \ref{sec:MainTheorems} we use this bound to prove Theorems \ref{theorem:Unconstrained} and \ref{theorem:Constrained}. 

\subsection{Lower-level problem}\label{sec:Lemma-Lowerlevel}
Denote the error in lower-level solution 
and error in its Jacobian matrix, respectively as
\[\txteps_h^{k,t}=\overline{\sfD}_\KL(h^{k,*},h^{k,t}),\qquad \txteps_\sfR^{k,t} = \|\sfR_{k,t}-\sfR_{k,*}\|^2.\]

The initialization error for $k>0$, $\txteps_h^{k,0}=\overline{\sfD}_\KL\big(h^*(\xUL^k),h^{k,0}\big):=\overline{\sfD}_\KL\big(h^*(\xUL^k),h^{k-1,D}\big)$ and for $k=0$, $\txteps_h^{0,0}=\overline{\sfD}_\KL\big(h^{0,*},h^{0,0}\big):=\overline{\sfD}_\KL\big(h^*(\xUL^0),h^{0}\big)$.

\subsubsection{Lipschitz continuity and boundedness results}
\begin{lemma}{(Globally bounding $\txteps_h^{k,0}$)}\label{lemma:GlobBndErr}
For each $k$, the initial error in lower-level solution is bounded i.e. $\txteps_h^{k,0}\leq \overline{\sfD}_\KL^{\max}=\sfW \ln(1/\nu^{\min}) + \nu^{\min} \sum_{i=1}^\sfW \sfq_\ell^i \ln(1/\sfq_\ell^i)$.
\end{lemma}
\proof{Proof of Lemma \ref{lemma:GlobBndErr}.}
See Appendix \ref{appendix:AuxProofs}. \Halmos
\endproof

\begin{lemma} \label{lemma:Lipschtz-ln_h}
    \begin{enumerate}[label=\normalfont(\alph*)]
        \item $\nabla_h\overline\psi_\sfH(h)$ is $1/\nu^{\min}$--Lipschitz continuous $\forall h\in\widetilde \calH$.
        \item For a fixed $\xUL$, $\nabla_h g^\eta(h,\xUL)$ is $L_g^\eta$--Lipschitz continuous $\forall h\in\widetilde \calH$ where $L_g^\eta=L_g+\eta/\nu^{\min}$.
    \end{enumerate} 
\end{lemma}
\proof{Proof of Lemma \ref{lemma:Lipschtz-ln_h}.}
See Appendix \ref{appendix:Dynamics-h^k,t}. \Halmos
\endproof

\subsubsection{Convergence rate of lower-level iterates}
\begin{lemma}\label{lemma:ConvgRate-Soln}
    Let the sequence $\{h^{k,t}\}$ be generated by Algorithm \ref{alg:Bilevel_ITD}. For a fixed $k$ and $0 < \alpha_k \leq \widetilde \alpha$,
    \begin{align} 
   \frac{1}{2} \|h^{k,*}-h^{k,t}\|^2 \leq \txteps_h^{k,t} \leq \txteps_h^{k,0}\ (1- \eta \alpha_k)^t \label{eq:Reln-3}
\end{align}
    where $\widetilde\alpha = 1/(L_g+\eta/\nu^{\min})<1/\eta$.
\end{lemma}
\proof{Proof of Lemma \ref{lemma:ConvgRate-Soln}.}
See Appendix \ref{appendix:ConvgRate-h^k,t}. \Halmos
\endproof

Define $\overline{\alpha}=\min\left\{\widetilde \alpha, {1}/{(2L_g^2/\eta+2\eta)}\right\}$.
\begin{lemma}\label{lemma:ConvgRate-Jacobian}
Let the sequence $\{\sfR_{k,t}\}$ be generated by Algorithm \ref{alg:Bilevel_ITD}. For a fixed $k$ and $0 < \alpha_k \leq \overline \alpha$, 
\begin{align}
    \txteps_\sfR^{k,t} \leq (\lambda_k)^{t}\ \big(\Gamma_{1k}+\Gamma_{2k}\ {\txteps_h^{k,0}}\big),\quad t\geq T^0_k
\end{align}
where $0<\lambda_k<1$ and constants $\lambda_k$, $\Gamma_{1k}$, $\Gamma_{2k}$, $T^0_k$ depend on $\eta$ \& $\alpha_k$.
 \end{lemma}
 \proof{Proof of Lemma \ref{lemma:ConvgRate-Jacobian}.}
See Appendix \ref{appendix:ConvgRate-R_k,t}. \Halmos
\endproof

\begin{remark}
     All results mentioned hereafter hold for $\alpha_k$ satisfying $0 < \alpha_k \leq \overline \alpha$.
\end{remark}

\subsection{Smoothness of $F(\xUL)$}\label{sec:Lemma-Smooth-F_x}
The smoothness of $F(\xUL)$ follows from key Lemmas \ref{lemma:Bnd-R*_Lip-h*} \& \ref{lemma:Lipschtz_matrixR^*}.

\begin{lemma}\label{lemma:Bnd-R*_Lip-h*}
    Denote $\sfR_*(\xUL)=\partial h^*(\xUL)/\partial x$.
    \begin{enumerate}[label=\normalfont(\alph*)]
        \item \label{lemma:Bounded_R*} $\exists$ $\sfC_0>0$ such that for any upper level decision $\xUL$, $\|\sfR_*(\xUL)\|\leq\sfC_0$. 
        \item \label{lemma:Lipschtz_h^*} Consequently, $h^*(\xUL)$ is $\sfC_0$--Lipschitz continuous w.r.t. $\xUL$.
    \end{enumerate}    
\end{lemma}

\begin{lemma}\label{lemma:Lipschtz_matrixR^*}
   The Jacobian matrix 
   $\sfR_*(\xUL)$ is $L_{\sfR_*}$--Lipschitz continuous w.r.t. $\xUL$.
\end{lemma}

\begin{lemma}\label{lemma:Lipschtz_BilevelGrad}
    $\nabla F(\xUL)$ is $L_F$--Lipschitz continuous w.r.t. $\xUL$ where $L_F=L_{\sfR_*} \Omega_f+  L_f (1+\sfC_0)^2$.
\end{lemma}

The proof of Lemmas \ref{lemma:Bnd-R*_Lip-h*}--\ref{lemma:Lipschtz_BilevelGrad} can be found in Appendix \ref{appendix:Lip_F}.

\subsection{Gradient estimation error.}\label{sec:GradEstErr}
The goal of this section is to derive Lemma \ref{lemma:GradEstErrBnd-2} which gives an expression for gradient estimation error that is easy to work with. Specifically, we bound the gradient estimation error in terms of quantities that decay exponentially w.r.t. $D$ (no. of lower-level iterations) plus the sum of deviations of negative gradient from the normal cone.

\begin{lemma}\label{lemma:GradEstErrBnd-1}
    Let the sequence $\{\xUL^k\}$ be generated by Algorithm \ref{alg:Bilevel_ITD}. For $D\geq T^0_k$, 
    \begin{align*}
    \big\|\widehat\nabla F(\xUL^k)-\nabla F(\xUL^k)\big\|^2 & \leq \sfC^D_{1k} + \sfC^D_{2k}\ \txteps_h^{k,0}
\end{align*}
where $\sfC^D_{1k}=\calO\left((\lambda_k)^D\right)$, $\sfC^D_{2k}=\calO\left((1- \eta \alpha_k)^{D}\right)+\calO((\lambda_k)^{D})$, and $\lambda_k$, $T^0_k$ 
as given in Lemma \ref{lemma:ConvgRate-Jacobian}.
\end{lemma}
\proof{Proof of Lemma \ref{lemma:GradEstErrBnd-1}.}
Recall $\nabla F(\xUL^k) = \nabla_\xUL f\big(h^{k,*},\xUL^k\big) + \sfR_{k,*}^\top \nabla_h f\big(h^{k,*},\xUL^k\big)$ and $\widehat\nabla F(\xUL^k) = \nabla_\xUL f\big(h^{k,D},\xUL^k\big) + \sfR_{k,D}^\top \nabla_h f\big(h^{k,D},\xUL^k\big)$. Therefore, we have
\begin{align*}
     \widehat\nabla F(\xUL^k)-\nabla F(\xUL^k) = \nabla_\xUL f\big(h^{k,D},\xUL^k\big) - \nabla_\xUL f\big(h^{k,*},\xUL^k\big) 
     & + \left(\sfR_{k,D}-\sfR_{k,*}\right)^\top \nabla_h f\big(h^{k,D},\xUL^k\big) \\
     & + \sfR_{k,*}^\top \left[\nabla_h f\big(h^{k,D},\xUL^k\big)-\nabla_h f\big(h^{k,*},\xUL^k\big)\right]
\end{align*}
which implies
\begin{align*}
    \big\|\widehat\nabla F(\xUL^k)-\nabla F(\xUL^k)\big\|^2 & \leq 3 L_f^2\ \big\|h^{k,D}-h^{k,*}\big\|^2 + 3 \Omega_f^2\ \|\sfR_{k,D}-\sfR_{k,*}\|^2 + 3 L_f^2 \sfC_0^2\ \big\|h^{k,D}-h^{k,*}\big\|^2
\end{align*}
Using $\|h^{k,D}-h^{k,*}\|^2\leq 2 \overline{\sfD}_\KL(h^{k,*},h^{k,D})=2\txteps_h^{k,D}$ gives
\begin{align*}
    \big\|\widehat\nabla F(\xUL^k)-\nabla F(\xUL^k)\big\|^2 
    & \leq 6 L_f^2 \left(1+\sfC_0^2\ \right) \txteps_h^{k,D} + 3 \Omega_f^2 \txteps_\sfR^{k,D} \\
    & \begin{aligned}[t]
            \leq 6 L_f^2 \left(1+\sfC_0^2\right)\ (1- \eta \alpha_k)^{D}\ \txteps_h^{k,0} + 3 \Omega_f^2\ (\lambda_k)^{D}\ \big(\Gamma_{1k}+\Gamma_{2k}\ \txteps_h^{k,0}\big) \\
             \left(\text{from Lemma \ref{lemma:ConvgRate-Soln} and Lemma \ref{lemma:ConvgRate-Jacobian}} \right)   
      \end{aligned}
\end{align*}
Taking $\sfC^D_{1k}=3 \Omega_f^2 (\lambda_k)^{D} \Gamma_{1k}$ and $\sfC^D_{2k}=6 L_f^2 \left(1+\sfC_0^2\right) (1- \eta \alpha_k)^{D}+3 \Omega_f^2 (\lambda_k)^{D} \Gamma_{2k}$ completes the proof.\Halmos
\endproof

\begin{remark}\label{remark:Const_alpha}
    For $\alpha_k=\alpha$, it can be shown that $\lambda_k=\lambda$,\; $T^0_k=T^0$,\; $\sfC^D_{1k}=\sfC^D_{1}:=\calO\left(\lambda^D\right)$,\; $\sfC^D_{2k}=\sfC^D_{2}:=\calO\left((1- \eta \alpha)^{D}\right)+\calO(\lambda^{D})$.
\end{remark}

\begin{lemma}{(Convergence of $h^{k,0}$ w.r.t. $k$)}\label{lemma:Convg_wrt-k}
Assume $\alpha_k=\alpha$ satisfying $0 < \alpha \leq \overline \alpha$ and 
$\beta_k$ satisfying $0<\beta_k\leq \frac{\nu^{\min} \gamma}{3\sqrt{2} \sfC_0}$.
Define $\gamma=1-(1-\eta\alpha)^D$ and $D^0=\{D^0\geq T^0 | \forall D\geq D^0: \gamma\geq 4/9,\ 0<\sfC^D_2\leq 1/2\}$. 
\newline For $D\geq D^0$,
\begin{align*}
    {\txteps_h^{k,0}} & \leq \left(1-\frac{\gamma}{4}\right)^k\ {\txteps_h^{0,0}} + \sfC^D_{1} \left(1-\frac{\gamma}{4}\right)^{k}
    + {\gamma}/{4} \sum_{\ell=0}^{k-1} \left(1-\frac{\gamma}{4}\right)^{k-1-\ell} \|\nabla F(\xUL^{\ell})+z^{\ell+1}\|^2
\end{align*}
where 
$z^{k+1}\in N_\calC(\xUL^{k+1})$ ($N_\calC(\xUL)$ is the normal cone to the set $\calC$ at point $\xUL$).
\end{lemma}
\proof{Proof of Lemma \ref{lemma:Convg_wrt-k}.}
See Appendix \ref{appendix:Proof-lemma:Convg_wrt-k}. \Halmos
\endproof

\begin{lemma}\label{lemma:GradEstErrBnd-2}
Assume  
$0 < \alpha_k=\alpha \leq \overline \alpha$ and $0<\beta_k\leq \frac{\nu^{\min} \gamma}{3\sqrt{2} \sfC_0}$. Let the sequence $\{\xUL^k\}$ be generated by Algorithm \ref{alg:Bilevel_ITD}.
For $D\geq D^0$,
    \begin{align}
        & \big\|\widehat\nabla F(\xUL^k)-\nabla F(\xUL^k)\big\|^2 \nonumber \\
        & \leq \sfC^D_{1} \big(1 + \sfC^D_{2}\big) + \sfC^D_{2}\ \left(1-\frac{\gamma}{4}\right)^k\ \txteps_h^{0,0} + \sfC^D_{2}\ \gamma/4\ \sum_{\ell=0}^{k-1} \left(1-\frac{\gamma}{4}\right)^{k-1-\ell} \|\nabla F(\xUL^{\ell})+z^{\ell+1}\|^2
    \end{align}
where $z^{k+1}\in N_\calC(\xUL^{k+1})$ and $\gamma$, $D^0$ as defined in Lemma \ref{lemma:Convg_wrt-k}.
\end{lemma}
\proof{Proof of Lemma \ref{lemma:GradEstErrBnd-2}.}
Substitute the bound on $\txteps_h^{k,0}$ from Lemma \ref{lemma:Convg_wrt-k} to upper bound the gradient estimation error in Lemma \ref{lemma:GradEstErrBnd-1} where $\sfC^D_{1k}=\sfC^D_1$, $\sfC^D_{2k}=\sfC^D_2$ when $\alpha_k=\alpha$ (Remark \ref{remark:Const_alpha}). The proof is  complete by using $\left(1-\frac{\gamma}{4}\right)^k\leq1$ for the term in multiplication with $\sfC^D_1\sfC^D_2$.
\Halmos
\endproof

\begin{remark}
     All results mentioned hereafter hold for $\alpha_k=\alpha$ satisfying $0 < \alpha \leq \overline \alpha$.
\end{remark}

\subsection{Proof of Theorems \ref{theorem:Unconstrained} and \ref{theorem:Constrained}}\label{sec:MainTheorems}
This section provides the proof of main theorems which follow from key Lemma \ref{lemma:Pre_ConvgRate}. In Lemma \ref{lemma:Pre_ConvgRate}, we bound the average deviation of negative gradient from normal cone
utilizing the expression for gradient estimation error from Lemma \ref{lemma:GradEstErrBnd-2}.
\begin{lemma}\label{lemma:Pre_ConvgRate}
Fix $\beta_k:=\beta$ satisfying $0<\beta\leq \overline \beta=\min\left\{\frac{\nu^{\min} \gamma}{3\sqrt{2} \sfC_0}, \frac{7}{48 L_F}\right\}$. 
For $D\geq D^0$,
\begin{align}\label{eq:Pre_ConvgRate}
    & \frac{1}{K} \sum_{k=0}^{K-1} \|\nabla F(\xUL^{k})+z^{k+1}\|^2 
    {
    \leq \frac{4 (F(\xUL^{0}) - \inf_\xUL F(\xUL)) + \left({\beta}+2L_F \beta^2\right) \sum_{k=0}^{K-1}  \left\{{3} \sfC^D_{1} + \left(1-\frac{\gamma}{4}\right)^k\txteps_h^{0,0}\right\}}{\left(\beta-6 L_F \beta^2\right) K}
    }
\end{align}
where $z^{k+1}\in N_\calC(\xUL^{k+1})$ and $\gamma$, $\sfC^D_{1}$, $D^0$ as defined in Lemma \ref{lemma:GradEstErrBnd-2}.  
\end{lemma}
\proof{Proof of Lemma \ref{lemma:Pre_ConvgRate}.}
From Lemma \ref{lemma:Lipschtz_BilevelGrad},we know $F(\cdot)$ is smooth. Therefore for any $\xUL$, $\xUL'$ following holds
\begin{align*}
    F(\xUL)\ \leq F(\xUL') + \langle \nabla F(\xUL'), \xUL-\xUL'\rangle + {L_F}/{2} \|\xUL-\xUL'\|^2 
\end{align*}
Take $\xUL:=\xUL^{k+1}$ and $\xUL':=\xUL^k$ to get
\begin{align*}
    F(\xUL^{k+1}) & \leq F(\xUL^{k}) + \left\langle \nabla F(\xUL^{k}), \xUL^{k+1}-\xUL^k\right\rangle + {L_F}/{2}\ \|\xUL^{k+1}-\xUL^k\|^2
\end{align*}
{Since $z^{k+1}\in N_{\mathcal{C}}(\xUL^{k+1})$, it implies $\langle z^{k+1}, \xUL^k-\xUL^{k+1} \rangle \leq 0$.  Add $\langle z^{k+1}, \xUL^{k+1}-\xUL^k \rangle \geq 0$ on RHS}
\begin{align*}
& F(\xUL^{k+1}) \\
& \leq F(\xUL^{k}) + \left\langle \nabla F(\xUL^{k})+z^{k+1}, \xUL^{k+1}-\xUL^k\right\rangle + {L_F}/{2}\ \|\xUL^{k+1}-\xUL^k\|^2 \\
& = F(\xUL^{k}) -\beta_k \big\langle \nabla F(\xUL^{k})+z^{k+1}, \widehat\nabla F(\xUL^{k})+z^{k+1}\big\rangle + {L_F}/{2}\ \beta_k^2 \|\widehat\nabla F(\xUL^{k})+z^{k+1}\|^2 \\
& \leq \begin{aligned}[t] F(\xUL^{k}) - \beta_k\|\nabla F(\xUL^{k})+z^{k+1}\|^2 -\beta_k \big\langle \nabla F(\xUL^{k})+z^{k+1}, \widehat\nabla F(\xUL^{k})-\nabla F(\xUL^{k})\big\rangle \\ 
+ L_F \beta_k^2 \|\widehat\nabla F(\xUL^{k})-\nabla F(\xUL^{k})\|^2 + L_F \beta_k^2 \|\nabla F(\xUL^{k})+z^{k+1}\|^2
\end{aligned}
\end{align*}
using $-a^\top b\leq \|a\|\cdot\|b\|\leq (\|a\|^2+\|b\|^2)/2$
\begin{align*}
& F(\xUL^{k+1}) \\ 
& \leq F(\xUL^{k}) - \left(\frac{\beta_k}{2}-L_F \beta_k^2\right) \|\nabla F(\xUL^{k})+z^{k+1}\|^2 + \left(\frac{\beta_k}{2}+L_F \beta_k^2\right) \|\widehat\nabla F(\xUL^{k})-\nabla F(\xUL^{k})\|^2\\
& \leq \begin{aligned}[t] & F(\xUL^{k}) - \left(\frac{\beta_k}{2}-L_F \beta_k^2\right) \|\nabla F(\xUL^{k})+z^{k+1}\|^2 \\ 
&  \begin{aligned}[t] + \left(\frac{\beta_k}{2}+L_F \beta_k^2\right) \Big\{\sfC^D_{1} \big(1 + \sfC^D_{2}\big) + \sfC^D_{2} \left(1-\frac{\gamma}{4}\right)^k \txteps_h^{0,0} + \sfC^D_{2} \gamma/4 \sum_{\ell=0}^{k-1} \left(1-\frac{\gamma}{4}\right)^{k-1-\ell} \|\nabla F(\xUL^{\ell})+z^{\ell+1}\|^2\Big\} \\
\left(\text{using Lemma \ref{lemma:GradEstErrBnd-2}}\right) \end{aligned}
\end{aligned}
\end{align*}
Summing from $k=0$ to $K-1$
\begin{align*}
& \sum_{k=0}^{K-1}\left(\frac{\beta_k}{2}-L_F \beta_k^2\right) \|\nabla F(\xUL^{k})+z^{k+1}\|^2 \\ 
& \begin{aligned}[t] \leq F(\xUL^{0}) - F(\xUL^K) + \sum_{k=0}^{K-1}\left(\frac{\beta_k}{2}+L_F \beta_k^2\right) & \Big\{\sfC^D_{1} (1 + \sfC^D_{2}) + \sfC^D_{2} \left(1-\frac{\gamma}{4}\right)^k \txteps_h^{0,0} \\ & \qquad + \sfC^D_{2} \frac{\gamma}{4}\sum_{\ell=0}^{k-1} \left(1-\frac{\gamma}{4}\right)^{k-1-\ell} \|\nabla F(\xUL^{\ell})+z^{\ell+1}\|^2 \Big\}
\end{aligned}\\
&  \begin{aligned}[t] \leq F(\xUL^{0}) - \inf_\xUL F(\xUL) + \sum_{k=0}^{K-1}\left(\frac{\beta_k}{2}+L_F \beta_k^2\right) \Big\{\sfC^D_{1} \big(1 + \sfC^D_{2}\big) + \sfC^D_{2} \left(1-\frac{\gamma}{4}\right)^k \txteps_h^{0,0} + \sfC^D_{2} \|\nabla F(\xUL^{k})+z^{k+1}\|^2 \Big\} \\
\begin{aligned}[t] \Biggl(\text{using } & \textstyle \sum_{k=1}^{K-1}\left(\frac{\beta_k}{2}+L_F \beta_k^2\right) \sum_{\ell=0}^{k-1} \left(1-\frac{\gamma}{4}\right)^{k-1-\ell} \|\nabla F(\xUL^{\ell})+z^{\ell+1}\|^2 \\
& \leq \textstyle \frac{1}{\gamma/4} \sum_{k=0}^{K-1} \left(\frac{\beta_k}{2}+L_F \beta_k^2\right) \|\nabla F(\xUL^{k})+z^{k+1}\|^2\quad \Biggr)
\end{aligned}
\end{aligned}
\end{align*}
Rearranging gives
\begin{align*}
& \sum_{k=0}^{K-1} \left((1-\sfC^D_{2})\frac{\beta_k}{2}-(1+\sfC^D_{2}) L_F \beta_k^2\right) \|\nabla F(\xUL^{k})+z^{k+1}\|^2 \\
& \leq F(\xUL^{0}) - \inf_\xUL F(\xUL) + \sum_{k=0}^{K-1} \left(\frac{\beta_k}{2}+L_F \beta_k^2\right)  \left\{\sfC^D_{1} \big(1 + \sfC^D_{2}\big) + \sfC^D_{2}\left(1-\frac{\gamma}{4}\right)^k\txteps_h^{0,0}\right\}
\end{align*}
Using $0<C^D_2\leq 1/2$ for $D\geq D^0$ from the definition of $D^0$, 
\begin{align*}
& \sum_{k=0}^{K-1} \left(\frac{\beta_k}{4}-\frac{3}{2} L_F \beta_k^2\right) \|\nabla F(\xUL^{k})+z^{k+1}\|^2 \\
& \leq F(\xUL^{0}) - \inf_\xUL F(\xUL) + \sum_{k=0}^{K-1} \left(\frac{\beta_k}{2}+L_F \beta_k^2\right)  \left\{\frac{3}{2} \sfC^D_{1} + \frac{1}{2}\left(1-\frac{\gamma}{4}\right)^k\txteps_h^{0,0}\right\}
\end{align*}
First multiply with $4$ both sides, then fix $\beta_k=\beta$ where $0<\beta<\frac{1}{6L_F}$ and divide both sides by the resulting constant term in LHS to get \eqref{eq:Pre_ConvgRate}. Take $\beta\leq\frac{7}{48 L_F}$ and combine this with the condition on $\beta$ from Lemma \ref{lemma:GradEstErrBnd-2} to complete the proof. \Halmos
\endproof

\proof{Proof. of Theorem \ref{theorem:Unconstrained}.}
For $\calC=\Real^{\sfq_u}$, the normal cone $N_\calC(\xUL)=\{0\}\ \forall \xUL\in\Real^{\sfq_u}$. Therefore, from Lemma \ref{lemma:Pre_ConvgRate} we have for $D\geq D^0$ and $0<\beta\leq \overline \beta$, 
\begin{align}\label{eq:Unconstr-Eq1}
    \frac{1}{K} \sum_{k=0}^{K-1} \|\nabla F(\xUL^{k})\|^2 
    & \leq \frac{4 (F(\xUL^{0}) - \inf_\xUL F(\xUL)) + \left({\beta}+2L_F \beta^2\right) \sum_{k=0}^{K-1}  \left\{{3} \sfC^D_{1} + \left(1-\frac{\gamma}{4}\right)^k\txteps_h^{0,0}\right\}}{\left(\beta-6 L_F \beta^2\right) K}
\end{align}

For $\beta\in\left(0,\frac{1}{6L_F}\right)$, the expression $\beta-6 L_F \beta^2$ is a concave function of $\beta$. Note that, $\exists$ $0<b<4$ such that $\left[\frac{4-b}{48L_F},\overline\beta\right]\subseteq\left[\frac{4-b}{48L_F},\frac{4+b}{48L_F}\right]$ and for $\beta\in\left[\frac{4-b}{48L_F},\frac{4+b}{48L_F}\right]$, we have
\begin{align}\label{eq:InterimBnds_beta}
    \beta-6 L_F \beta^2 \geq \frac{16-b^2}{384 L_F},\quad {\beta}+2L_F \beta^2 \leq \frac{(4+b)(28+b)}{1152 L_F}
\end{align}

Substitute above bounds and the bound on $\txteps_h^{0,0}$ (from Lemma \ref{lemma:GlobBndErr}) in \eqref{eq:Unconstr-Eq1},
\begin{align*}
    & \frac{1}{K} \sum_{k=0}^{K-1} \|\nabla F(\xUL^{k})\|^2 \\
    & \leq \frac{1536 L_F (F(\xUL^{0}) - \inf_\xUL F(\xUL)) + {(4+b)(28+b)} \overline{\sfD}_\KL^{\max} \frac{1-(1-\gamma/4)^K}{3\gamma/4}}{(16-b^2)\cdot K} + \frac{(4+b)(28+b)}{16-b^2} \sfC^D_{1} 
\end{align*}
The proof is complete by taking $\underline{\beta}=\frac{4-b}{48L_F}$ and noting $\sfC^D_{1}=\calO(\lambda^D)$. \Halmos
\endproof

\proof{Proof of Theorem \ref{theorem:Constrained}.}
\begin{align*}
    & \|\nabla F(\xUL^{k+1})+z^{k+1}\|^2 \\
    & =  \|\nabla F(\xUL^{k+1})-\nabla F(\xUL^{k})+\nabla F(\xUL^{k})+z^{k+1}\|^2 \\
    & \leq 2 \|\nabla F(\xUL^{k})+z^{k+1}\|^2 + 2 \|\nabla F(\xUL^{k+1})-\nabla F(\xUL^{k})\|^2 \\
    & \leq  2 \|\nabla F(\xUL^{k})+z^{k+1}\|^2 + 2\ L_F^2 \|\xUL^{k+1}-\xUL^k\|^2 \\
    & =  2 \|\nabla F(\xUL^{k})+z^{k+1}\|^2 + 2\ \beta^2 L_F^2 \|\widehat\nabla F(\xUL^{k})+z^{k+1}\|^2 \\
    & \leq  (4\ \beta^2 L_F^2 + 2) \|\nabla F(\xUL^{k})+z^{k+1}\|^2 + 4\ \beta^2 L_F^2 \|\widehat\nabla F(\xUL^{k})-\nabla F(\xUL^{k})\|^2 \\
    & \leq \begin{aligned}[t]
        & (4\ \beta^2 L_F^2 + 2) \|\nabla F(\xUL^{k})+z^{k+1}\|^2 \\ 
        & \begin{aligned}[t] + 4\ \beta^2 L_F^2 \Big\{\sfC^D_{1} \big(1 + \sfC^D_{2}\big) + \sfC^D_{2}\ \left(1-\frac{\gamma}{4}\right)^k\ \txteps_h^{0,0} + \sfC^D_{2}\ \gamma/4\ \sum_{\ell=0}^{k-1} \left(1-\frac{\gamma}{4}\right)^{k-1-\ell} \|\nabla F(\xUL^{\ell})+z^{\ell+1}\|^2 \Big\}\\
    \left(\text{using Lemma \ref{lemma:GradEstErrBnd-2}}\right) \end{aligned}
    \end{aligned}
\end{align*}
Summing from $k=0$ to $K-1$
\begin{align*}
     & \sum_{k=0}^{K-1} \|\nabla F(\xUL^{k+1})+z^{k+1}\|^2 \\
     & \begin{aligned}[t]\leq \left(4 \beta^2 L_F^2 (1+\sfC^D_{2}) + 2\right) \sum_{k=0}^{K-1}\|\nabla F(\xUL^{k})+z^{k+1}\|^2 + 4\ \beta^2 L_F^2 \sum_{k=0}^{K-1} \Big\{\sfC^D_{1} \big(1 + \sfC^D_{2}\big) + \sfC^D_{2}\ \left(1-\frac{\gamma}{4}\right)^k \txteps_h^{0,0} \Big\} \\
    \Big(\text{using } \textstyle \sum_{k=1}^{K-1} \sum_{\ell=0}^{k-1} \left(1-\frac{\gamma}{4}\right)^{k-1-\ell} \|\nabla F(\xUL^{\ell})+z^{\ell+1}\|^2 \leq \frac{1}{\gamma/4} \sum_{k=0}^{K-1} \|\nabla F(\xUL^{k})+z^{k+1}\|^2 \Big) 
    \end{aligned} \\
    & \begin{aligned}[t]\leq \left(6 \beta^2 L_F^2  + 2\right) \sum_{k=0}^{K-1}\|\nabla F(\xUL^{k})+z^{k+1}\|^2 + \beta^2 L_F^2 \sum_{k=0}^{K-1} \Big\{6\sfC^D_{1}+ 2 \left(1-\frac{\gamma}{4}\right)^k \txteps_h^{0,0} \Big\} \\
    \Big(\text{using $0<C^D_2\leq 1/2$ for $D\geq D^0$ from the definition of $D^0$}\Big) 
    \end{aligned}
\end{align*}
Dividing by $K$ both sides and using Lemma \ref{lemma:Pre_ConvgRate} to bound $\frac{1}{K}\sum_{k=0}^{K-1}\|\nabla F(\xUL^{k})+z^{k+1}\|^2$  gives
\begin{align*}
& \frac{1}{K} \sum_{k=0}^{K-1} \|\nabla F(\xUL^{k+1})+z^{k+1}\|^2 \\
& \leq \frac{(6 \beta^2 L_F^2  + 2)}{(\beta-6 L_F \beta^2)} \frac{4 (F(\xUL^{0}) - \inf_\xUL F(\xUL)) + \left({\beta}+2L_F \beta^2\right) \sum_{k=0}^{K-1}  \left\{{3} \sfC^D_{1} + \left(1-\frac{\gamma}{4}\right)^k\txteps_h^{0,0}\right\}}{ K} \\
& \quad + \beta^2 L_F^2 \frac{\sum_{k=0}^{K-1} \Big\{6\sfC^D_{1} + 2 \left(1-\frac{\gamma}{4}\right)^k \txteps_h^{0,0} \Big\}}{K}
\end{align*}
Similar to proof of Theorem \ref{theorem:Unconstrained},
$\exists$ $0<b<4$ such that $\left[\frac{4-b}{48L_F},\overline\beta\right]\subseteq\left[\frac{4-b}{48L_F},\frac{4+b}{48L_F}\right]$ and for $\beta\in\left[\frac{4-b}{48L_F},\frac{4+b}{48L_F}\right]$, we have
\begin{align*}
    \beta-6 L_F \beta^2 \geq \frac{16-b^2}{384 L_F},\quad {\beta}+2L_F \beta^2 \leq \frac{(4+b)(28+b)}{1152 L_F},\quad 6 \beta^2 L_F^2  + 2\leq \frac{(4+b)^2}{384} + 2
\end{align*}
Substitute above bounds and the bound on $\txteps_h^{0,0}$ (from Lemma \ref{lemma:GlobBndErr}) in \eqref{eq:Unconstr-Eq1},
\begin{align*}
& \frac{1}{K} \sum_{k=0}^{K-1} \|\nabla F(\xUL^{k+1})+z^{k+1}\|^2 \\
& \leq \frac{{(4+b)^2} + 768}{16-b^2} \frac{4 L_F (F(\xUL^{0}) - \inf_\xUL F(\xUL)) + \frac{(4+b)(28+b)}{1152} \sum_{k=0}^{K-1}  \left\{{3} \sfC^D_{1} + \left(1-\frac{\gamma}{4}\right)^k\overline{\sfD}_\KL^{\max}\right\}}{ K} \\
& \quad + \frac{(4+b)^2}{48^2} \frac{\sum_{k=0}^{K-1} \Big\{6\sfC^D_{1} + 2 \left(1-\frac{\gamma}{4}\right)^k \overline{\sfD}_\KL^{\max} \Big\}}{K}
\end{align*}

Simplifying and rearranging gives
\begin{align*}
& \frac{1}{K} \sum_{k=0}^{K-1} \|\nabla F(\xUL^{k+1})+z^{k+1}\|^2 \\
& \leq \frac{4 \widetilde b L_F (F(\xUL^{0}) - \inf_\xUL F(\xUL))}{K} 
+ \frac{(4+b)\big((4+b)+\widetilde b (28+b) \big)}{1152}  \left\{\overline{\sfD}_\KL^{\max} \frac{1-(1-\gamma/4)^K}{\gamma/4\cdot K} + 3 \sfC^D_{1}\right\}
\end{align*}

where $\widetilde b=\frac{{(4+b)^2} + 768}{16-b^2}$. The proof is complete by taking $\underline{\beta}=\frac{4-b}{48L_F}$ and noting $\sfC^D_{1}=\calO(\lambda^D)$. \Halmos
\endproof



\bibliographystyle{informs2014} 
\bibliography{Traffic-Equilibria_Applications.bib,Gradient-based_Bilevel.bib} 


\section*{Appendices.}

\begin{APPENDICES}

The proof details of various intermediate results used in proving the main results can be found in Appendices 
\ref{appendix:Proof-lemma:Convg_wrt-k}--\ref{appendix:AuxProofs}. Specifically, Appendix \ref{appendix:Proof-lemma:Convg_wrt-k} provdies the proof of Lemma \ref{lemma:Convg_wrt-k}. Appendix \ref{appendix:Dynamics-lowerlevel} provides the dynamics of lower-level iterates and Appendix \ref{appendix:ConvgRate-LowerLevel} derives their convergence rate to a fixed point. 
Appendix \ref{appendix:ProofsAppendixC} provides proof of various lemmas in Appendix \ref{appendix:ConvgRate-LowerLevel}. In Appendix \ref{appendix:Lip_F}, we prove the smoothness of $F(\xUL)$ and related lemmas. Finally, Appendix \ref{appendix:AuxProofs} proves various boundedness and Lipschitz continuity results.

\section{Proof of Lemma \ref{lemma:Convg_wrt-k}.}\label{appendix:Proof-lemma:Convg_wrt-k}
\begin{align*}
{\txteps_h^{k+1,0}} - {\txteps_h^{k,D}} & = \overline{\sfD}_\KL\big(h^{k+1,*},h^{k+1,0}\big) - \overline{\sfD}_\KL\big(h^{k,*},h^{k,D}\big) \\ & = \overline{\sfD}_\KL\big(h^{k+1,*},h^{k,D}\big) - \overline{\sfD}_\KL\big(h^{k,*},h^{k,D}\big) \quad (\text{using $h^{k+1,0}=h^{k,D}$ from Algorithm \ref{alg:Bilevel_ITD}})\\
& = \psi_\mathsf{H}\big(h^{k+1,*}\big) -\psi_\mathsf{H}\big(h^{k,*}\big) + \big\langle \nabla \psi_\mathsf{H}\big(h^{k,D}\big),\ h^{k,*}-h^{k+1,*} \big\rangle  \\
& = - \overline{\sfD}_\KL\big(h^{k,*},h^{k+1,*}\big) + \big\langle \nabla \psi_\mathsf{H}\big(h^{k,D}\big)-{\nabla \psi_\mathsf{H}\big(h^{k+1,*}\big)},\ h^{k,*}-h^{k+1,*} \big\rangle \\
& \leq - \frac{1}{2} \|h^{k+1,*}-h^{k,*}\|^2 + \|\nabla \psi_\mathsf{H}\big(h^{k+1,*}\big)-\nabla \psi_\mathsf{H}\big(h^{k,D}\big)\|\cdot \|h^{k+1,*}-h^{k,*}\| \\
& \leq - \frac{1}{2} \|h^{k+1,*}-h^{k,*}\|^2 + 1/\nu^{\min}\ \|h^{k+1,*}-h^{k,D}\|\cdot \|h^{k+1,*}-h^{k,*}\| \quad (\text{using Lemma \ref{lemma:Lipschtz-ln_h}})
\end{align*}
{Consider some $\varphi_{k+1}>0$ and using $2\|a\|\cdot \|b\|\leq \|a\|^2+\|b\|^2$,}
\begin{align*}
& {\txteps_h^{k+1,0}} - {\txteps_h^{k,D}} \\
& \leq - \frac{1}{2} \|h^{k+1,*}-h^{k,*}\|^2 + \frac{1}{2(1+\varphi_{k+1})}\|h^{k+1,*}-h^{k,D}\|^2+\frac{(1+\varphi_{k+1})/(\nu^{\min})^2}{2}\|h^{k+1,*}-h^{k,*}\|^2 \\
& \leq \frac{1}{1+\varphi_{k+1}} {\txteps_h^{k+1,0}} + \frac{(1+\varphi_{k+1})/(\nu^{\min})^2-1}{2}\sfC_0^2\ \|\xUL^{k+1}-\xUL^k\|^2 \quad (\text{using Lemma \ref{lemma:Bnd-R*_Lip-h*}\ref{lemma:Lipschtz_h^*}})
\end{align*}
Rearranging the above inequality gives
\begin{align}\label{eq:Rearrng-ineq}
    {\txteps_h^{k+1,0}} & \leq \left(1+\frac{1}{\varphi_{k+1}}\right) \left({\txteps_h^{k,D}} + \frac{(1+\varphi_{k+1})^2/(\nu^{\min})^2-1}{2} \sfC_0^2\ \|\xUL^{k+1}-\xUL^k\|^2\right)
\end{align}
Recall the following upper level update 
\begin{align*}
    & \xUL^{k+1} = \underset{\xUL\in\mathcal{C}}{\arg\min} \frac{1}{2} \|\xUL - \xUL^k + \beta_k\ \widehat\nabla F(\xUL^k)\|^2 
\end{align*}
For a closed convex set $\calC$, the optimality condition of the projection operator is
\begin{align}
    & -\big(\xUL^{k+1} - \xUL^k + \beta_k\ \widehat\nabla F(\xUL^k)\big) \in N_\mathcal{C}(\xUL^{k+1})\nonumber \\
    \implies & \exists \text{ unique } z^{k+1}\in N_\mathcal{C}(\xUL^{k+1}) \text{ such that } -\xUL^{k+1} + \xUL^k - \beta_k\ \widehat\nabla F(\xUL^k) = \beta_k z^{k+1}  \label{eq:ProjOpt_Cond}
\end{align}
where $N_C(\xUL)=\left\{p\in\Real^{\sfq_u} |\ \langle p ,\ \xUL' - \xUL \rangle \leq 0\; \forall \xUL'\in \mathcal{C}\right\}$ is the normal cone to the set $\mathcal{C}$ at point $\xUL$. As a result, \eqref{eq:ProjOpt_Cond} implies $\|\xUL^{k+1}-\xUL^k\|=\beta_{k}\|\widehat\nabla F(\xUL^{k})+z^{k+1}\|$ where $z^{k+1}\in N_C(\xUL^{k+1})$. Substituting this in \eqref{eq:Rearrng-ineq} gives
\begin{align*}
{\txteps_h^{k+1,0}} & \leq \left(1+\frac{1}{\varphi_{k+1}}\right) \left({\txteps_h^{k,D}} + \frac{(1+\varphi_{k+1})/(\nu^{\min})^2}{2} \sfC_0^2\ \beta_{k}^2\|\widehat\nabla F(\xUL^{k})+z^{k+1}\|^2\right) \\
 & \leq \left(1+\frac{1}{\varphi_{k+1}}\right)\ \left({\txteps_h^{k,D}} + \big[\sfC^D_{1k} + \sfC^D_{2k}\ {\txteps_h^{k,0}} + \|\nabla F(\xUL^{k})+z^{k+1}\|^2\big]\ \frac{\sfC_0^2 \ (1+\varphi_{k+1})\ \beta_{k}^2}{(\nu^{\min})^2} \right) \text{ for } D\geq T^0_k \\
&\qquad \big(\text{using $\|\widehat\nabla F(\xUL^{k})+z^{k+1}\|^2 \leq 2\|\widehat\nabla F(\xUL^{k})-\nabla F(\xUL^{k})\|^2 + 2\|\nabla F(\xUL^{k})+z^{k+1}\|^2$ and Lemma \ref{lemma:GradEstErrBnd-1}\big)}
\end{align*}
{Introduce $\gamma_{k}$ as $0<\gamma_{k}=1-\left(1-\eta \alpha_{k}\right)^{D}<1$ and 
using it with Lemma \ref{lemma:ConvgRate-Soln} to bound ${\txteps_h^{k,D}}$ gives}
\begin{align*}
 {\txteps_h^{k+1,0}} \leq \left(1+\frac{1}{\varphi_{k+1}}\right) & \left(\Big[1-\gamma_{k}+ \sfC^D_{2k} \frac{\sfC_0^2 (1+\varphi_{k+1}) \beta_{k}^2}{(\nu^{\min})^2}\Big] {\txteps_h^{k,0}} \right.\\
&\qquad \left. + \big[\sfC^D_{1k} + \|\nabla F(\xUL^{k})+z^{k+1}\|^2\big] \frac{\sfC_0^2 (1+\varphi_{k+1}) \beta_{k}^2}{(\nu^{\min})^2} \right)
\end{align*}
Take $\varphi_{k+1} = \frac{1-\gamma_k/2}{\gamma_k/4}$ and say following is satisfied
\begin{align}\label{eq:Req_gamma}
    \sfC^D_{2k}\ \frac{\sfC_0^2\ (1+\varphi_{k+1})\ \beta_{k}^2}{(\nu^{\min})^2}\leq \frac{\gamma_k}{2}
\end{align} 
which leads to
\begin{align}
    {\txteps_h^{k+1,0}} & \leq \left(1-\frac{\gamma_k}{4}\right)\ {\txteps_h^{k,0}} + \frac{\sfC_0^2 \ \beta_{k}^2}{(\nu^{\min})^2}\ \frac{(1-\gamma_k/4)^2}{(1-\gamma_k/2)\ \gamma_k/4}\ \big[\sfC^D_{1k} + \|\nabla F(\xUL^{k})+z^{k+1}\|^2\big] \nonumber\\
    & \leq \left(1-\frac{\gamma_k}{4}\right)\ {\txteps_h^{k,0}} + \frac{9 \sfC_0^2}{2(\nu^{\min})^2}\ {\beta_{k}^2}/{\gamma_k}\ \big[\sfC^D_{1k} + \|\nabla F(\xUL^{k})+z^{k+1}\|^2\big] \label{eq:err_h^k+1,0}\\
    & \qquad \big(\text{using that $1<{(1-\gamma_k/4)^2}/{(1-\gamma_k/2)}<{9}/{8}$ for $0<\gamma_k<1$}\big)\nonumber
\end{align}
Doing recursion on \eqref{eq:err_h^k+1,0} w.r.t. $k$ gives
\begin{align}
    {\txteps_h^{k,0}} & \leq \prod_{\ell=0}^{k-1} \left(1-\frac{\gamma_\ell}{4}\right)\ {\txteps_h^{0,0}} + \frac{9 \sfC_0^2}{2(\nu^{\min})^2}\ \sum_{\ell=0}^{k-1} {\beta_{\ell}^2}/{\gamma_\ell}\ \big[\sfC^D_{1\ell} + \|\nabla F(\xUL^{\ell})+z^{\ell+1}\|^2\big] \prod_{j=\ell+1}^{k-1} \left(1-\frac{\gamma_j}{4}\right) \label{eq:Recur-err_h^k,0}
\end{align}
\newline Consider a fixed lower-level step size i.e. $\alpha_k=\alpha$ where $0<\alpha\leq\overline \alpha$. Then $\gamma_k=\gamma$, $T^0_k=T^0$, 
$\sfC^D_{1k}=\sfC^D_1$, and $\sfC^D_{2k}=\sfC^D_2$. Additionally if upper-level step size $\beta_k$ satisfies $\beta_k\leq \frac{\nu^{\min} \gamma}{3\sqrt{2} \sfC_0}$, \eqref{eq:Recur-err_h^k,0} finally becomes
\begin{align}
    {\txteps_h^{k,0}} & \leq \left(1-\frac{\gamma}{4}\right)^k\ {\txteps_h^{0,0}} + \sfC^D_{1} \left(1-\frac{\gamma}{4}\right)^{k}
    + {\gamma}/{4} \sum_{\ell=0}^{k-1} \big[\|\nabla F(\xUL^{\ell})+z^{\ell+1}\|^2\big] \left(1-\frac{\gamma}{4}\right)^{k-1-\ell} \label{eq:Recur-err_h^k,0_constSS}
\end{align}
and the requirement \eqref{eq:Req_gamma} becomes -- there should exist $D\geq T^0$ such that $\gamma$ satisfies
\begin{subequations}
\begin{align}\label{eq:Req_gamma-a}
\frac{1-\gamma/4}{\gamma/4} \leq \frac{2\gamma}{\sfC^D_{3}} \iff \frac{\gamma^2}{2} + \frac{\sfC^D_{3} \gamma}{4} - \sfC^D_{3} \geq 0 \iff \gamma \geq -\sfC^D_{3}/4+\sqrt{(\sfC^D_{3})^2/16+2\sfC^D_{3}}
\end{align}
where $\sfC^D_{3} = {4 \sfC_0^2\sfC^D_{2}\beta_{k}^{2}}/{(\nu^{\min})^2}\leq {2\sfC^D_{2}\gamma^2}/{9}$. 
Note the following necessary condition needs to be satisfied be satisfied by $\sfC^D_{3}>0$
\begin{align}\label{eq:Req_gamma-b}
-\sfC^D_{3}/4+\sqrt{(\sfC^D_{3})^2/16+2\sfC^D_{3}} < 1 \implies \sfC^D_{3} < 2/3
\end{align}
\end{subequations}
Now, since $\sfC^D_{2}$ decays exponentially fast to 0 w.r.t. $D$ then $\exists$ $\widehat D$ such that $0<\sfC^D_{2}\leq 1/2$ for $D\geq \widehat D$. This implies \eqref{eq:Req_gamma-b} is satisfied for $D\geq \widehat D$ as $\sfC^D_{3}\leq \gamma^2/9<1/9$ since $0<\gamma<1$. Observe that the RHS in \eqref{eq:Req_gamma-a} is increasing in $\sfC^D_3$ so if $D\geq \widehat D$ it suffices for $\gamma$ to satisfy
$\gamma\geq -(1/9)/4+\sqrt{(1/9)^2/16+2\cdot 1/9}=4/9$ which by definition of $\gamma$ means choosing $D$ such that $\left(1-\eta \alpha\right)^D\leq 5/9\iff D\geq \widetilde D = \ln(5/9)/\ln(1-\eta\alpha)$. 
\\ \\
Define $D^0=\max\{T^0,\widehat D,\widetilde D\}$. Then the bound \eqref{eq:Recur-err_h^k,0_constSS} holds for $D\geq D^0$. This completes the proof.\Halmos
\endproof

\section{Preliminaries of mirror descent.}\label{appendix:Prelim_MD} In this appendix, we state the concepts related to mirror descent and some properties of Bregman divergence.
\begin{enumerate}
\item For $\psi(\cdot)$ which is $\mu_\psi$-strongly convex, 
    \begin{enumerate}
    \item $\overline{\psi}(h)=\sum_{i=1}^\sfW \psi(h_i)$ is $\mu_\psi$-strongly convex.     
    \item The Bregman divergence for $i\in[\sfW]$ w.r.t. $\psi$ is given by \[\sfD_\psi(h_i,h'_i)=\psi(h_i)-\psi(h'_i)-\langle \nabla \psi(h'_i), h_i-h'_i\rangle \;\; \text{ where } \;\; h_i, h'_i\in\dom(\psi)\]
    \item Define the overall Bregman divergence as $\overline{\sfD}_\psi(h,h')=\sum_{i=1}^\sfW\sfD_\psi(h_i,h'_i)$ then
    \[\overline{\sfD}_{\psi}(h,h') = \overline{\psi}(h) - \overline{\psi}(h') - \langle \nabla_h \overline{\psi}(h') , h-h'\rangle \;\; \text{ where } \;\; h, h'\in\dom(\overline{\psi})\]
    i.e. $\overline{\sfD}_{\psi}$ is the Bregman divergence w.r.t. $\overline{\psi}$.
     \item The PMD update w.r.t. $\overline{\sfD}_{\psi}$ for a function $q(\cdot)$ is given by
    \begin{align*}
    & \mathsf{z}^{t+1} = (\nabla \overline{\psi})^{-1}\left[\nabla \overline{\psi}(h^{t})-\alpha\nabla_h q(h^{t}) \right]\\
    & h^{t+1} = \underset{h\in H}{\arg\min}\ \overline{\sfD}_\psi\left(h,\ \mathsf{z}^{t+1}\right)
    \end{align*}
    where $H\subseteq\dom(\overline{\psi})\bigcap \dom(q)$.
    \item \label{propty:StrongCnvx-1} For a convex function $q(\cdot)$, the function $\widetilde q(h) = q(h) + \eta\ \overline{\psi}(h)$ is $\eta$-strongly convex w.r.t. $\overline{\sfD}_{\psi}$ i.e.
    \[\widetilde q(h) \geq \widetilde q(h') + \langle \nabla_h \widetilde q(h'), h-h'\rangle + \eta \overline{\sfD}_{\psi}(h,h')\quad \forall h,h'\in\dom(\overline{\psi})\bigcap \dom(q) \]
    Since $\overline{\sfD}_{\psi}(h,h')\geq \mu_\psi/2\ \big\|h-h'\big\|_2^2$ then $\widetilde q$ is $\eta\mu_\psi$-strongly convex w.r.t. Euclidean norm $\|\cdot\|_2$.   
    \end{enumerate}
    
\item For the case of negative Shanon entropy i.e. $\psi(h_i):=\psi_\sfH(h_i)=\sum_{j=1}^{\sfq_\ell^i} \left(h_{i,j} \ln(h_{i,j})-h_{i,j}\right)$ 
    \begin{enumerate}
    \item The resulting Bregman divergence is known as KL divergence given by
    \[\sfD_\KL(h_i,h'_i)=\sum_{j=1}^{\sfq_\ell^i} h_{i,j} \ln\frac{h_{i,j}}{h'_{i,j}} - \left(h_{i,j}-h'_{i,j}\right) \quad \text{and} \quad \overline{\sfD}_\KL(h,h')=\sum_{i=1}^\sfW \sfD_\KL(h_i,h'_i)\]
    \item \label{propty:StrongCnvx} $\overline{\psi}_\sfH(h) = \sum_{i=1}^\sfW \psi_\sfH(h_i)$ is 1-strongly convex $\forall h\in\calH$. Using \ref{propty:StrongCnvx-1} above we have for any given $\xUL$, the regularized lower-level objective $g^\eta(\cdot,\xUL)$ is $\eta$-strongly convex w.r.t. both KL divergence $\overline{\sfD}_\KL$ and Euclidean norm $\|\cdot\|_2$ for all $h\in\calH$.
    \end{enumerate}

\end{enumerate}

\section{Dynamics of lower-level iterates}\label{appendix:Dynamics-lowerlevel} 
\label{appendix:LowerLevel}

\subsection{Dynamics of $h^{k,t}$.}\label{appendix:Dynamics-h^k,t}
 In this appendix, we provide the dynamics of $h^{k,t}$ and its fixed point in Lemma \ref{lemma:Dynamics-h^k,t}--\ref{lemma:FixedPoint-h^k,t}. Further, we explicitly characterize the interior $\widetilde\calH$ of simplex to which iterates $h^{k,t}$ belong in Lemma \ref{lemma:Iters_SmplxIntr} and use it to establish the Lipschitz continuity of $g^\eta(\cdot,\xUL)$ over $\widetilde\calH$ in Lemma \ref{lemma:Lipschtz-ln_h}.

\begin{lemma}\label{lemma:Dynamics-h^k,t}
For a fixed upper-level iterate $\xUL^k$, the PMD step \eqref{eq:PMDStep_Opt} w.r.t. $\overline{\sfD}_\KL$ on the lower-level constraint set $\calH$ is given by
\begin{align}\label{eq:Proj_MD}
    h_i^{k,t+1}=\frac{h_i^{k,t}\circ \exp\left(-\alpha_k \nabla_{h_i}g^\eta(h^{k,t},\xUL^k)\right)} {{h_i^{k,t}}^\top \exp\left(-\alpha_k \nabla_{h_i}g^\eta(h^{k,t},\xUL^k)\right)} \quad \forall i\in[\sfW]
\end{align}
\end{lemma}

\begin{lemma}\label{lemma:FixedPoint-h^k,t}
\begin{enumerate}[label=\normalfont(\alph*)]
    \item \label{lemma:FixedPoint-h^k,t_a} For a fixed upper-level iterate $\xUL^k$, the optimal solution $h^{k,*}:=h^*(\xUL^k)$ to lower-level problem \eqref{eq:LowerLevel-EntropyReg} is a fixed point of the dynamics \eqref{eq:Proj_MD}, i.e.
 \begin{align}\label{eq:FixedPoint}
    h_i^{k,*}=\frac{h_i^{k,*}\circ \exp\left(-\alpha_k \nabla_{h_i}g^\eta(h^{k,*},\xUL^k)\right)} {{h_i^{k,*}}^\top \exp\left(-\alpha_k \nabla_{h_i}g^\eta(h^{k,*},\xUL^k)\right)} \quad \forall i\in[\sfW]
\end{align}
\item \label{lemma:FixedPoint-h^k,t_b} For a fixed upper-level iterate $\xUL^k$, if the fixed point $h^{k,*}:=h^*(\xUL^k)$ of dynamics \eqref{eq:Proj_MD} satisfies $h^{k,*}\in\intr(\calH)$, then it is the optimal solution to lower-level problem \eqref{eq:LowerLevel-EntropyReg}. 
\end{enumerate}
\end{lemma}
\proof{Proof of Lemma \ref{lemma:FixedPoint-h^k,t}.}
See Appendix \ref{appendix:Proof-lemma:FixedPoint-h^k,t}.
\Halmos \endproof


\begin{lemma}{($h^{k,t}$ in the interior of $\calH$)} \label{lemma:Iters_SmplxIntr}
    Suppose Assumption \ref{assump:Bounded_LLGrad} holds. Let the sequence $\{h^{k,t}\}$ be generated by Algorithm \ref{alg:Bilevel_ITD}. 
    Then 
    following holds
    \begin{enumerate}[label=\normalfont(\alph*)]
        \item If $h^{0,0}:=h^{0}\in \widetilde \calH$ then $h^{k,t}\in \widetilde \calH\;\; \forall k,t$.
        \item The optimal $h^{k,*}:=h^*(\xUL^k)\in \widetilde \calH\;\; \forall k$.
    \end{enumerate}    
    where $\widetilde\calH=\widetilde\calH_1\times\cdots\times\widetilde\calH_\sfW\subset\intr(\calH)$, $\widetilde \calH_i=\big\{h_i\in\Real^{\sfq_\ell^i}: h_i\geq \nu^{\min} \cdot \mathbf{1}_{\sfq_\ell^i},\ \mathbf{1}^\top h_i = 1\big\}\subset\intr(\calH_i)$, $0 < \nu^{\min}:=e^{-\frac{2\Omega_g}{\eta}}/\sfq_\ell^{\max}<1$ and $\sfq_\ell^{\max}=\max_{i\in[\sfW]} \sfq_\ell^{\max}$.
\end{lemma}
\proof{Proof of Lemma \ref{lemma:Iters_SmplxIntr}.}
See Appendix \ref{appendix:Proof-lemma:Iters_SmplxIntr}
\Halmos \endproof

\begin{remark}\label{remark:FP-indp_alpha}
    For a fixed upper-level iterate $\xUL^k$, the optimal solution $h^{k,*}$ to lower-level problem \eqref{eq:LowerLevel-EntropyReg} is unique and is independent of step size $\alpha_k$. Since, $h^{k,*}\in\intr(\calH)$ from Lemma \ref{lemma:Iters_SmplxIntr} then using Lemma \ref{lemma:FixedPoint-h^k,t} it holds that -- $\exists$ a unique fixed point of \eqref{eq:Proj_MD} in the interior of $\calH$ which is optimal to \eqref{eq:LowerLevel-EntropyReg} and \textbf{does not} depend on $\alpha_k$. 
\end{remark}

\proof{Proof of Lemma \ref{lemma:Lipschtz-ln_h}.}
\begin{enumerate}[label=\normalfont(\alph*)]
\item Since $\nabla_h\psi_\sfH(h)=\ln(h)$ which implies $\nabla_h^2\psi_\sfH(h) = \Diag(1/h)$ and as a result 
\begin{align*}
\|\nabla_h^2\psi_\sfH(h)\| = 1/\min_{i,j} h_{i,j} \leq 1/\nu^{\min}\;\; \forall h\in\widetilde \calH 
\end{align*}
For any $h, h' \in \widetilde \calH$,
\begin{align*}
& \nabla_h \psi_\sfH(h)-\nabla_h \psi_\sfH(h') = \ln(h)-\ln(h') = \int_{t=0}^1 \nabla_h^2\psi_\sfH(h' + t(h-h')) \cdot (h-h') dt \\
\implies & \|\nabla_h \psi_\sfH(h)-\nabla_h \psi_\sfH(h')\| \leq 1/\nu^{\min} \|h-h'\|. 
\end{align*}

\item Since $\nabla_h g^\eta(h,\xUL) = \nabla_h g(h,\xUL) + \eta\nabla_h \psi_\sfH(h)$ then for $h, h' \in \widetilde \calH$, 
    \begin{align*}
        \|\nabla_h g^\eta(h,\xUL)-\nabla_h g^\eta(h',\xUL)\| \leq (L_g+\eta/\nu^{\min})\|h-h'\|
    \end{align*}
\end{enumerate} 
\Halmos
\endproof

\subsection{Dynamics of $\sfR_{k,t}$.}\label{appendix:Dynamics-R_k,t} In this appendix, we derive the Jacobian dynamics as a linear time-varying dynamic system in Lemma \ref{lemma:Dynamics-R_k,t} and argue existence of its  fixed point in Lemma \ref{lemma:FixedPoint-R_k,t}. The uniqueness of fixed point (Lemma \ref{lemma:unique_R*}) follows from intermediate Lemma \ref{lemma:M_k,*} the proof of which relies on Lemma \ref{lemma:B_k,t factor}.
\begin{lemma}\label{lemma:Dynamics-R_k,t}
For each upper-level iteration $k$, the dynamics \eqref{eq:Alg_JacbStep} of $\sfR_{k,t}:=\partial h^{k,t}/\partial \xUL^k$  in Algorithm \ref{alg:Bilevel_ITD} is given by
\begin{align} 
     \sfR_{k,t+1} & = \Phi\big(\sfR_{k,t},h^{k,t},h^{k,t+1},\xUL^k\big) = M_{k,t} \sfR_{k,t} + U_{k,t} \label{eq:Jacb_Dyn_v2}
\end{align}
where $M_{k,t} = B_{k,t}\big[(1-\eta\alpha_k)\Diag\left(1/h^{k,t}\right)-\alpha_k \nabla^2_{h} g(h^{k,t},\xUL^k)\big]$, 
$U_{k,t} = -\alpha_k B_{k,t} \nabla_\xUL \nabla_{h} g(h^{k,t},\xUL^k)$, $B_{k,t}=\blk\left(B^1_{k,t},\hdots,B^\sfW_{k,t}\right)$, $B^i_{k,t} = \Diag\big(h^{k,t+1}_i\big) - h^{k,t+1}_i {h^{k,t+1}_i}^\top$, and $\sfR_{k,0}=0$. 
\end{lemma}
\proof{Proof of Lemma \ref{lemma:Dynamics-R_k,t}.}
See Appendix \ref{appendix:Proof-lemma:Dynamics-R_k,t}.
\Halmos \endproof

\begin{lemma}\label{lemma:FixedPoint-R_k,t}
For each $k$, $\exists$ a fixed point $\sfR_{k,*}:=\partial h^{k,*}/\partial \xUL^k$ of $\Phi(\cdot)$ in Lemma \ref{lemma:Dynamics-R_k,t} as, 
\begin{align} 
     \sfR_{k,*} & = \Phi\big(\sfR_{k,*},h^{k,*},h^{k,*},\xUL^k\big) = M_{k,*} \sfR_{k,*} + U_{k,*} \label{eq:Jacb_Dyn-FixedPoint_v2}
\end{align}
where $M_{k,*} = B_{k,*}\big[(1-\eta\alpha_k)\Diag\left(1/h^{k,*}\right)-\alpha_k \nabla^2_{h} g(h^{k,*},\xUL^k)\big]$, $U_{k,*} = -\alpha_k B_{k,*} \nabla_\xUL \nabla_{h} g(h^{k,*},\xUL^k)$, $B_{k,*}=\blk\big(B_{k,*}^1,\hdots,B_{k,*}^\sfW\big)$ and $B^i_{k,*} = \Diag\big(h^{k,*}_i\big) - h^{k,*}_i {h^{k,*}_i}^\top$.
\end{lemma}
\proof{Proof of Lemma \ref{lemma:FixedPoint-R_k,t}.}
See Appendix \ref{appendix:Proof-lemma:FixedPoint-R_k,t}.
\Halmos \endproof

\begin{lemma}\label{lemma:B_k,t factor} Define $B^i = \Diag\big(h_i\big) - h_i {h_i}^\top$. Then for each $i\in[\sfW]$, $\|B^i\|\leq 1$ and $\exists$ matrix $\sfV^i\in\Real^{\sfq_\ell^i\times(\sfq_\ell^i-1)}$ such that $B^i = {\Lambda^i}{\Lambda^i}^\top$ where 
${\Lambda^i} = \Diag(\sqrt{h_i}) {\sfV^i}\in\Real^{\sfq_\ell^i \times (\sfq_\ell^i-1)}$ 
and $\begin{bmatrix} {\sqrt{h_i}} & {\sfV^i} \end{bmatrix}_{\sfq_\ell^i\times\sfq_\ell^i}$ is an orthonormal matrix.
\newline Define $B=\blk\left(B^1,\hdots,B^\sfW\right)$,  $\sfV=\blk\big({\sfV^1},\hdots,{\sfV^\sfW}\big)\in\Real^{_{\sfq_\ell\times(\sfq_\ell-\sfW)}}$, $\sqrt{H}=\blk\big(\sqrt{h_1},\hdots,\sqrt{h_\sfW})$\\ $\in\Real^{_{\sfq_\ell\times\sfW}}$. As a consequence,  $\|B\|\leq1$ and $B = \Lambda \Lambda^\top$ where
$\Lambda= \blk\big({\Lambda^1},\hdots,{\Lambda^\sfW}\big)=\Diag(\sqrt{h}) \sfV$ and $\begin{bmatrix}
\sqrt{H} & \sfV \end{bmatrix}$ is an orthonormal matrix.
\end{lemma}
\proof{Proof of Lemma \ref{lemma:B_k,t factor}.}
See Appendix \ref{appendix:Proof-lemma:B_k,t factor}.
\Halmos \endproof

\begin{remark}\label{remark:JacbFP+B_k,*}\ 
    \begin{enumerate}[label=\normalfont(\alph*)]
    \item \label{remark:JacbFP-indp_alpha}
    From Remark \ref{remark:FP-indp_alpha}, it follows naturally that $\sfR_{k,*}:={\partial h^{k,*}}/{\partial \xUL^k}$ also does not depend on $\alpha_k$. 
   
    \item \label{remark:B_k,* factor}
    As a consequence of Lemma \ref{lemma:B_k,t factor}, 
    $\exists$ matrix $\Lambda^i_{k,*}=\Diag(\sqrt{h^{k,*}_i}) {\sfV^i_{k,*}}\in\Real^{\sfq_\ell^i\times(\sfq_\ell^i-1)}$ corresponding to $B^i_{k,*}$ such that $B^i_{k,*} = {\Lambda^i_{k,*}}^\top\Lambda^i_{k,*}$ where $\begin{bmatrix} {\sqrt{h^{k,*}_i}} & {\sfV^i_{k,*}} \end{bmatrix}$ is an orthonormal matrix. 

    Define matrices $\sfV_{k,*}=\blk\big({\sfV^1_{k,*}},\hdots,{\sfV^\sfW_{k,*}}\big)$, $\sqrt{H_{k,*}}=\blk\big(\sqrt{h^{k,*}_1},\hdots,\sqrt{h^{k,*}_\sfW})$. Then
    $B_{k,*} = \Lambda_{k,*}\Lambda_{k,*}^\top$  where $\Lambda_{k,*}= \blk\big(\Lambda^1_{k,*},\hdots,\Lambda^\sfW_{k,*}\big)=\Diag(\sqrt{h^{k,*}}) \sfV_{k,*}$ and $\begin{bmatrix}
    \sqrt{H_{k,*}} & \sfV_{k,*} \end{bmatrix}$ is an orthonormal matrix. Therefore, $M_{k,*}$ can also be written as
    \begin{align}\label{eq:M_k,*-v2}
        M_{k,*} = \Lambda_{k,*}\Lambda_{k,*}^\top\big[(1-\eta\alpha_k)\Diag\left(1/h^{k,*}\right)-\alpha_k \nabla^2_{h} g(h^{k,*},\xUL^k)\big]
    \end{align}
    \end{enumerate}
\end{remark}

Define matrices $\widehat M_{k,*} = (1-\eta\alpha_k) I - \alpha_k \Lambda_{k,*}^\top \nabla^2_{h} g(h^{k,*},\xUL^k) \Lambda_{k,*}$ and $H_{k,*}=\blk\big(h^{k,*}_1,\hdots,h^{k,*}_\sfW)$.
\begin{lemma}\label{lemma:M_k,*}
    The following results hold
    \begin{enumerate}[label=\normalfont(\alph*)]
    \item \label{lemma:eigM*_result-1} $\nu^{\min}\cdot I \preceq \Lambda_{k,*}^\top \Lambda_{k,*} \preceq I,\quad \Lambda_{k,*}^\top \Diag\left(1/h^{k,*}\right) \Lambda_{k,*}=I,\quad I\preceq \Lambda_{k,*}^\top \Diag^2\left(1/h^{k,*}\right) \Lambda_{k,*} \preceq 1/\nu^{\min} I,\newline H_{k,*}^\top H_{k,*}\preceq I,\;\; \Lambda_{k,*}^\top \Diag\left(1/h^{k,*}\right) H_{k,*}=0$.
    
    \item \label{lemma:eigM*_result-2} 
    $0\preceq (1-(L_g+\eta)\alpha_k) I \preceq \widehat M_{k,*} \preceq (1-\eta\alpha_k) I$ for $0<\alpha_k \leq {1}/(L_g+\eta)$.
    
    \item \label{lemma:eigM*_result-3} 
    $0 \leq 1-(L_g+\eta)\alpha_k\leq \eig(M_{k,*}) \leq 1-\eta\alpha_k < 1$ for $0<\alpha_k \leq {1}/(L_g+\eta)$.     
    \end{enumerate}
\end{lemma}
\proof{Proof of Lemma \ref{lemma:M_k,*}.}
See Appendix \ref{appendix:Proof-lemma:M_k,*}.
\Halmos \endproof

\begin{lemma} \label{lemma:unique_R*}
 For each $k$, the Jacobian matrix $\sfR_{k,*}$ given by \eqref{eq:Jacb_Dyn-FixedPoint_v2} is unique.
\end{lemma}
\proof{Proof of Lemma \ref{lemma:unique_R*}.}
Firstly note that for a fixed $\alpha_k$, the matrices $M_{k,*}$ and $U_{k,*}$ are unique as they depend only on lower-level optimal solution $h^{k,*}$ which is unique.
Secondly, $I-M_{k,*}$ is invertible using Lemma \ref{lemma:M_k,*}\ref{lemma:eigM*_result-3} which then implies
\begin{align}\label{eq:Jacb_Dyn-FixedPoint_v3}
   \sfR_{k,*} = (I - M_{k,*})^{-1} U_{k,*} 
\end{align}
Therefore, $\sfR_{k,*}$ is uniquely determined. \Halmos
\endproof

\section{Convergence rate of lower-level iterates.}\label{appendix:ConvgRate-LowerLevel} The goal of this appendix is to derive key lemmas -- Lemma \ref{lemma:ConvgRate-Soln}, Lemma \ref{lemma:ConvgRate-Jacobian} which are related to the lower-level problem and required in the proof of Lemmas \ref{lemma:Convg_wrt-k} and \ref{lemma:GradEstErrBnd-1}. 

\subsection{Convergence rate of $h^{k,t}$.}\label{appendix:ConvgRate-h^k,t} This section proves the convergence rate of $\txteps_h^{k,t}$ w.r.t. $t$ for a fixed $k$ using standard analysis of mirror descent for strongly convex functions. The error in lower-level solution, $\txteps_h^{k,t}=\overline{\sfD}_\KL(h^{k,*},h^{k,t})$ is well-defined as $h^{k,*}$ is unique for each $k$.
\proof{Proof of Lemma \ref{lemma:ConvgRate-Soln}.}
Firstly, the optimality condition of the projection problem \eqref{eq:PMDStep_Opt} is
\begin{align}\label{eq:PMDStep_OptCond}
    \big\langle \nabla \overline \psi_\sfH(h^{k,t+1})-\nabla \overline \psi_\sfH(h^{k,t})+\alpha_k \nabla_h g^\eta(h^{k,t},\xUL^k),\ h-h^{k,t+1}\big\rangle \geq 0\;\; \forall h\in\mathcal{H}
\end{align}
Using property of Bregman divergence: $\sfD_\psi(u,w)-\sfD_\psi(u,v)-\sfD_\psi(v,w)=\langle \nabla \psi(v)-\nabla \psi(w),\ u-v\rangle$
\begin{align}\label{eq:ConvgRate-h_ineq1}
& \overline{\sfD}_\KL(h^{k,*},h^{k,t+1}) \nonumber\\
& = \overline{\sfD}_\KL(h^{k,*},h^{k,t}) - \overline{\sfD}_\KL(h^{k,t+1},h^{k,t}) - \big\langle \nabla \overline \psi_\sfH(h^{k,t+1})-\nabla \overline\psi_\sfH(h^{k,t}),\ h^{k,*}-h^{k,t+1}\big\rangle \nonumber\\
& \leq \overline{\sfD}_\KL(h^{k,*},h^{k,t})+\alpha_k \big\langle \nabla_h g^\eta(h^{k,t},\xUL^k),\ h^{k,*}-h^{k,t+1}\big\rangle-\overline{\sfD}_\KL(h^{k,t+1},h^{k,t}) \quad \text{(using \eqref{eq:PMDStep_OptCond})} \nonumber\\
& = \begin{aligned} \overline{\sfD}_\KL(h^{k,*},h^{k,t})+\alpha_k \big\langle \nabla_h g^\eta(h^{k,t},\xUL^k),\ h^{k,*}-h^{k,t}\big\rangle+\alpha_k \big\langle \nabla_h g^\eta(h^{k,t},\xUL^k),\ h^{k,t}-h^{k,t+1}\big\rangle \\ -\overline{\sfD}_\KL(h^{k,t+1},h^{k,t}) \end{aligned}\nonumber \\
& \begin{aligned}
\ \leq (1- \eta \alpha_k )\ \overline{\sfD}_\KL(h^{k,*},h^{k,t})+\alpha_k \left[g^\eta(h^{k,*},\xUL^k) - g^\eta(h^{k,t},\xUL^k)\right]+\alpha_k \big\langle \nabla_h g^\eta(h^{k,t},\xUL^k),\ h^{k,t}-h^{k,t+1}\big\rangle  \\ -\overline{\sfD}_\KL(h^{k,t+1},h^{k,t}) \\ \text{(using strong convextiy of $g^\eta(\cdot,\xUL)$ w.r.t. KL divergence from \ref{propty:StrongCnvx} in Appendix \ref{appendix:Prelim_MD})}
\end{aligned}
\end{align}
From Lemma \ref{lemma:Lipschtz-ln_h}, $\nabla_h g^\eta(h)$ is $L_g^\eta$--Lipschitz continuous for $h\in\widetilde\calH$ where $L_g^\eta$. This implies following holds 
\begin{align}\label{eq:DescentLemma-g}
    g^\eta(h,\xUL^k) \leq g^\eta(h',\xUL^k) + \langle\nabla g^\eta(h',\xUL^k),h-h'\rangle + \frac{L_g^\eta}{2} \|h-h'\|^2\quad \forall h,h'\in\widetilde\calH.
\end{align}
Since iterates $h^{k,t}\in\widetilde\calH$ (from Lemma \ref{lemma:Iters_SmplxIntr}), using \eqref{eq:DescentLemma-g} with the fact that $\overline \psi_\sfH(\cdot)$ is $1-$strongly convex on $\calH$ \big(i.e. $\frac{1}{2}\|h-h'\|^2\leq \overline{\sfD}_\KL(h,h')$\big) implies
\begin{align}\label{eq:ConvgRate-h_ineq2}
    \big\langle \nabla_h g^\eta(h^{k,t},\xUL^k),\ h^{k,t}-h^{k,t+1}\big\rangle \leq g^\eta(h^{k,t},\xUL^k)-g^\eta(h^{k,t+1},\xUL^k) + L_g^\eta\  \overline{\sfD}_\KL(h^{k,t+1},h^{k,t})
\end{align}
Substituting \eqref{eq:ConvgRate-h_ineq2} into \eqref{eq:ConvgRate-h_ineq1} gives,
\begin{align}
& \overline{\sfD}_\psi(h^{k,*},h^{k,t+1}) \nonumber \\
& \leq (1- \eta \alpha_k )\ \overline{\sfD}_\KL(h^{k,*},h^{k,t})+\alpha_k \big(g^\eta(h^{k,*},\xUL^k) - g^\eta(h^{k,t+1},\xUL^k)\big) +\left({\alpha_k  L_g^\eta}-1\right)\ \overline{\sfD}_\KL(h^{k,t+1},h^{k,t}) \label{eq:Reln-1}
\end{align}
For a fixed $k$, $g^\eta(h^{k,*},\xUL^k) \leq g^\eta(h^{k,t+1},\xUL^k)$ as $h^{k,*}$ is the optimal solution to lower-level probem \eqref{eq:LowerLevel-EntropyReg} and considering $\alpha_k$ such that $0 < \alpha_k \leq \widetilde \alpha=1/L_g^\eta$, the inequality \eqref{eq:Reln-1} implies
\begin{align}
    \overline{\sfD}_\KL(h^{k,*},h^{k,t+1}) \leq (1- \eta \alpha_k)\ \overline{\sfD}_\KL(h^{k,*},h^{k,t}) \label{eq:Reln-2}
\end{align}
where $L_g^\eta=L_g+\eta/\nu^{\min}>\eta$. Doing recursion on \eqref{eq:Reln-2} w.r.t. $t$ completes the proof. \Halmos
\endproof

\subsection{Stability analysis and convergence of $\sfR_{k,t}$.}\label{appendix:StabConvg-R_k,t}
In this appendix, we analyze the LMI stability of Jacobian dynamics in Appendix \ref{appendix:Stab-R_k,t} using concepts from robust control which allows to derive novel results on convergence rate of $\sfR_{k,t}$ in Appendix \ref{appendix:ConvgRate-R_k,t}.

The Jacobian error matrix, $\rmX_{k,t} = \sfR_{k,t} - \sfR_{k,*}$ is well defined using Lemma \ref{lemma:unique_R*}. Then using \eqref{eq:Jacb_Dyn_v2} and \eqref{eq:Jacb_Dyn-FixedPoint_v2}
\begin{align}\label{eq:JacbErr_Dyn}
    \rmX_{k,t+1} = M_{k,*} \rmX_{k,t} + \Delta_{k,t} \rmX_{k,t} + \rmW_{k,t}
\end{align}
where $\Delta_{k,t}=M_{k,t}-M_{k,*}$, $\rmW_{k,t} = U_{k,t}-U_{k,*}+ \Delta_{k,t}\ \sfR_{k,*}$, $\rmX_{k,0}=-\sfR_{k,*}$ and $\rmX_{k,*}=0$.

\subsubsection{LMI for robust stability.}\label{appendix:Stab-R_k,t}
Introducing additional variables $\Psi_{k,t}$, $\rmY_{k,t}$ to consider the following extension of dynamic system \eqref{eq:JacbErr_Dyn} for some $\epsilon_k>0$.
\begin{align}\label{eq:dynamic_sys-sep}
    \begin{bmatrix} \rmX_{k,t+1} \\ \Psi_{k,t} \\ \rmY_{k,t} \end{bmatrix} = 
    \begin{bmatrix}
    M_{k,*} & \epsilon_k I & I \\ I & 0 & 0 \\ I & 0 & 0
    \end{bmatrix}
    \begin{bmatrix}
    \rmX_{k,t} \\ \Omega_{k,t} \\ \rmW_{k,t} 
    \end{bmatrix} \quad
    \Omega_{k,t} = \frac{1}{\epsilon_k } \Delta_{k,t} \Psi_{k,t}
\end{align}
This allows to consider analysis techniques used for robust stability and performance analysis of discrete-time uncertain system \citet{skelton2013unified}.
\begin{lemma} \label{lemma:Robust_Q-stability}
For a fixed $k$ and $0<\alpha_k\leq \frac{1}{2(L_g^2/\eta+\eta)}$, following LMI \eqref{eq:LMI_eps} holds for the dynamic system \eqref{eq:dynamic_sys-sep}.
\begin{align}\label{eq:LMI_eps}
    \begin{bmatrix} \lambda_k\ \bfP_k & 0 & 0 \\ 0 & s_k I & 0 \\ 0 & 0 & \omega_k I\end{bmatrix} \succeq
    \begin{bmatrix} M_{k,*} & \epsilon_k \  I & I \\ I & 0 & 0 \\ I & 0 & 0\end{bmatrix}^\top \begin{bmatrix} \bfP_k & 0 & 0\\ 0 & s_k I & 0 \\ 0 & 0 & I \end{bmatrix} 
    \begin{bmatrix} M_{k,*} & \epsilon_k \  I & I \\ I & 0 & 0 \\ I & 0 & 0\end{bmatrix}
\end{align}
where $\lambda_k = {1-\eta\alpha_k/2}$,\; $\bfP_k = C^{\bfP}_k \ \Diag(1/h^{k,*})$,\; $0< \epsilon_k \leq \overline \epsilon_k$,\; $s_k=4$,\;\; $\omega_k={3 C^\bfP_k}/{\nu^{\min}}+{(C^\bfP_k)^2 (1-\eta\alpha_k)^2}/{(\nu^{\min})^3}$,\; $\overline \epsilon_k = \min\left\{\sqrt{\frac{2 \nu^{\min}}{C^\bfP_k}},\frac{(\nu^{\min})^{3/2}}{2 C^\bfP_k (1-\eta\alpha_k)}\right\}$,\; $C^\bfP_k=(\frac{3}{2}s_k+1)/{\widetilde \mu_k}$ and $\widetilde \mu_k = \frac{{\eta\alpha_k/4}}{1+{16} (1-\eta\alpha_k)^2 L_g^2/{\eta^2}}>0$.
\end{lemma}
\proof{Proof of Lemma \ref{lemma:Robust_Q-stability}.}
See Appendix \ref{appendix:Proof-lemma:Robust_Q-stability}.
\Halmos \endproof

\begin{remark}
    Recall $\overline{\alpha}=\min\left\{\widetilde \alpha, {1}/{(2L_g^2/\eta+2\eta)}\right\}$.
     All results proved hereafter hold for $\alpha_k$ satisfying $0 < \alpha_k \leq \overline \alpha$.
\end{remark}
\subsubsection{Convergence rate of $\sfR_{k,t}$.}\label{appendix:ConvgRate-R_k,t} 
In this appendix, the aim is to derive the convergence rate of $\txteps_\sfR^{k,t}$ w.r.t. $t$ for a fixed $k$. This is given in Lemma \ref{lemma:ConvgRate-Jacobian} whose proof follows directly from key Lemma \ref{lemma:Stab_dyn} which in turn use Lemma \ref{lemma:ConvgRate-W_k,t} for its proof.

\begin{lemma}\label{lemma:ConvgRate-W_k,t}
For a fixed $k$, $\|\Delta_{k,t}\|$ and $\|\rmW_{k,t}\|$ in \eqref{eq:dynamic_sys-sep} decays exponentially fast to $0$ w.r.t. $t$ as
\begin{align*}
    \|\Delta_{k,t}\| \leq C^\Delta_k\ (\txteps_h^{k,0})^{1/2}\ (1- \eta \alpha_k)^{t/2},\qquad 
    \|\rmW_{k,t}\| \leq (C^U_k + \sfC_0 C^\Delta_k)\ (\txteps_h^{k,0})^{1/2}\ (1- \eta \alpha_k)^{t/2}
\end{align*}
where 
\newline $C^U_k = \sqrt{2} \alpha_k (\tau_g^h+\Lambda^{\xUL h}_g L_B (1- \eta \alpha_k)^{1/2})$,\ $C^\Delta_k = \sqrt{2} \left(L_{1/h} (1-\eta\alpha_k) + \rho_g^h \alpha_k + L_B (1-\eta\alpha_k)^{3/2}/\nu^{\min}\right)$.
\end{lemma}
\proof{Proof of Lemma \ref{lemma:ConvgRate-W_k,t}.}
See Appendix \ref{appendix:Proof-lemma:ConvgRate-W_k,t}.
\Halmos \endproof

\begin{lemma}\label{lemma:Stab_dyn}
For a fixed $k$, $\exists$ $T^0_k\geq0$ such that $\|\Delta_{k,t}\|\leq\epsilon_k \;\; \forall t\geq T^0_k$ and following bounds hold for the dynamic system \eqref{eq:dynamic_sys-sep}.
 \begin{enumerate}[label=\normalfont(\alph*)]     
 \item \label{lemma:Stab_dyn-a} $\begin{aligned}[t]
        \|\rmX_{k,t}\|^2 \leq {(\lambda_k)^{t-T^0_k}}/{\nu^{\min}}\ \|\rmX_{k,T^0_k}\|^2 + \frac{\theta_k\ \txteps_h^{k,0}}{1-\phi_k}\ \big(1-\phi_k^{t-T^0_k}\big)\ (\lambda_k)^{t-1},\quad t\geq T^0_k
       \end{aligned}$
   \item \label{lemma:Stab_dyn-b} $\begin{aligned}[t]
         \|\rmX_{k,t}\|^2 \leq {(\tlambda_k)^{t}}/{\nu^{\min}}\ \sfC_0^2 + \frac{\theta_k\ \txteps_h^{k,0}}{1-\tphi_k}\ \big(1-\tphi_k^{t}\big)\ (\tlambda_k)^{t-1},\quad 0\leq t \leq T^0_k
       \end{aligned}$       
 \end{enumerate}
where $\epsilon_k = \min\left\{C^\Delta_k(\sfD_\KL^{\max})^{1/2},\ \overline \epsilon_k\right\}$, $0<\lambda_k = 1-{\eta \alpha_k}/2<1$, $\tlambda_k=\lambda_k + \frac{(C^\Delta_k)^2 \sfD_\KL^{\max}-\epsilon_k^2}{\epsilon_k^2/4\ C^{\bfP}_k}$, $\theta_k = (C^U_k + \sfC_0 C^\Delta_k)^2\omega_k/C^{\bfP}_k$, $\phi_k={(1- \eta \alpha_k)}/{\lambda_k}$, $\tphi_k={(1- \eta \alpha_k)}/{\tlambda_k}$
 and\; $0<\tphi_k\leq\phi_k<1$. 
\end{lemma}
\proof{Proof of Lemma \ref{lemma:Stab_dyn}.}
See Appendix \ref{appendix:Proof-lemma:Stab_dyn}.
\Halmos \endproof

\proof{Proof of Lemma \ref{lemma:ConvgRate-Jacobian}.}
 Note the following bound on $\|\rmX_{k,T^0_k}\|$ from Lemma \ref{lemma:Stab_dyn}\ref{lemma:Stab_dyn-b}
 \[
 \|\rmX_{k,T^0_k}\|^2 \leq {(\tlambda_k)^{T^0_k}}/{\nu^{\min}}\ \sfC_0^2 + \frac{\theta_k\ \txteps_h^{k,0}}{1-\tphi_k}\ \big(1-\tphi_k^{T^0_k}\big)\ (\tlambda_k)^{T^0_k-1}
 \]
Substitute the above bound on $\|\rmX_{k,T^0_k}\|$ in Lemma \ref{lemma:Stab_dyn}\ref{lemma:Stab_dyn-a} and rearrange to get for $t\geq T^0_k$
\begin{align*}
& {\txteps_\sfR^{k,t}} \\
& \leq (\lambda_k)^{t-T^0_k} \frac{\sfC_0^2\ (\tlambda_k)^{T^0_k}}{(\nu^{\min})^2} + (\lambda_k)^{t-T^0_k}\ \theta_k\ \left((\tlambda_k)^{T^0_k-1}/\nu^{\min}\ \frac{1-(\tphi_k)^{T^0_k}}{1-\tphi_k} + (\lambda_k)^{T^0_k-1}\ \frac{1-(\phi_k)^{t-T^0_k}}{1-\phi_k}\right)\ {\txteps_h^{k,0}} \\ 
& \leq (\lambda_k)^{t-T^0_k} \frac{\sfC_0^2\ (\tlambda_k)^{T^0_k}}{(\nu^{\min})^2} + (\lambda_k)^{t-T^0_k}\ \theta_k\ \left(\frac{(\tlambda_k)^{T^0_k-1}/\nu^{\min}}{1-\tphi_k} + \frac{(\lambda_k)^{T^0_k-1}}{1-\phi_k}\right)\ {\txteps_h^{k,0}} 
\end{align*}
Take $\begin{aligned}[t]
\Gamma_{1k}=\frac{\sfC_0^2 (\tlambda_k/\lambda_k)^{T^0_k}}{(\nu^{\min})^2} \text{ and }  \Gamma_{2k}=\frac{\theta_k}{\lambda_k} \left(\frac{(\tlambda_k/\lambda_k)^{T^0_k-1}/\nu^{\min}}{1-\tphi_k} + \frac{1}{1-\phi_k}\right) 
\end{aligned}$ to complete the proof. \Halmos
\endproof

\begin{remark}
    Lemma \ref{lemma:Stab_dyn} implies that for a fixed $k$, the dynamic system \eqref{eq:dynamic_sys-sep} is not guaranteed to be stable for $0\leq t \leq T^0_k$ since $\tlambda_k$ may be greater than $1$ but converges linearly for $t\geq T^0_k$. This phenomenon is demonstrated in Example \ref{example:ToyExample} below.
\end{remark}

\begin{example}\label{example:ToyExample}
Consider $\mathcal{C}=\{\xUL\in\Real^{\sfq_u}:\|\xUL\|\leq1\}$, $\sfq_u=\sfq_\ell=\sum_{i\in[\sfW]}\sfq_\ell^i$, the lower level objective $g(h,\xUL) = \frac{1}{2}\ (h^\top \xUL)^2 = \frac{1}{2}\ h^\top \xUL \xUL^\top h  = \frac{1}{2}\ \xUL^\top h h^\top \xUL$ where $h=(h_i, i\in[\sfW])\in\Real^{\sfq_\ell}$ and $h_i\in\calH_i$.
\begin{align*}
     & \nabla_h g(h,\xUL) = (\xUL^\top h) \xUL,\quad \nabla^2_h g(h,\xUL) = \xUL\xUL^\top,\quad \nabla_\xUL \nabla_h g(h,\xUL) = \xUL^\top h \cdot I + h\xUL^\top 
\end{align*}
Note $\|\nabla_h g(h,\xUL)\| \leq \left\|\xUL\right\|^2\cdot \|h\|\leq 1$ and $\|\nabla_h g(h,\xUL)-\nabla_h g(h',\xUL)\|\leq \|\xUL\xUL^\top\|\cdot \|h-h'\| \leq \|h-h'\|$. This means $\Omega_g=1$, $L_g=1$. 
For some $\xUL^k$ chosen randomly on the unit circle, Figure \ref{fig:ToyExample} shows the plot of error dynamics ${\txteps_h^{k,t}}= \overline{\sfD}_\KL(h^{k,*},h^{k,t})$ and ${\txteps_\sfR^{k,t}}= \|\sfR_{k,t}-\sfR_{k,*}\|^2$ along with magnitude of maximum eigen value of $M_{k,t}$ and $\|M_{k,t}\|$. From Figure \ref{fig:ToyExample}(b), we observe that the error in $\sfR_{k,t}$ first increases until $t\leq T^0_k=460$ iterations and then decays to $0$. This is coherent with the stability analysis of dynamic system \eqref{eq:Jacb_Dyn_v2} which is unstable for $t\leq T^0_k$ as shown in Figure \ref{fig:ToyExample}(c) by eigen value(s) of $M_{k,t}$ outside unit circle for $t\leq T^0_k$ . Since $M_{k,t}$ is not necessarily a symmetric matrix, its spectral radius may not equal (and in general, is smaller than) its spectral norm. In fact for the above example, Figure \ref{fig:ToyExample}(c)--(d) show that the spectral radius is always smaller than the spectral norm of $M_{k,t}$. Moreover $\|M_{k,t}\|>1$ and as a result, the norm-based techniques to bound $\|\sfR_{k,t}-\sfR_{k,*}\|$ in order to analyze bilevel optimization algorithms do not apply more generally as they rely on $\|M_{k,t}\|$ to be less than $1$ (for instance \citet{ji2021bilevel} which considers unconstrained upper- and lower-level problems). Instead, we need to resort to LMI-based techniques as discussed in this paper. Specifically, Lemma \ref{lemma:Robust_Q-stability} and Lemma \ref{lemma:Stab_dyn} use concepts from robust stability \& control (\citet{skelton2013unified}) to analyze the convergence behavior of $\sfR_{k,t}$.

\begin{figure}[htbp]
   \begin{tabular}{ccrr}
     \includegraphics[width=0.24\textwidth]{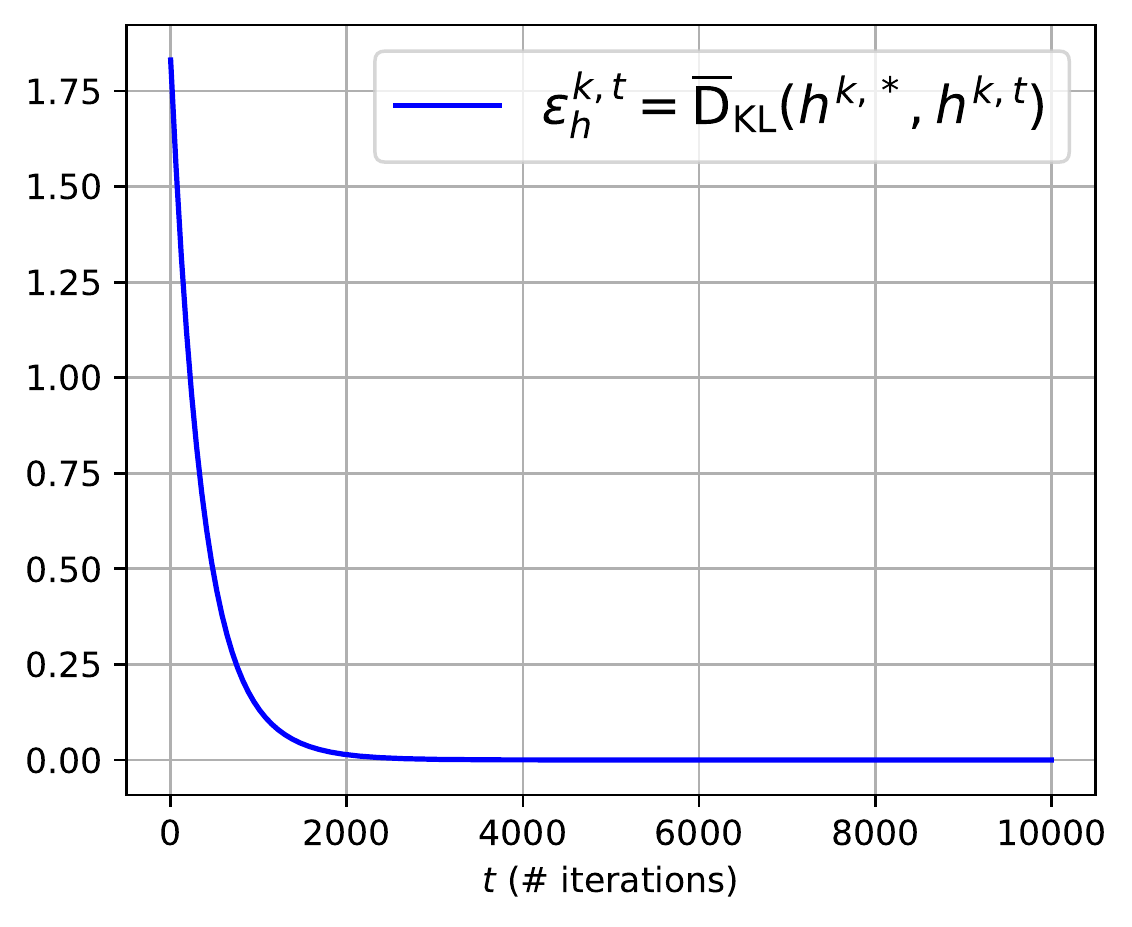} &  
     \includegraphics[width=0.225\textwidth]{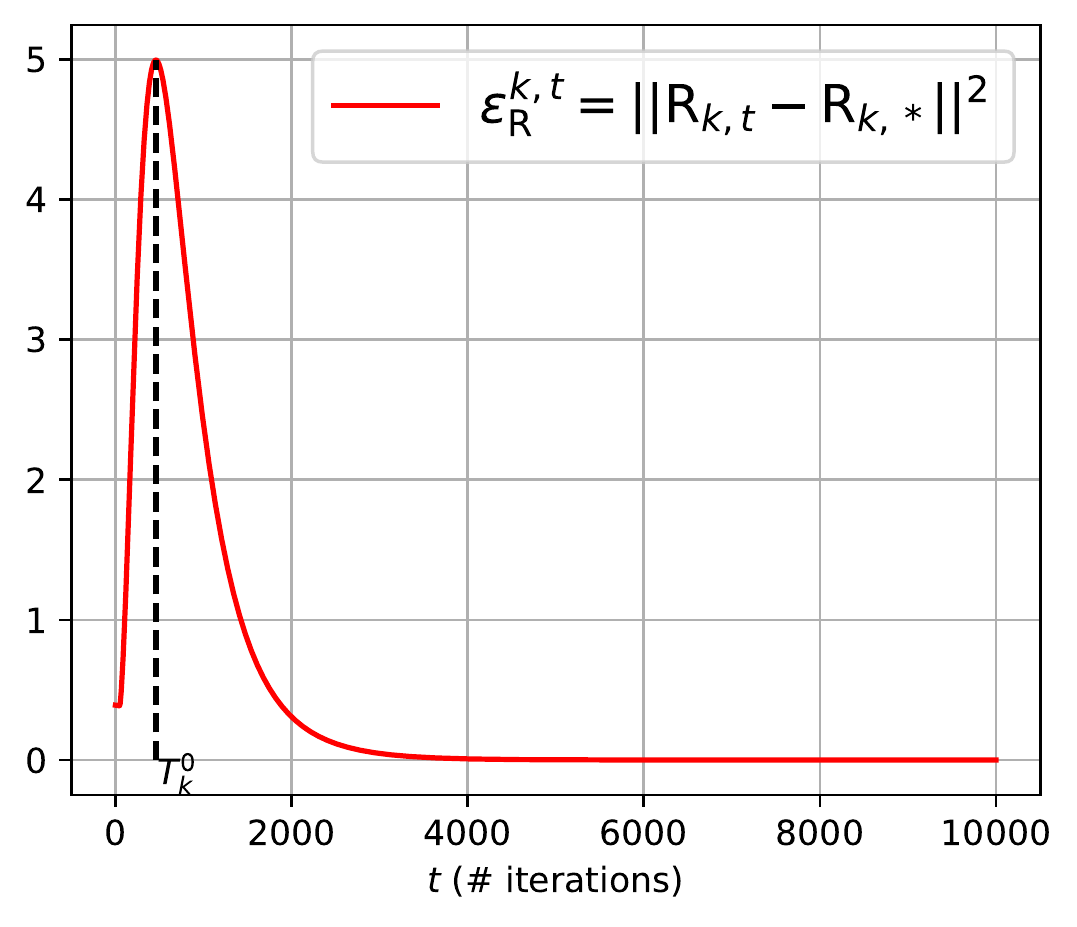} &
     \includegraphics[width=0.24\textwidth]{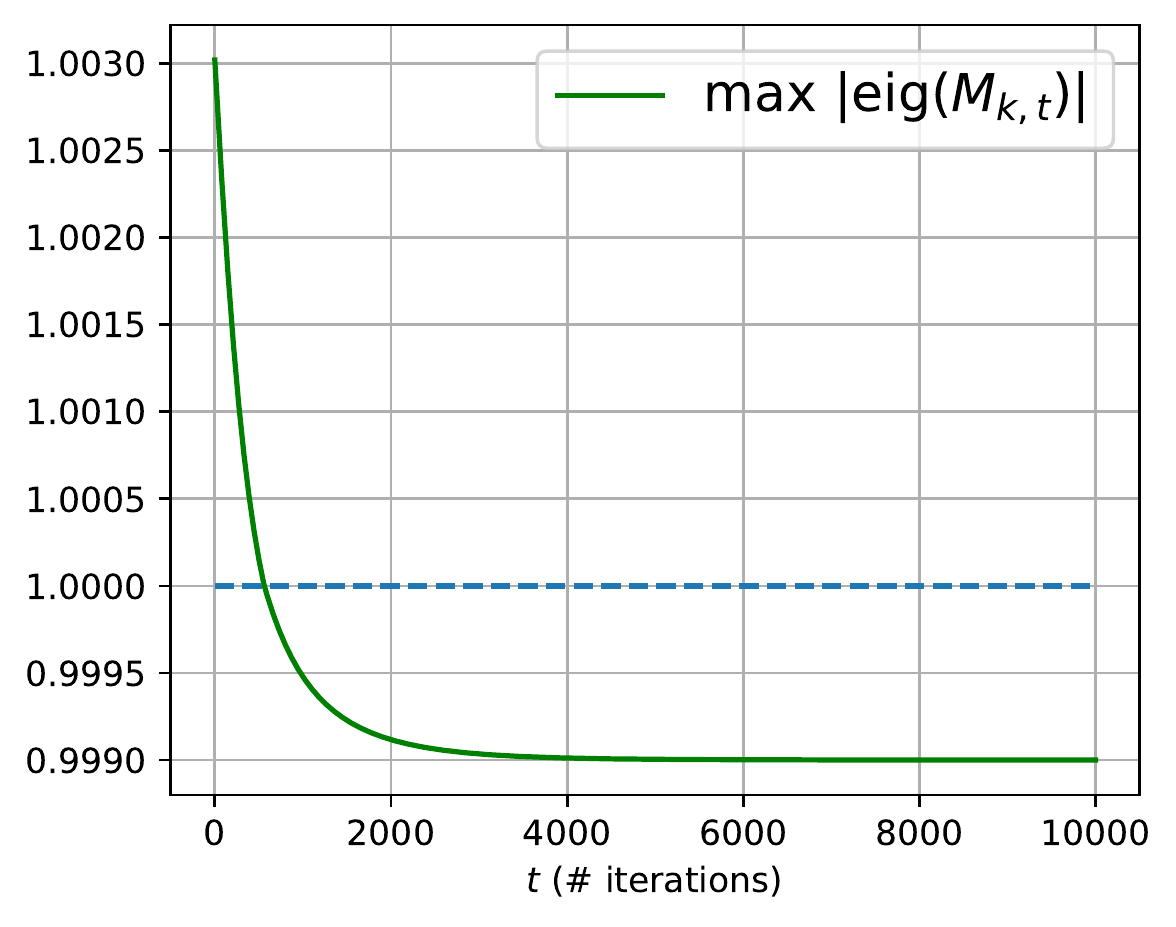} &
     \includegraphics[width=0.225\textwidth]{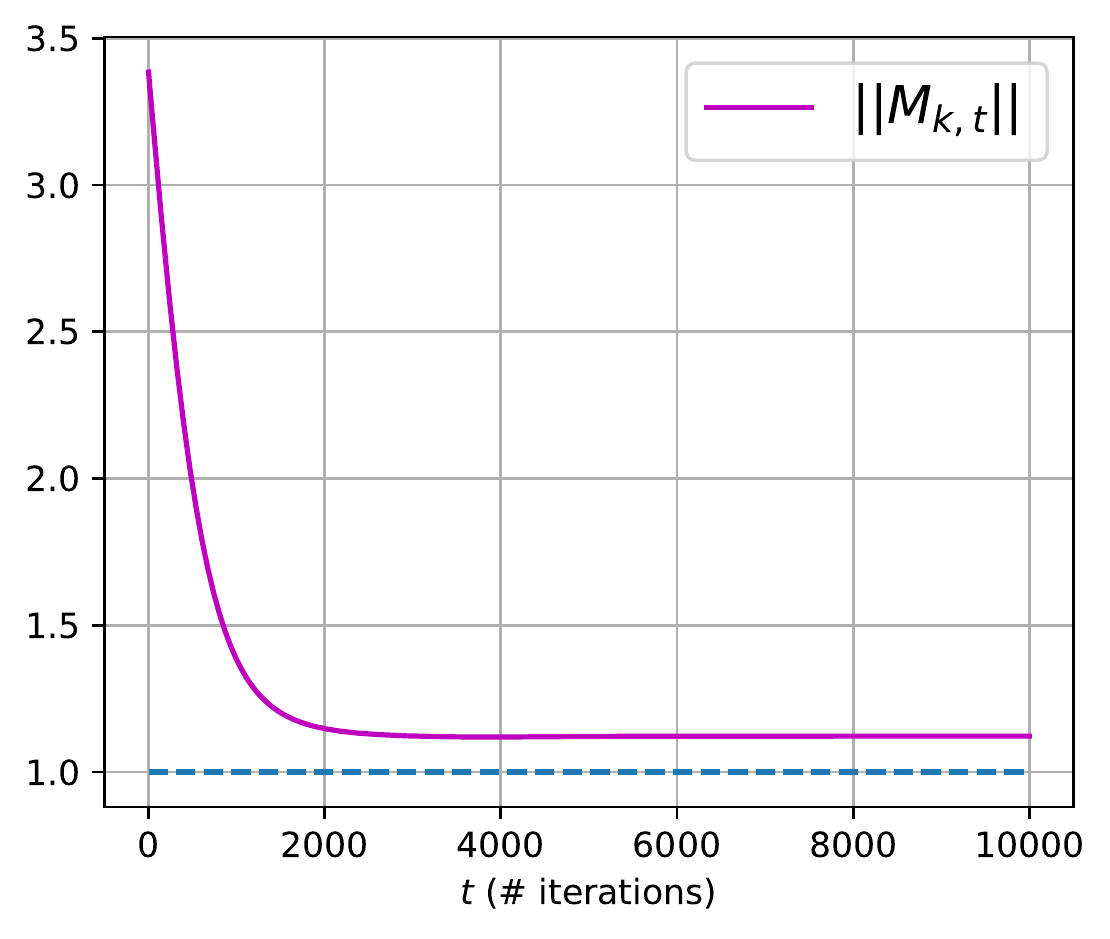}
     \\
      \scriptsize (a) & \scriptsize (b) & \scriptsize{(c) Spectral radius of $M_{k,t}$} & \scriptsize {(d) Spectral norm of $M_{k,t}$}
\end{tabular}
\caption{$\sfW=2$, \ $\sfq^{\max}_\ell=\sfq^{1}_\ell=\sfq^{2}_\ell=30$, \ $\eta=0.02$, \ $\alpha_k =0.05$} 
\label{fig:ToyExample}
\end{figure}
\end{example}

\section{Proof of various lemmas in Appendix \ref{appendix:ConvgRate-LowerLevel}.} \label{appendix:ProofsAppendixC}
\subsection{Proof of Lemma \ref{lemma:FixedPoint-h^k,t}.} \label{appendix:Proof-lemma:FixedPoint-h^k,t}
\begin{enumerate}[label=\normalfont(\alph*)]
\item
The Lagrangian of the lower level problem \eqref{eq:LowerLevel-EntropyReg}
\begin{align*}
    \mathcal{L}(h;\lambda,\sigma) = g^\eta(h,\xUL^k) + \sum_{i=1}^{\sfW}\lambda_i^k \Big(1-\sum_{j=1}^{\sfq_\ell^i} h_{i,j}\Big) - \sum_{i=1}^{\sfW}{\sigma_i^k}^\top h_i
\end{align*}
where $\sigma^k_i\geq0$. Since the lower level problem is convex and Slater's condition holds, the following KKT optimality conditions are necessary and sufficient. For each $i\in [\sfW]$
\begin{subequations}
\begin{align}
    & \sum_{j=1}^{\sfq_\ell^i} h^{k,*}_{i,j} = 1 \label{eq:KKT-1}\\
    & \nabla_{h_{i,j}}g^\eta(h^{k,*},\xUL^k) - \lambda^{k,*}_i - \sigma^{k,*}_{i,j} = 0 \quad \forall j\in [\sfq_\ell^i] \label{eq:KKT-2}\\
    & \sigma^{k,*}_{i,j}\cdot h^{k,*}_{i,j} = 0,\;\; h^{k,*}_{i,j}\geq 0,\;\; \sigma^{k,*}_{i,j}\geq 0 \quad \forall j\in [\sfq_\ell^i] \label{eq:KKT-3}
\end{align}
\end{subequations}
Let $\mathsf{S}_k(i)=\{j\in [\sfq_\ell^i]:\ h^{k,*}_{i,j}>0\}$ be the support of optimal solution $h^{k,*}_i$. This implies that for each $i\in [\sfW]$, $h^{k,*}_{i,j} = 0\;\forall j\notin \mathsf{S}_k(i)$ and by complementarity condition \eqref{eq:KKT-3}, $\sigma^{k,*}_{i,j} = 0\;\forall j\in \mathsf{S}_k(i)$. Plugging \eqref{eq:KKT-1}-\eqref{eq:KKT-3} into the RHS of \eqref{eq:FixedPoint}
\begin{align*}
    \frac{h_{i,j}^{k,*}\exp\left(-\alpha_k\nabla_{h_{i,j}}g^\eta(h^{k,*},\xUL^k) \right)}{\sum_{j'=1}^{\sfq_\ell^i} h_{i,j'}^{k,*}\exp\big(-\alpha_k\nabla_{h_{i,j'}}g^\eta(h^{k,*},\xUL^k) \big)} & \overset{\eqref{eq:KKT-2}}{=} \frac{h_{i,j}^{k,*}\exp\left(-\alpha_k (\lambda^{k,*}_i + \sigma^{k,*}_{i,j}) \right)}{\sum_{j'=1}^{\sfq_\ell^i} h_{i,j'}^{k,*}\exp\left(-\alpha_k(\lambda^{k,*}_i + \sigma^{k,*}_{i,j'}) \right)}\\
    & \overset{\eqref{eq:KKT-3}}{=}\begin{cases}
    0 & \textrm{if}\; j\notin \mathsf{S}_k(i)\\
    \frac{\displaystyle h_{i,j}^{k,*}\exp\big(-\alpha_k \lambda^{k,*}_i \big)}{\displaystyle{\mathsmaller \sum_{j'\in \mathsf{S}_k(i)}} h_{i,j'}^{k,*}\exp\big(-\alpha_k \lambda^{k,*}_i \big)} & \textrm{if}\; j\in \mathsf{S}_k(i)
    \end{cases}\\
    & \overset{\eqref{eq:KKT-1}}{=}\begin{cases}
    0 & \textrm{if}\; j\notin \mathsf{S}_k(i)\\
    {h_{i,j}^{k,*}} & \textrm{if}\; j\in \mathsf{S}_k(i) 
    \end{cases} \quad = \quad h^{k,*}_{i,j}
\end{align*}
\item Start with fixed point equation $h_i^{k,*}={h_i^{k,*}\circ \exp\left(-\alpha_k \nabla_{h_i}g^\eta(h^{k,*},\xUL^k)\right)}/{Z_i^{k,*}}$ for each $i\in[\sfW]$ where $Z_i^{k,*} = {h_i^{k,*}}^\top \exp\left(-\alpha_k \nabla_{h_i}g^\eta(h^{k,*},\xUL^k)\right)$. If $h^{k,*}\in\intr(\calH)$ then the fixed point equation implies
\begin{align}
    1 =  \exp\left(-\alpha_k \nabla_{h_{i,j}}g^\eta(h^{k,*},\xUL^k)\right)/{Z_i^{k,*}}\quad \forall j\in[\sfq^i_\ell] \label{eq:FixedPoint-reln}
\end{align}
Taking dual multipliers as $\lambda^{k,*}_i=-\ln(Z_i^{k,*})/{\alpha_k}$ and $\sigma^{k,*}_{i,j}=0$, the KKT conditions \eqref{eq:KKT-1}-\eqref{eq:KKT-3} can be satisfied.
\begin{align*}
    \nabla_{h_{i,j}}g^\eta(h^{k,*},\xUL^k) - \lambda^{k,*}_i - \sigma^{k,*}_{i,j} & = \nabla_{h_{i,j}}g^\eta(h^{k,*},\xUL^k) + \frac{1}{\alpha_k}\ln(Z_i^{k,*})\\
    & \overset{\eqref{eq:FixedPoint-reln}}{=} \nabla_{h_{i,j}}g^\eta(h^{k,*},\xUL^k) + \frac{1}{\alpha_k}\ln(e^{-\alpha_k \nabla_{h_{i,j}}g^\eta(h^{k,*},\xUL^k)})\\
    & =\ 0
\end{align*}
\end{enumerate}
This completes the proof. \Halmos
\endproof

\subsection{Proof of Lemma \ref{lemma:Iters_SmplxIntr}.}\label{appendix:Proof-lemma:Iters_SmplxIntr}
\begin{enumerate}[label=\normalfont(\alph*)]
\item  Need to prove for each $i$: if $h^{0,0}_i\in \widetilde \calH$ then $h^{k,t}_i\in \widetilde \calH_i\;\; \forall k,t$. The proof is by induction. For a given $k$, assume that for some $t$, $\exists$ $\varepsilon_i>0$ such that $h_i^{k,t}\geq \varepsilon_i\cdot \mathbf{1}_{\sfq_\ell^i}$ $\forall i\in[\sfW]$. Then for each $i\in[\sfW]$ and for every $j\in[\sfq_\ell^i]$
\begin{align*}
    h_{i,j}^{k,t+1}:=\frac{h_{i,j}^{k,t} e^{-\alpha_k \nabla_{h_{i,j}}g^\eta(h^{k,t},\xUL^k)}} {\sum_{j'=1}^{\sfq_\ell^i} h_{i,j'}^{k,t} e^{-\alpha_k \nabla_{h_{i,j'}}g^\eta(h^{k,t},\xUL^k)}}
\end{align*}
Using $\nabla_{h_{i,j}}g^\eta(h^{k,t},\xUL^k) = \nabla_{h_{i,j}}g(h^{k,t},\xUL^k) + \eta\ \ln(h^{k,t}_{i,j})$ and the assumption $h^{k,t}_i>0$,
\begin{align} \label{eq:Interm-LB-h^{k,t}}
    & h_{i,j}^{k,t+1} \nonumber\\ & = \frac{(h_{i,j}^{k,t})^{1-\eta\alpha_k}\ e^{-\alpha_k\nabla_{h_{i,j}}g(h^{k,t},\xUL^k)}} {\sum_{j'=1}^{\sfq_\ell^i} (h_{i,j'}^{k,t})^{1-\eta\alpha_k}\ e^{-\alpha_k\nabla_{h_{i,j'}}g(h^{k,t},\xUL^k)}} \geq e^{-2\alpha_k\Omega_g} \frac{(h_{i,j}^{k,t})^{1-\eta\alpha_k}} {\sum_{j'=1}^{\sfq_\ell^i} (h_{i,j'}^{k,t})^{1-\eta\alpha_k}} \geq e^{-2\alpha_k\Omega_g} \frac{(\varepsilon_i)^{1-\eta\alpha_k}} {\sum_{j'=1}^{\sfq_\ell^i} (h_{i,j'}^{k,t})^{1-\eta\alpha_k}}
\end{align}
We can globally upper bound the term $\sum_{j'=1}^{\sfq_\ell^i} (h_{i,j'}^{k,t})^{1-\eta\alpha_k}$ in the denominator over the simplex $\mathcal{H}_i$. Note the fact that for $0<\eta\alpha_k<1$, the function $f(\rmx)=\sum_{j=1}^{\sfq_\rmx} (\rmx_{j})^{1-\eta\alpha_k}$ is concave for $\rmx\geq0$. Therefore its maximum over a simplex can be found by solving
\[\max\ \sum_{j=1}^{\sfq_\rmx} (\rmx_{j})^{1-\eta\alpha_k} \quad \st \quad \sum_{j=1}^{\sfq_\rmx} \rmx_{j} = 1, \;\; \rmx \geq 0 \]
Relaxing non-negativity constraint and writing the Lagrangian we have
\begin{align*}
    &\mathcal{L}(\rmx,\kappa)=\sum_{j=1}^{\sfq_\rmx} (\rmx_{j})^{1-\eta\alpha_k} + \kappa\Big(1-\sum_{j=1}^{\sfq_\rmx}\rmx_{j}\Big),\qquad  \nabla_{\rmx_{j}} \mathcal{L}(\rmx^{*},\kappa^*) = 0 \implies (1-\eta\alpha_k)\ (\rmx^{*}_{j})^{-\eta\alpha_k} = \kappa^*
\end{align*}
This means $\rmx^{*}_{j}$ need to be equal $\forall\ j\in[\sfq_\ell^i]$ and plugging it into simplex constraint implies $\rmx^{*}_{j}=\frac{1}{\sfq_\rmx}$ which satisfies the non-negativity constraint. Replace $\rmx\to h_i$ and $\sfq_\rmx\to\sfq_\ell^i$ to get
\[\textstyle{\sum_{j=1}^{\sfq_\ell^i}} (h_{i,j})^{1-\eta\alpha_k} \leq (\sfq_\ell^i)^{\eta\alpha_k} \;\; \forall\ h_i\in\mathcal{H}_i\]
Substituting this result back into \eqref{eq:Interm-LB-h^{k,t}} gives
\begin{align} \label{eq:Interm2-LB-h^{k,t}}
    h_{i,j}^{k,t+1} \geq e^{-2\alpha_k\Omega_g} \frac{(\varepsilon_i)^{1-\eta\alpha_k}}{(\sfq_\ell^i)^{\eta\alpha_k}}
\end{align}
Now this enables to find a condition on $\varepsilon_i$ such that if
\[e^{-2\alpha_k\Omega_g} \frac{(\varepsilon_i)^{1-\eta\alpha_k}}{(\sfq_\ell^i)^{\eta\alpha_k}} \geq \varepsilon_i \implies \varepsilon_i \leq e^{-\frac{2\Omega_g}{\eta}}/\sfq_\ell^i\]
then $h_i^{k,t+1} \geq \varepsilon_i\cdot \mathbf{1}_{\sfq_\ell^i}$. For any given $k$, doing induction on $t$ gives: \textit{if $h_i^{k,0}\geq\varepsilon_i\cdot \mathbf{1}_{\sfq_\ell^i}$ then $h_i^{k,t}\geq \varepsilon_i\cdot \mathbf{1}_{\sfq_\ell^i}\ \forall t$}. Since $h^{k,0}=h^{k-1,D}$ for $k>0$ in Algorithm \ref{alg:Bilevel_ITD}, then combining previous induction on $t$ with doing induction on $k$ gives:  \textit{if $h_i^{0,0}\geq\varepsilon_i\cdot \mathbf{1}_{\sfq_\ell^i}$ then $h_i^{k,t}\geq \varepsilon_i\cdot \mathbf{1}_{\sfq_\ell^i}\ \forall k,t$}. The proof is complete by taking $\varepsilon_i=e^{-\frac{2\Omega_g}{\eta}}/\sfq_\ell^i$ and defining $\nu^{\min}=\min_{i\in[\sfW]} \varepsilon_i$.

\item Note that $\nabla_{h_{i,j}}g^\eta(h^{k,*},\xUL^k) = \nabla_{h_{i,j}}g(h^{k,*},\xUL^k) + \eta\ln(h^{k,*}_{i,j})$. For any $h=(h_i,\ i\in[\sfW])$ where $h_i\in\mathcal{H}_i$, following holds by optimality condition
\[
\sum_{i=1}^{\sfW} \sum_{j=1}^{\sfq_\ell^i} \left[\nabla_{h_{i,j}}g(h^{k,*},\xUL^k) + \eta\ln(h^{k,*}_{i,j})\right]\cdot(h_{i,j}-h_{i,j}^{k,*}) \geq 0
\]
Assume there exists $j_1, j_2\in[\sfq_\ell^i]$ for some $i$ such that $h_{i,j_1}^{k,*}=0$, $h_{i,j_2}^{k,*}>0$. Set $h_{i,j_1}=h_{i, j_2}^{k,*}>0$, $h_{i,j_2}=0$, $h_{i,j}=h_{i,j}^{k,*}\ j\neq j_1,j_2$. Then LHS in the optimality condition can be made to go $-\infty$ since $\nabla_{h_{i,j}}g(h^{k,*},\xUL^k)$ is assumed to be bounded and the product $\ln{0}\cdot(h_{i,j_1}-0)$ $\to-\infty$. We have a contradiction here implying all the entries $h_i^{k,*}$ must be strictly greater than 0 i.e. $h_i^{k,*}>0$.

Since $h_i^{k,*}>0$, from \eqref{eq:KKT-3} it holds $\sigma^{k,*}_{i,j}=0\ \forall j\in [\sfq_\ell^i]$ and using \eqref{eq:KKT-2},
\[\nabla_{h_{i,j}}g(h^{k,*},\xUL^k) + \eta\ln(h^{k,*}_{i,j}) - \lambda^{k,*}_i = 0 \implies h^{k,*}_{i,j} = e^{\frac{1}{\eta} \lambda^{k,*}_i} \cdot e^{-\frac{1}{\eta} \nabla_{h_{i,j}}g(h^{k,*},\xUL^k)}\]

Now using $\sum_{j=1}^{\sfq_\ell^i} h^{k,*}_{i,j}=1$ to lower bound $\lambda^{k,*}_i$.
\[
1 = e^{\frac{1}{\eta} \lambda^{k,*}_i}\ \textstyle{\sum_{j=1}^{\sfq_\ell^i}} e^{-\frac{1}{\eta} \nabla_{h_{i,j}}g(h^{k,*},\xUL^k)} \leq e^{\frac{1}{\eta} \lambda^{k,*}_i}\ \sfq_\ell^i e^{\frac{\Omega_g}{\eta} } \implies -\eta \ln \sfq_\ell^i - \Omega_g \leq \lambda_i^*
\]

As a result
\[h^{k,*}_{i,j} \geq e^{\frac{1}{\eta} \lambda^{k,*}_i} \cdot e^{-\frac{\Omega_g}{\eta} } \geq e^{-\frac{2\Omega_g}{\eta} }/\sfq_\ell^i\geq \nu^{\min}\]
\end{enumerate}    
The proof is complete by noting that $h^{k,*}_{i}=(h^{k,*}_{i,j},j\in[\sfq_\ell^i])\in\widetilde\calH_i$ for each $i$.  \Halmos
\endproof

\subsection{Proof of Lemma \ref{lemma:Dynamics-R_k,t}.}\label{appendix:Proof-lemma:Dynamics-R_k,t}\ \\
Differentiate both sides of dynamics \eqref{eq:Proj_MD} (rewritten as \eqref{eq:Proj_MDv2}) w.r.t $\xUL^k$
\begin{align}\label{eq:Proj_MDv2}
    h_i^{k,t+1}={h_i^{k,t}\circ \exp\left(-\alpha_k \nabla_{h_i}g^\eta(h^{k,t},\xUL^k)\right)}/{Z_i^{k,t}} \quad \forall i\in[\sfW]
\end{align}
where $Z_i^{k,t} =\mathbf{1}^\top(h_i^{k,t}\circ \exp\left(-\alpha_k \nabla_{h_i}g^\eta(h^{k,t},\xUL^k)\right))$.
Then for each $i\in[\sfW]$, we have\\
\begin{align}\label{eq:JacbDynm_i}
    & \frac{\partial h_i^{k,t+1}}{\partial \xUL^k} = \frac{1}{Z_i^{k,t}} \frac{\partial }{\partial \xUL^k} (h_i^{k,t}\circ \exp\left(-\alpha_k \nabla_{h_i}g^\eta(h^{k,t},\xUL^k)\right) - \frac{h_i^{k,t}\circ \exp\left(-\alpha_k \nabla_{h_i}g^\eta(h^{k,t},\xUL^k)\right)}{(Z_i^{k,t})^2} \frac{\partial Z_i^{k,t}}{\partial \xUL^k}
\end{align}
Calculating the derivatives needed in the expression \eqref{eq:JacbDynm_i}
\begin{align*}
    & \frac{\partial }{\partial \xUL^k} (h_i^{k,t}\circ \exp\left(-\alpha_k \nabla_{h_i}g^\eta(h^{k,t},\xUL^k)\right) \\
    & = \Diag(e^{-\alpha_k \nabla_{h_i}g^\eta(h^{k,t},\xUL^k)})\frac{\partial h_i^{k,t}}{\partial \xUL^k} \\
    & \quad -\alpha_k \Diag(h_i^{k,t}\circ e^{-\alpha_k \nabla_{h_i}g^\eta(h^{k,t},\xUL^k)})\left(\nabla_\xUL\nabla_{h_i} g^\eta(h^{k,t},\xUL^k)+\nabla_h\nabla_{h_i}g^\eta(h^{k,t},\xUL^k)\frac{\partial h^{k,t}}{\partial \xUL^k}\right)
\end{align*}
\begin{align*}
    \frac{\partial Z_i^{k,t}}{\partial \xUL^k} & = \mathbf{1}^\top \frac{\partial }{\partial \xUL^k} (h_i^{k,t}\circ \exp\left(-\alpha_k \nabla_{h_i}g^\eta(h^{k,t},\xUL^k)\right) \\
    & = {e^{-\alpha_k \nabla_{h_i}g^\eta(h^{k,t},\xUL^k)}}^\top\frac{\partial h_i^{k,t}}{\partial \xUL^k} \\
    & \quad -\alpha_k ({h_i^{k,t}\circ e^{-\alpha_k \nabla_{h_i}g^\eta(h^{k,t},\xUL^k)}})^\top\left(\nabla_\xUL\nabla_{h_i} g^\eta(h^{k,t},\xUL^k)+\nabla_h\nabla_{h_i}g^\eta(h^{k,t},\xUL^k)\frac{\partial h^{k,t}}{\partial \xUL^k}\right)
\end{align*}
Substituting above expressions in \eqref{eq:JacbDynm_i} and then using $\exp\left(-\alpha_k \nabla_{h_i}g^\eta(h^{k,t},\xUL^k)\right)=h_i^{k,t+1}/h_i^{k,t}$ from \eqref{eq:Proj_MDv2} as $h^{k,t}\in\widetilde\calH$ (from Lemma \ref{lemma:Iters_SmplxIntr}) 
to get
\begin{align*}
    & \frac{\partial h_i^{k,t+1}}{\partial \xUL^k} = \begin{aligned}[t] & \Diag(h_i^{k,t+1}/h_i^{k,t})\frac{\partial h_i^{k,t}}{\partial \xUL^k}
    - h_i^{k,t+1} (h_i^{k,t+1}/h_i^{k,t})^\top\frac{\partial h_i^{k,t}}{\partial \xUL^k} \\
    & - \alpha_k \left(\Diag(h_i^{k,t+1})- h_i^{k,t+1} {h_i^{k,t+1}}^\top\right)\left(\nabla_\xUL\nabla_{h_i} g^\eta(h^{k,t},\xUL^k)+\nabla_h\nabla_{h_i}g^\eta(h^{k,t},\xUL^k)\frac{\partial h^{k,t}}{\partial \xUL^k}\right)
    \end{aligned}
\end{align*}
Noting $\Diag(h_i^{k,t+1}/h_i^{k,t})=\Diag(h_i^{k,t+1})\Diag(1/h_i^{k,t})$ and $(h_i^{k,t+1}/h_i^{k,t})^\top={h_i^{k,t+1}}^\top \Diag(1/h_i^{k,t})$ gives
\begin{align}\label{eq:JacbDynm_i-v2}
    \frac{\partial h_i^{k,t+1}}{\partial \xUL^k} = B^i_{k,t} \Diag(1/h_i^{k,t}) \frac{\partial h_i^{k,t}}{\partial \xUL^k} - \alpha_k B^i_{k,t} \left(\nabla_\xUL\nabla_{h_i} g^\eta(h^{k,t},\xUL^k)+\nabla_h\nabla_{h_i}g^\eta(h^{k,t},\xUL^k)\frac{\partial h^{k,t}}{\partial \xUL^k}\right)
\end{align}
where $B^i_{k,t}=\Diag(h_i^{k,t+1})- h_i^{k,t+1} {h_i^{k,t+1}}^\top$. To find $\sfR_{k,t+1}:={\partial h^{k,t+1}}/{\partial \xUL^k}=\begin{bmatrix} {\partial h_1^{k,t+1}}/{\partial \xUL^k}\\ \vdots\\ {\partial h_\sfW^{k,t+1}}/{\partial \xUL^k}\end{bmatrix}$ is equivalent to stacking \eqref{eq:JacbDynm_i-v2} as rows for $i=1,\hdots,\sfW$ and using $\nabla^2_h g^\eta(h^{k,t},\xUL^k) = \nabla^2_h g(h^{k,t},\xUL^k) + \eta \Diag(1/h)$, $\nabla_\xUL \nabla_{h} g^\eta(h,\xUL)=\nabla_\xUL \nabla_{h} g(h,\xUL)$ leads to \eqref{eq:Jacb_Dyn_v2}. 
This completes the proof. \Halmos
\endproof

\subsection{Proof of Lemma \ref{lemma:FixedPoint-R_k,t}.}\label{appendix:Proof-lemma:FixedPoint-R_k,t}\ \\
To find $\sfR_{k,*}$ we need to differentiate following fixed point equation w.r.t $\xUL^k$
\begin{align}
    h_i^{k,*}={h_i^{k,*}\circ \exp\left(-\alpha_k \nabla_{h_i}g^\eta(h^{k,*},\xUL^k)\right)}/{Z_i^{k,*}} \quad \forall i\in[\sfW]
\end{align}
which leads to \eqref{eq:Jacb_Dyn-FixedPoint_v2} following essentially the same steps as in the proof of Lemma \ref{lemma:Dynamics-R_k,t}. The existence of a matrix $\sfR_{k,*}$ satisfying \eqref{eq:Jacb_Dyn-FixedPoint_v2} can be argued using implicit function theorem. This completes the proof. \Halmos
\endproof

\subsection{Proof of Lemma \ref{lemma:B_k,t factor}.}\label{appendix:Proof-lemma:B_k,t factor} \ \\
The matrix $B^i = \Diag\big(h_i\big) - h_i {h_i}^\top$
can be rewritten as 
\begin{align} \label{eq:B_it-form1}
    B^i = \Diag(\sqrt{h_i})\ \big(I - \sqrt{h_i} {\sqrt{h_i}}^\top\big)\ \Diag(\sqrt{h_i})
\end{align}
where $\sqrt{h_i}$ is a unit vector since $\big\|\sqrt{h_i}\big\|_2^2=\sqrt{h_i}^\top\sqrt{h_i}=\|h_i\|_1=1$. Then
\[0\preceq I - \sqrt{h_i} {\sqrt{h_i}}^\top\preceq I \;\; \implies \;\; 0\preceq \Diag(\sqrt{h_i})\ \big(I - \sqrt{h_i} {\sqrt{h_i}}^\top\big)\ \Diag(\sqrt{h_i})\preceq \Diag(h_i) \preceq I\]
which proves $\|B^i\|\leq 1$. Further by Gram-Schmidt orthonormalization process, one can find vectors $\sfv^{i,j}\in\Real^{\sfq_\ell^i}$, $j=1,\hdots,\sfq_\ell^i-1$ or equivalently matrix $\sfV^i = [\sfv^{i,1},\hdots,\sfv^{i,\sfq_\ell^i-1}]\in\Real^{\sfq_\ell^i\times(\sfq_\ell^i-1)}$ such that ${\sfV^i}^\top \sfV^i = I$ and ${\sqrt{h_i}}^\top \sfV^i = 0$. Compactly,
\begin{align*}
    \begin{bmatrix} {\sqrt{h_i}} & \sfV^i \end{bmatrix}^\top \begin{bmatrix} {\sqrt{h_i}} & \sfV^i \end{bmatrix} =
    \begin{bmatrix} {\sqrt{h_i}}^\top \sqrt{h_i} &  {\sqrt{h_i}}^\top \sfV^i \\
    {\sfV^i}^\top \sqrt{h_i} &  {\sfV^i}^\top \sfV^i
    \end{bmatrix} = 
    \begin{bmatrix} 1 & 0 \\ 0 & I \end{bmatrix} = I_{\sfq_\ell^i\times\sfq_\ell^i}
\end{align*}
By construction $\begin{bmatrix} {\sqrt{h_i}} & \sfV^i \end{bmatrix}$ is invertible. As a result, $\begin{bmatrix} {\sqrt{h_i}} & \sfV^i \end{bmatrix}^\top=\begin{bmatrix} {\sqrt{h_i}} & \sfV^i \end{bmatrix}^{-1}$ implying \[\begin{bmatrix} {\sqrt{h_i}} & \sfV^i \end{bmatrix}\begin{bmatrix} {\sqrt{h_i}} & \sfV^i \end{bmatrix}^\top=\sqrt{h_i} {\sqrt{h_i}}^\top + \sfV^i {\sfV^i}^\top = I \quad \text{or} \quad I-\sqrt{h_i} {\sqrt{h_i}}^\top = \sfV^i {\sfV^i}^\top\]
Therefore, from \eqref{eq:B_it-form1}
\[B^i = \Diag(\sqrt{h_i})\ \sfV^i {\sfV^i}^\top\ \Diag(\sqrt{h_i})={\Lambda^i}{\Lambda^i}^\top\]
Recall $B=\blk\left(B^1,\hdots,B^\sfW\right)$ which readily gives $\|B\|\leq1$ and 
\begin{align*}
    B & = \blk\Big(\Diag(\sqrt{h_1}) \sfV^1 {\sfV^1}^\top \Diag(\sqrt{h_1})\;\;,\;\;\hdots \;\; ,\;\;\Diag(\sqrt{h_\sfW}) \sfV^\sfW {\sfV^\sfW}^\top \Diag(\sqrt{h_\sfW})\Big) \\
    & = \Diag(\sqrt{h})\ \blk\big( \sfV^1,\hdots,\sfV^\sfW\big)\ \blk\big({\sfV^1}^\top,\hdots,{\sfV^\sfW}^\top\big)\ \Diag(\sqrt{h}) \\
    & = \Lambda \Lambda^\top
\end{align*}
The proof is complete. \Halmos
\endproof

\subsection{Proof of Lemma \ref{lemma:M_k,*}.}\label{appendix:Proof-lemma:M_k,*} \ \\ 
\begin{enumerate}[label=\normalfont(\alph*)]
\item $\begin{aligned}[t]
& \Lambda_{k,*}^\top \Lambda_{k,*} \nonumber = \sfV_{k,*}^\top\ \Diag^2(\sqrt{h^{k,*}}) \sfV_{k,*} 
\end{aligned}$\\
Since $\Diag^2(\sqrt{h^{k,*}})=\Diag(h^{k,*})$, $\nu^{\min}\cdot I\preceq \Diag(h^{k,*}) \preceq I$ and $\sfV_{k,*}^\top\sfV_{k,*}=I$, therefore it holds that $\nu^{\min}\cdot I\preceq \Lambda_{k,*}^\top \Lambda_{k,*}\preceq I$. \\ \\
$\begin{aligned}[t]
& \Lambda_{k,*}^\top \Diag\left(1/h^{k,*}\right) \Lambda_{k,*} = \sfV_{k,*}^\top \Diag(\sqrt{h^{k,*}}) \Diag\left(1/h^{k,*}\right)\Diag(\sqrt{h^{k,*}})\sfV_{k,*}
\end{aligned}$ \\
Since $\Diag(\sqrt{h^{k,*}}) \Diag\left(1/h^{k,*}\right)\Diag(\sqrt{h^{k,*}})=I$ and $\sfV_{k,*}^\top\sfV_{k,*}=I$, therefore it holds that $\Lambda_{k,*}^\top \Diag\left(1/h^{k,*}\right) \Lambda_{k,*} = I$.\\ \\
$\begin{aligned}[t]
& \Lambda_{k,*}^\top \Diag^2\left(1/h^{k,*}\right) \Lambda_{k,*} = \sfV_{k,*}^\top \Diag(\sqrt{h^{k,*}}) \Diag^2\left(1/h^{k,*}\right)\Diag(\sqrt{h^{k,*}})\sfV_{k,*}
\end{aligned}$\\
Since $\Diag(\sqrt{h^{k,*}}) \Diag^2(1/h^{k,*})\Diag(\sqrt{h^{k,*}})=\Diag(1/h^{k,*})$, $ I\preceq \Diag(1/h^{k,*}) \preceq 1/\nu^{\min}\cdot  I$ and $\sfV_{k,*}^\top\sfV_{k,*}=I$, therefore it holds that $I\preceq \Lambda_{k,*}^\top \Diag^2\left(1/h^{k,*}\right) \Lambda_{k,*}\preceq 1/\nu^{\min}\cdot  I$.\\ \\
$\begin{aligned}[t]
&H_{k,*}^\top H_{k,*} = \blk(\|h^{*}_1\|^2,\hdots,\|h^{*}_\sfW\|^2) \preceq I
\end{aligned}$ \\ \\
$\begin{aligned}[t]
& \Lambda_{k,*}^\top \Diag\left(1/h^{k,*}\right) H_{k,*} = \sfV_{k,*}^\top \Diag(\sqrt{h^{k,*}}) \Diag\left(1/h^{k,*}\right)\Diag(\sqrt{h^{k,*}})\sqrt{H_{k,*}}
\end{aligned}$ \\
Since $\Diag(\sqrt{h^{k,*}}) \Diag\left(1/h^{k,*}\right)\Diag(\sqrt{h^{k,*}})=I$ and $\sfV_{k,*}^\top\sqrt{H_{k,*}}=0$, therefore it holds that $\Lambda_{k,*}^\top \Diag\left(1/h^{k,*}\right) H_{k,*} = 0$. 
\item 
Since $0\preceq \nabla^2_{h} g(h^{k,*},\xUL^k) \preceq L_g \cdot I$ it implies $0 \preceq \Lambda_{k,*}^\top\nabla^2_{h} g(h^{k,*},\xUL^k)\Lambda_{k,*} \preceq L_g \Lambda_{k,*}^\top \Lambda_{k,*}$,
Using the result for $\Lambda_{k,*}^\top \Lambda_{k,*}$ in part (a) above, following holds
\[
(1-(L_g+\eta)\alpha_k) I \preceq \widehat M_{k,*} = (1-\eta\alpha_k) I - \alpha_k \Lambda_{k,*}^\top \nabla^2_{h} g(h^{k,*},\xUL^k) \Lambda_{k,*} \preceq (1-\eta\alpha_k) I
\]

\item Using Weinstein–Aronszajn identity given in \citet{pozrikidis2014introduction}
on form \eqref{eq:M_k,*-v2} of $M_{k,*}$: For non-zero eigen values
\begin{align*}
\eig(M_{k,*}) 
& = \eig\left(\Lambda_{k,*}^\top \big[(1-\eta\alpha_k) \Diag\left(1/h^{k,*}\right)-\alpha_k \nabla^2_{h} g(h^{k,*},\xUL^k)\big] \Lambda_{k,*}\right) \overset{\text{part (a)}}{=} \eig(\widehat M_{k,*}) 
\end{align*}
The proof is complete by directly applying the result in part (b). \Halmos
\end{enumerate}
\endproof

\subsection{Proof of Lemma \ref{lemma:Robust_Q-stability}.}\label{appendix:Proof-lemma:Robust_Q-stability}
The dynamic system \eqref{eq:dynamic_sys-sep} is robust $\mathcal{Q}$-stable if $\exists$ matrix $\bfP_k\succ0$, scalars $\epsilon_k$, $s_k$, $\omega_k>0$ and $\lambda_k\in(0,1)$ such that the following LMI holds \citet{skelton2013unified}. 
\begin{subequations}\label{eq:LMI_Qstability}
\begin{align}
    \eqref{eq:LMI_eps} \iff & \begin{bmatrix} \lambda_k \bfP_k - M_{k,*}^\top \bfP_k M_{k,*} - (s_k+1) I  & \quad & - \epsilon_k M_{k,*}^\top \bfP_k  & \quad &  -M_{k,*}^\top \bfP_k \\ -\epsilon_k \bfP_k M_{k,*} & \quad & s_k I-\epsilon_k^2 \bfP_k  & \quad &  -\epsilon_k \bfP_k \\ -\bfP_k M_{k,*} & \quad & -\epsilon_k \bfP_k  & \quad &  \omega_k I-\bfP_k \end{bmatrix} \succeq 0 \nonumber \\ 
    \overset{\substack{\text{Schur}\\ 
    \text{complement}}}{\iff} & \mathcal{M}_k 
    = \begin{bmatrix} \lambda_k \bfP_k - M_{k,*}^\top \bfP_k M_{k,*} - (s_k+1) I  & & - \epsilon_k M_{k,*}^\top \bfP_k \\ -\epsilon_k \bfP_k M_{k,*} & & s_k I-\epsilon_k^2 \bfP_k \end{bmatrix} \succ 0 \label{eq:LMI_Qstability-a} \\
    & \omega_k I-\bfP_k - \bfP_k \begin{bmatrix} M_{k,*} & \epsilon_k I\end{bmatrix} \mathcal{M}_k^{-1} \begin{bmatrix} M_{k,*}^\top  \\ \epsilon_k I \end{bmatrix} \bfP_k \succeq 0 \label{eq:LMI_Qstability-b}
\end{align}
\end{subequations}
\underline{Proving that there exists $\bfP_k \succ 0$ and $\lambda_k\in(0,1)$ such that $\lambda_k \bfP_k - M_{k,*}^\top \bfP_k M_{k,*} \succ 0$}: \quad \\
\newline Recall $\Lambda_{k,*}=\Diag(\sqrt{h^{k,*}})\ \sfV_{k,*}$ where $\sfV_{k,*}=\blk\big({\sfV^1_{k,*}},\hdots,{\sfV^\sfW_{k,*}}\big)$,  $H_{k,*}=\blk(h^{k,*}_1,\hdots,h^{k,*}_\sfW)$, and $\begin{bmatrix} \sfV_{k,*}  & \sqrt{H_{k,*}} \end{bmatrix}$ is an orthonormal matrix. Then $\begin{bmatrix} \Lambda_{k,*} & H_{k,*}     
\end{bmatrix}\in\Real^{\sfq_\ell\times \sfq_\ell}$ is invertible since $\begin{bmatrix} \Lambda_{k,*} & H_{k,*} \end{bmatrix}=\Diag(\sqrt{h^{k,*}})\begin{bmatrix} \sfV_{k,*} &  \sqrt{H_{k,*}} \end{bmatrix}$ is the product of two invertible matrices. Therefore
\begin{align}\label{eq:LMI1_stability}
    \lambda_k \bfP_k - M_{k,*}^\top \bfP_k M_{k,*} \succ 0 \iff
    & \begin{bmatrix} \Lambda_{k,*}^\top \\ H_{k,*}^\top \end{bmatrix} \left(\lambda_k \bfP_k - M_{k,*}^\top \bfP_k M_{k,*} \right) \begin{bmatrix} \Lambda_{k,*} & H_{k,*} \end{bmatrix}:=\begin{bmatrix} A_{11} & A_{12} \\ A_{12}^\top & A_{22}  \end{bmatrix} \succ 0
\end{align}
where for $\bfP_k=C^\bfP_k \Diag(1/h^{k,*})$ with $C^\bfP_k>0$, following holds using Lemma \ref{lemma:M_k,*}
\begin{align*}
    & A_{11} := \lambda_k \Lambda_{k,*}^\top \bfP_k\Lambda_{k,*} - \Lambda_{k,*}^\top M_{k,*}^\top \bfP_k M_{k,*} \Lambda_{k,*} = \lambda_k C^\bfP_k I - C^\bfP_k \widehat M_{k,*}^2, 
    \\
    & A_{12} := \lambda_k \Lambda_{k,*}^\top \bfP_k H_{k,*} - \Lambda_{k,*}^\top M_{k,*}^\top \bfP_k M_{k,*} H_{k,*} = \alpha_k C^\bfP_k  \widehat{M}_{*} [\Lambda_{k,*}^\top \nabla^2_{h} g(h^{k,*},\xUL^k) H_{k,*}],
    \\
    & A_{22} 
    \begin{aligned}[t]
    & := \lambda_k H_{k,*}^\top \bfP_k H_{k,*} - H_{k,*}^\top M_{k,*}^\top \bfP_k M_{k,*} H_{k,*} \\
    & = \lambda_k  C^\bfP_k I- \alpha_k^2 C^\bfP_k\ [H_{k,*}^\top \nabla^2_{h} g(h^{k,*},\xUL^k) \Lambda_{k,*}] [\Lambda_{k,*}^\top \nabla^2_{h} g(h^{k,*},\xUL^k) H_{k,*}].
    \end{aligned}
\end{align*} 
Using Schur complement 
\begin{align}\label{eq:LMI1-v2_stability}
    \eqref{eq:LMI1_stability} \iff A_{11}\succ 0,\ A_{22} - A_{12}^\top A_{11}^{-1} A_{12} \succ 0
\end{align}
Choose $\lambda_k={1-\eta\alpha_k/2}$ which gives,
\begin{align*}
    & A_{11} \succeq \lambda_k C^\bfP_k I - C^\bfP_k (1-\eta\alpha_k)^2 I = C^\bfP_k {(1-\eta\alpha_k/2)}  I - C^\bfP_k (1-\eta\alpha_k)^2 I \succeq {\frac{\eta\alpha_k}{2}}(1-\eta\alpha_k)C^\bfP_k I \succ 0 \\ \\
    & A_{22} - A_{12}^\top A_{11}^{-1} A_{12}\\
    & = \begin{aligned}[t]& \lambda_k  C^\bfP_k I- \alpha_k^2 C^\bfP_k\ [H_{k,*}^\top \nabla^2_{h} g(h^{k,*},\xUL^k) \Lambda_{k,*}] [\Lambda_{k,*}^\top \nabla^2_{h} g(h^{k,*},\xUL^k) H_{k,*}]\\
    & \qquad - \alpha_k^2  (C^\bfP_k)^2 [H_{k,*}^\top \nabla^2_{h} g(h^{k,*},\xUL^k) \Lambda_{k,*}]  \widehat{M}_{*} A_{11}^{-1} \widehat{M}_{*} [\Lambda_{k,*}^\top \nabla^2_{h} g(h^{k,*},\xUL^k) H_{k,*}] \end{aligned}\\
    & = \lambda_k  C^\bfP_k I- \alpha_k^2 C^\bfP_k\ [H_{k,*}^\top \nabla^2_{h} g(h^{k,*},\xUL^k) \Lambda_{k,*}] \left(I + C^\bfP_k \widehat{M}_{*} A_{11}^{-1} \widehat{M}_{*}\right) [\Lambda_{k,*}^\top \nabla^2_{h} g(h^{k,*},\xUL^k) H_{k,*}] \\
    & \succeq \lambda_k  C^\bfP_k I- \alpha_k^2 C^\bfP_k\ [H_{k,*}^\top \nabla^2_{h} g(h^{k,*},\xUL^k) \Lambda_{k,*}] \left(1 + C^\bfP_k \frac{(1-\eta\alpha_k)^2}{C^\bfP_k {\frac{\eta\alpha_k}{2}} (1-\eta\alpha_k)}\right) [\Lambda_{k,*}^\top \nabla^2_{h} g(h^{k,*},\xUL^k) H_{k,*}] \\
    & = C^\bfP_k \left({(1-\eta\alpha_k/2)} I- \frac{\alpha_k^2{(1-\eta\alpha_k/2)}}{{\frac{\eta\alpha_k}{2}}} [H_{k,*}^\top \nabla^2_{h} g(h^{k,*},\xUL^k) \Lambda_{k,*}] [\Lambda_{k,*}^\top \nabla^2_{h} g(h^{k,*},\xUL^k) H_{k,*}]\right)
\end{align*}
Using Lemma \ref{lemma:M_k,*}, $H_{k,*}^\top \nabla^2_{h} g(h^{k,*},\xUL^k) \Lambda_{k,*}\Lambda_{k,*}^\top \nabla^2_{h} g(h^{k,*},\xUL^k) H_{k,*} \preceq H_{k,*}^\top (\nabla^2_{h} g(h^{k,*},\xUL^k))^2 H_{k,*} \preceq L_g^2 I$, then
\begin{align} \label{eq:LMI1-v2-b_stability}
    & A_{22} - A_{12}^\top A_{11}^{-1} A_{12} \succeq C^\bfP_k {\left(1-\frac{\eta\alpha_k}{2}\right)} \left(1- \frac{2\alpha_k}{\eta} L_g^2\right) I
\end{align}
From \eqref{eq:LMI1-v2-b_stability}, we have $A_{22} - A_{12}^\top A_{11}^{-1} A_{12}\succ0$ for $\alpha_k<\frac{\eta}{2L_g^2}$. 
Fix $\alpha_k\leq\frac{\eta}{2(L_g^2+\eta^2)}$, then \eqref{eq:LMI1-v2_stability} or equivalently \eqref{eq:LMI1_stability}  holds as  
\begin{subequations}\label{eq:LMI1-v3_stability}    
\begin{align}
    & A_{11} \succeq C^\bfP_k {\frac{\eta\alpha_k}{2}} (1-\eta\alpha_k) I \succeq C^\bfP_k {\frac{\eta\alpha_k}{2}} \left(1-\frac{\eta^2}{2(L_g^2+\eta^2)}\right) I \succeq {\frac{\eta\alpha_k}{4}} C^\bfP_k I \\
    & \begin{aligned}[t] A_{22} - A_{12}^\top A_{11}^{-1} A_{12}
    & \succeq C^\bfP_k {\left(1-\frac{\eta\alpha_k}{2}\right)} \left(1- \frac{2\alpha_k}{\eta} L_g^2\right) I\\
    & \succeq C^\bfP_k {\left(1-\frac{\eta^2}{4(L_g^2+\eta^2)}\right)} \left(1- \frac{2}{\eta}\frac{\eta}{4L_g^2} L_g^2\right) I \succeq \frac{3}{8} C^\bfP_k I 
    \end{aligned}
\end{align}
\end{subequations}
where $C^\bfP_k>0$. Since $0<\eta\alpha_k<1$ then ${3}/{8}>{{\eta\alpha_k}/{4}}$. \\
\newline \underline{Lower bounding $\lambda_k \bfP_k - M_{k,*}^\top \bfP_k M_{k,*} \succ 0$ for $\bfP_k=C^\bfP_k \Diag(1/h^{k,*}) \succ 0$ and $0<\lambda_k={1-\eta\alpha_k/2}<1$}:
\begin{align*}
    \begin{bmatrix} A_{11} & A_{12} \\ A_{12}^\top & A_{22}  \end{bmatrix}
    & = \begin{bmatrix}
        I & 0 \\ A_{12}^\top A_{11}^{-1} & I
    \end{bmatrix}
    \begin{bmatrix}
        A_{11} & 0 \\ 0 & A_{22} - A_{12}^\top A_{11}^{-1} A_{12}
    \end{bmatrix}
    \begin{bmatrix}
        I & A_{11}^{-1} A_{12} \\ 0 & I
    \end{bmatrix} \overset{\eqref{eq:LMI1-v3_stability}}{\succeq} {\frac{\eta\alpha_k}{4}} C^\bfP_k \begin{bmatrix}
        I & 0 \\ A_{12}^\top A_{11}^{-1} & I
    \end{bmatrix}
    \begin{bmatrix}
        I & A_{11}^{-1} A_{12} \\ 0 & I
    \end{bmatrix}
\end{align*}
Since $\begin{bmatrix} I & A_{11}^{-1} A_{12} \\ 0 & I \end{bmatrix}^{-1} = \begin{bmatrix} I & & \ -A_{11}^{-1} A_{12} \\ 0 & & \ I
    \end{bmatrix}$ \text{and} 
    $\left\|\begin{bmatrix}
        I & & \ -A_{11}^{-1} A_{12} \\ 0 & & \ I
    \end{bmatrix}\right\|\leq 1 + \left\|A_{11}^{-1} A_{12} \right\|\leq 1 + \|A_{11}^{-1}\|\cdot\|A_{12}\|$ therefore 
    $\begin{bmatrix}
        I & 0 \\ A_{12}^\top A_{11}^{-1} & I \end{bmatrix} \begin{bmatrix}
    I & A_{11}^{-1} A_{12} \\ 0 & I \end{bmatrix} \succeq \frac{1}{(1 + \|A_{11}^{-1}\|\cdot\|A_{12}\|)^2} I \succeq \frac{1}{1 + (\|A_{11}^{-1}\|\cdot\|A_{12}\|)^2} I$. 
    From \eqref{eq:LMI1-v3_stability}, $\|A_{11}^{-1}\| \leq {{4}/(\eta\alpha_k C^\bfP_k)}$. Using $\|A_{12}\|=\sqrt{\|A_{12}^\top A_{12}\|}$ below to bound $\|A_{12}\|$.
    $$\begin{aligned}
A_{12}^\top A_{12} & =\alpha_k^2  (C^\bfP_k)^2 [H_{k,*}^\top \nabla^2_{h} g(h^{k,*},\xUL^k) \Lambda_{k,*}]  \widehat{M}_{*} \widehat{M}_{*} [\Lambda_{k,*}^\top \nabla^2_{h} g(h^{k,*},\xUL^k) H_{k,*}] \\
& \preceq \alpha_k^2  (C^\bfP_k)^2 (1-\eta\alpha_k)^2 [H_{k,*}^\top \nabla^2_{h} g(h^{k,*},\xUL^k) \Lambda_{k,*}] [\Lambda_{k,*}^\top \nabla^2_{h} g(h^{k,*},\xUL^k) H_{k,*}] \\
& \preceq \alpha_k^2  (C^\bfP_k)^2 (1-\eta\alpha_k)^2 L_g^2 I \\
\implies \|A_{12}\| & \leq \alpha_k (1-\eta\alpha_k) C^\bfP_k L_g
\end{aligned}$$
As a result, $\begin{bmatrix} A_{11} & A_{12} \\ A_{12}^\top & A_{22}  \end{bmatrix} 
\succeq \widetilde \mu_k C^\bfP_k I$ where  
$\widetilde \mu_k = \frac{{\eta\alpha_k/4}}{1+{16} (1-\eta\alpha_k)^2 L_g^2/{\eta^2}}>0$. Recall from \eqref{eq:LMI1_stability}
\[
\begin{bmatrix} \Lambda_{k,*}^\top \\ H_{k,*}^\top \end{bmatrix} \left(\lambda_k \bfP_k - M_{k,*}^\top \bfP_k M_{k,*} \right) \begin{bmatrix} \Lambda_{k,*} & H_{k,*} \end{bmatrix}=\begin{bmatrix} A_{11} & A_{12} \\ A_{12}^\top & A_{22}  \end{bmatrix}
\]
Using $\begin{bmatrix} \Lambda_{k,*} & H_{k,*} \end{bmatrix}=\Diag(\sqrt{h^{k,*}})\begin{bmatrix} \sfV_{k,*} &  \sqrt{H_{k,*}} \end{bmatrix}$, the fact $\begin{bmatrix} \sfV_{k,*} &  \sqrt{H_{k,*}} \end{bmatrix}$ is an orthonormal matrix, and $h^{k,*}_i$ lies in a simplex, we have
\begin{align}\label{eq:LMI2_stability}
    \lambda_k\ \bfP_k - M_{k,*}^\top \bfP_k M_{k,*}  & = 
    \Diag\Big(\frac{1}{\sqrt{h^{k,*}}}\Big) \begin{bmatrix} \sfV_{k,*} & \sqrt{H_{k,*}} \end{bmatrix} \begin{bmatrix} A_{11} & A_{12} \\ A_{12}^\top & A_{22}  \end{bmatrix} \begin{bmatrix} {\sfV_{k,*}}^\top \\ \sqrt{H_{k,*}}^\top \end{bmatrix} \Diag\Big(\frac{1}{\sqrt{h^{k,*}}}\Big) \succeq \widetilde \mu_k\ C^\bfP_k I
\end{align}
\underline{\smash{Analyzing robust $\mathcal{Q}$-stability requirement \eqref{eq:LMI_Qstability-a}}}:
\newline To prove $\exists$ $C^\bfP_k, \epsilon_k, s_k>0$ such that for $\bfP_k=C^\bfP_k \Diag(1/h^{k,*})$ and $\lambda_k={1-\eta\alpha_k/2}$ following holds
\begin{align*}
    \mathcal{M}_k 
    = \begin{bmatrix} \lambda_k \bfP_k - M_{k,*}^\top \bfP_k M_{k,*} - (s_k+1) I  & & - \epsilon_k M_{k,*}^\top \bfP_k \\ -\epsilon_k \bfP_k M_{k,*} & & s_k I-\epsilon_k^2 \bfP_k \end{bmatrix} \succ 0 
\end{align*}
Choose  $\epsilon_k$ such that $\epsilon_k^2 \bfP_k\preceq \frac{s_k}{2} I$ and choose $C^\bfP_k=(\frac{3}{2}s_k+1)/{\widetilde \mu_k}$
\begin{align*}
    \mathcal{M}_k \overset{\eqref{eq:LMI2_stability}}{\succeq} \begin{bmatrix} \widetilde \mu_k\ C^\bfP_k I - (s_k+1) I  & & - \epsilon_k M_{k,*}^\top \bfP_k \\ -\epsilon_k \bfP_k M_{k,*} & & s I-\epsilon_k^2 \bfP_k \end{bmatrix} 
    \succeq \begin{bmatrix} \frac{s_k}{2} I  & & - \epsilon_k M_{k,*}^\top \bfP_k \\ 
    -\epsilon_k \bfP_k M_{k,*} & & \frac{s_k}{2} I \end{bmatrix}
\end{align*}
Say $\epsilon_k$ also satisfies $\begin{bmatrix} {s_k}/{4}\ I  & & 0 \\ 0 & & {s_k}/{4}\ I \end{bmatrix} \succeq \begin{bmatrix} 0  & & \epsilon_k M_{k,*}^\top \bfP_k \\ \epsilon_k \bfP_k M_{k,*} & & 0 \end{bmatrix}$ which implies $\mathcal{M} \succeq \begin{bmatrix} \frac{s_k}{4} I  & & 0 \\ 0 & & \frac{s_k}{4} I \end{bmatrix}$.
Since $\left\|\begin{bmatrix} 0  & & \epsilon_k M_{k,*}^\top \bfP_k \\ 
\epsilon_k \bfP_k M_{k,*} & & 0 \end{bmatrix} \right\| \leq 2\epsilon_k\|M_{k,*}^\top \bfP_k\|$ then it is sufficient that $\epsilon_k$ satisfies: $\epsilon_k^2 \bfP_k\preceq \frac{s_k}{2} I$,\ $2\epsilon_k\|M_{k,*}^\top \bfP_k\| \leq \frac{s_k}{4}$.\quad \\
\newline We further bound $\|M_{k,*}^\top \bfP_k\|$ as
\begin{align*}
     & M_{k,*}^\top \bfP_k \bfP_k M_{k,*} \\
     & \begin{aligned}[t] = (C^\bfP_k)^2 \big[(1-\eta\alpha_k)\Diag\left(1/h^{k,*}\right)-\alpha_k \nabla^2_{h} g(h^{k,*},\xUL^k)\big]\Lambda_{k,*}\Lambda_{k,*}^\top \Diag^2\left(1/h^{k,*}\right) \Lambda_{k,*}\Lambda_{k,*}^\top \\ \big[(1-\eta\alpha_k)\Diag\left(1/h^{k,*}\right)-\alpha_k \nabla^2_{h} g(h^{k,*},\xUL^k)\big] \end{aligned}\\
     & \preceq (C^\bfP_k)^2/\nu^{\min} (1-\eta\alpha_k)^2/(\nu^{\min})^2 \quad \text{(using Lemma \ref{lemma:M_k,*})}\\
     & \implies \|M_{k,*}^\top \bfP_k\| \leq C^\bfP_k (1-\eta\alpha_k)/(\nu^{\min})^{3/2}
\end{align*}
Choose $s_k=4$ so that $\mathcal{M}\succeq I$ and since $\bfP_k\preceq C^\bfP_k/\nu^{\min}\ I$ it suffices to choose $\epsilon_k$ satisfying
\begin{align*}
    \epsilon_k^2  C^\bfP_k/\nu^{\min}\leq \frac{4}{2},\quad 2\epsilon_k\frac{C^\bfP_k (1-\eta\alpha_k)}{(\nu^{\min})^{3/2}} \leq \frac{4}{4}
\end{align*}
As a result, choose $\epsilon_k \leq \overline \epsilon_k = \min\left\{\sqrt{\frac{2 \nu^{\min}}{C^\bfP_k}},\frac{(\nu^{\min})^{3/2}}{2 C^\bfP_k (1-\eta\alpha_k)}\right\}$. \\
\newline \underline{\smash{Analyzing robust $\mathcal{Q}$-stability requirement \eqref{eq:LMI_Qstability-b}}}: 
\newline To prove for $\bfP_k=C^\bfP_k \Diag(1/h^{k,*})$, $\lambda_k={1-\eta\alpha_k/2}$, $0<\epsilon_k\leq \overline \epsilon_k$, and $s_k=4$ there exists $\omega_k>0$ such that following holds
\begin{align*}
    & \omega_k I-\bfP_k - \bfP_k \begin{bmatrix} M_{k,*} & \epsilon_k I\end{bmatrix} \mathcal{M}^{-1} \begin{bmatrix} M_{k,*}^\top  \\ \epsilon_k I \end{bmatrix} \bfP_k \succeq 0
\end{align*}
Now,
\begin{align*}
    & \omega_k I-\bfP_k - \bfP_k \begin{bmatrix} M_{k,*} \\ & \epsilon_k I\end{bmatrix} \mathcal{M}^{-1} \begin{bmatrix} M_{k,*}^\top  \\ \epsilon_k I \end{bmatrix} \bfP_k \\
    & \succeq \omega_k I-\bfP_k - \bfP_k \begin{bmatrix} M_{k,*} & \epsilon_k I\end{bmatrix} \begin{bmatrix} M_{k,*}^\top  \\ \epsilon_k I \end{bmatrix} \bfP_k \quad \text{(using $\mathcal{M}\succeq I$)}\\
    & = \omega_k I-\bfP_k - \bfP_k M_{k,*}M_{k,*}^\top \bfP_k - \epsilon_k^2 \bfP_k^2 \\
    & \succeq \omega_k I-\bfP_k - \frac{(C^\bfP_k)^2 (1-\eta\alpha_k)^2}{(\nu^{\min})^3} I - 2 \bfP_k \\
    & \quad \text{\Big(using $\|\bfP_k M_{k,*}M_{k,*}^\top \bfP_k\|=\|M_{k,*}^\top \bfP_k \bfP_k M_{k,*}\|\leq (C^\bfP_k)(1-\eta\alpha_k)^2/(\nu^{\min})^3$ and $\epsilon_k^2 \bfP_k\preceq \frac{s_k}{2} I= 2I$\Big)} \\
    & \succeq \omega_k I- \frac{3 C^\bfP_k}{\nu^{\min}} I - \frac{(C^\bfP_k)^2 (1-\eta\alpha_k)^2}{(\nu^{\min})^3} I \quad \text{(using $\bfP_k\preceq C^\bfP_k/\nu^{\min}\ I$)}
\end{align*}
It suffices to choose $\omega_k={3 C^\bfP_k}/{\nu^{\min}}+{(C^\bfP_k)^2 (1-\eta\alpha_k)^2}/{(\nu^{\min})^3}$. The proof is complete.
\Halmos
\endproof

\subsection{Proof of Lemma \ref{lemma:ConvgRate-W_k,t}.}\label{appendix:Proof-lemma:ConvgRate-W_k,t}
Note that Lemma \ref{lemma:Bnd-R*_Lip-h*}\ref{lemma:Bounded_R*} implies 
\[\|\rmW_{k,t}\| \leq \big\|U_{k,*}-U_{k,t}\big\| + \sfC_0 \|\Delta_{k,t}\|\]
\begin{enumerate}
\item Bounding $\big\|U_{k,*}-U_{k,t}\big\|$.
\begin{align*}
    & \big\|U_{k,*}-U_{k,t}\big\| \\ & = \alpha_k \big\|B_{k,*} \nabla_\xUL \nabla_{h} g(h^{k,*},\xUL^k)-B_{k,t} \nabla_\xUL \nabla_{h} g(h^{k,t},\xUL^k)\big\| \\
                        & \leq \alpha_k \|B_{k,*}\| \big\|\nabla_\xUL \nabla_{h} g(h^{k,*},\xUL^k)-\nabla_\xUL \nabla_{h} g(h^{k,t},\xUL^k)\big\|+ \big\|B_{k,*}-B_{k,t}\big\| \big\|\nabla_\xUL \nabla_{h} g(h^{k,t},\xUL^k)\big\| \\
                        & \leq \alpha_k\ \tau_g^h\ \|h^{k,*}-h^{k,t}\|+\alpha_k\ \Lambda^{\xUL h}_g L_B\ \|h^{k,*}-h^{k,t+1}\|\\
                        & \qquad \text{(using Lemma \ref{lemma:Lipschtz_matrixB}, $\|B_{k,*}\|\leq1$ from Lemma \ref{lemma:B_k,t factor})} \\
                        & \leq C^U_k\ ({\txteps_h^{k,0}})^{1/2}\ (1- \eta \alpha_k)^{t/2} \quad \text{(using Lemma \ref{lemma:ConvgRate-Soln})}
\end{align*}
where $C^U_k = \sqrt{2}\ \alpha_k\ (\tau_g^h+\Lambda^{\xUL h}_g L_B (1- \eta \alpha_k)^{1/2})$. 

\item Bounding $\|\Delta_{k,t}\|$.
\newline Recall $\Delta_{k,t} = M_{k,t}-M_{k,*} = B_{k,t}[(1-\eta\alpha_k)\Diag\left(1/h^{k,t}\right)-\alpha_k \nabla^2_{h} g(h^{k,t},\xUL^k)\big]-B_{k,*}\big[(1-\eta\alpha_k)\Diag\left(1/h^{k,*}\right)-\alpha_k \nabla^2_{h} g^\eta(h^{k,*},\xUL^k)]$.
\begin{align*}
    \|\Delta_{k,t}\|
    & \leq \|B_{k,t}\|\cdot[(1-\eta\alpha_k)\|\Diag\left(1/h^{k,t}\right)-\Diag\left(1/h^{k,*}\right)\|+\alpha_k \|\nabla^2_h g(h^{k,t},\xUL^k) - \nabla^2_h g(h^{k,*},\xUL^k)\|] \\ 
    & \qquad + \|B_{k,t}-B_{k,*}\| \|(1-\eta\alpha_k)\Diag\left(1/h^{k,*}\right)-\alpha_k \nabla^2_{h} g(h^{k,*},\xUL^k)\|
\end{align*}
{For $0<\alpha_k\leq\overline \alpha$, it holds $0\preceq(1-(L_g+\eta)\alpha_k) I \preceq (1-\eta\alpha_k)\Diag\left(1/h^{k,*}\right)-\alpha_k \nabla^2_{h} g(h^{k,*},\xUL^k)\preceq (1-\eta\alpha_k)/\nu^{\min} I \implies \|(1-\eta\alpha_k)\Diag\left(1/h^{k,*}\right)-\alpha_k \nabla^2_{h} g(h^{k,*},\xUL^k)\|\leq (1-\eta\alpha_k)/\nu^{\min}$.}
\begin{align*}
\|\Delta_{k,t}\| & \leq \begin{aligned}[t] \left(L_{1/h} (1-\eta\alpha_k) + \rho_g^h \alpha_k\right) \|h^{k,*}-h^{k,t}\| + L_B (1-\eta\alpha_k)/\nu^{\min} \|h^{k,*}-h^{k,t+1}\| \\
    \text{(using $\|B_{k,t}\|\leq1$ from Lemma \ref{lemma:B_k,t factor}, Lemma \ref{lemma:Lipschtz_matrixB}, Lemma \ref{lemma:Lipschtz_1/h})}
    \end{aligned}\\
    & \leq C^\Delta_k\ ({\txteps_h^{k,0}})^{1/2}\ (1- \eta \alpha_k)^{t/2} \quad \text{(using Lemma \ref{lemma:ConvgRate-Soln})}
\end{align*}
where $C^\Delta_k = \sqrt{2}\ \left(L_{1/h} (1-\eta\alpha_k) + \rho_g^h \alpha_k + L_B (1-\eta\alpha_k)^{3/2}/\nu^{\min}\right)$.
\end{enumerate}
Therefore, $\begin{aligned}[t]
    \|\rmW_{k,t}\| & \leq \big\|U_{k,*}-U_{k,t}\big\| + \sfC_0 \|\Delta_{k,t}\| \leq (C^U_k + \sfC_0 C^\Delta_k)\ ({\txteps_h^{k,0}})^{1/2}\ (1- \eta \alpha_k)^{t/2}
\end{aligned}$. \Halmos
\endproof

\subsection{Proof of Lemma \ref{lemma:Stab_dyn}.}\label{appendix:Proof-lemma:Stab_dyn}
Using $\|\Delta_{k,t}\|\leq C^\Delta_k\ ({\txteps_h^{k,0}})^{1/2}\ (1- \eta \alpha_k)^{t/2}$ from Lemma \ref{lemma:ConvgRate-W_k,t} and $\txteps_h^{k,0}\leq \overline \sfD_\KL^{\max}$ from Lemma \ref{lemma:GlobBndErr}, $\exists\ T^0_k\geq0$ 
such that $\|\Delta_{k,t}\|\leq \epsilon_k\;\;\forall t\geq T^0_k$. Doing
\begin{align*}
 & C^\Delta_k\ (\sfD_\KL^{\max})^{1/2} (1- \eta \alpha_k)^{t/2} \leq \epsilon_k \;\; \iff \;\; t \geq \frac{\ln(C^\Delta_k (\sfD_\KL^{\max})^{1/2}/\epsilon_k)}{\ln(1/\sqrt{1- \eta \alpha_k})}
\end{align*}
and noting that $\ln(1/\sqrt{1- \eta \alpha_k})=-\frac{1}{2}\ln(1- \eta \alpha_k)\geq \frac{1}{2} \eta \alpha_k$, it suffices to take $T^0_k = \left\lceil \frac{2\ln(C^\Delta_k (\sfD_\KL^{\max})^{1/2}/\epsilon_k)}{\eta \alpha_k} \right\rceil$.

\begin{enumerate}[label=\normalfont(\alph*)]
\item From Lemma \ref{lemma:Robust_Q-stability}, multiply LMI \eqref{eq:LMI_eps} with $\begin{bmatrix} \rmX_{k,t} & \Omega_{k,t} & \rmW_{k,t} \end{bmatrix}^\top$ from left and with $\begin{bmatrix} \rmX_{k,t} \\ \Omega_{k,t} \\ \rmW_{k,t} \end{bmatrix}$ from right to get \eqref{eq:proof2-ineq1} which holds for $\epsilon_k \leq \overline \epsilon_k$.
\begin{align}\label{eq:proof2-ineq1}
 \rmX_{k,t+1}^\top \bfP_k \rmX_{k,t+1} + s_k \Psi_{k,t}^\top \Psi_{k,t} + \rmY_{k,t}^\top \rmY_{k,t} \ \preceq \ \lambda_k\ \rmX_{k,t}^\top \bfP_k \rmX_{k,t} + s_k \Omega_{k,t}^\top \Omega_{k,t} + \omega_k\rmW_{k,t}^\top \rmW_{k,t} 
\end{align}

Substitute $\Omega_{k,t} = \frac{1}{\epsilon_k} \Delta_{k,t} \Psi_{k,t}$ from \eqref{eq:dynamic_sys-sep} and using the fact $\rmY_{k,t}^\top \rmY_{k,t}\succeq 0$,
\begin{align}\label{eq:proof2-ineq2}
 \rmX_{k,t+1}^\top \bfP_k \rmX_{k,t+1} + s_k \Psi_{k,t}^\top (I-\frac{1}{\epsilon_k^2}\Delta_{k,t}^\top \Delta_{k,t}) \Psi_{k,t} \ \preceq \ \lambda_k\ \rmX_{k,t}^\top \bfP_k \rmX_{k,t} + \omega_k\ \rmW_{k,t}^\top \rmW_{k,t} 
\end{align}

For $0 < \epsilon_k \leq C^\Delta_k (\sfD_\KL^{\max})^{1/2}$ and $t\geq T^0_k$, it holds $\Delta_{k,t}^\top \Delta_{k,t}\preceq \epsilon_k^2\ I$ which implies $\Psi_{k,t}^\top (I-\frac{1}{\epsilon_k^2}\Delta_{k,t}^\top \Delta_{k,t}) \Psi_{k,t}\succeq0$ so that the term can be dropped from \eqref{eq:proof2-ineq2}. \newline Then for $\epsilon_k = \min\left\{C^\Delta_k (\sfD_\KL^{\max})^{1/2},\ \overline \epsilon_k\right\}$, following holds
\begin{align}\label{eq:proof2-ineq3}
 & \rmX_{k,t+1}^\top \bfP_k \rmX_{k,t+1} \preceq \lambda_k\ \rmX_{k,t}^\top \bfP_k \rmX_{k,t} + \omega_k\ \rmW_{k,t}^\top \rmW_{k,t},\quad t\geq T^0_k 
\end{align}
Doing recursion on \eqref{eq:proof2-ineq3}, we have for $t\geq T^0_k$ 
\begin{align*}
\rmX_{k,t}^\top \bfP_k \rmX_{k,t} \preceq (\lambda_k)^{t-T^0_k}\ \rmX_{k,T^0_k}^\top \bfP_k \rmX_{k,T^0_k} + \omega_k\ \sum_{\ell=T^0_k}^{t-1} (\lambda_k)^{\ell}\ \rmW_{k,t-1-\ell}^\top \rmW_{k,t-1-\ell}
\end{align*}
Use the fact, $C^{\bfP}_k \ I \ \preceq \ \bfP_k = C^{\bfP}_k\ \Diag(1/h^{k,*}) \ \preceq \ C^{\bfP}_k/\nu^{\min}\ I$ \ to get
\begin{align*}
& C^{\bfP}_k \ \rmX_{k,t}^\top \rmX_{k,t} \preceq C^{\bfP}_k/\nu^{\min}\ (\lambda_k)^{t-T^0_k}\ \rmX_{k,T^0_k}^\top \rmX_{k,T^0_k} + \omega_k\ \sum_{\ell=T^0_k}^{t-1} (\lambda_k)^{\ell}\ \rmW_{k,t-1-\ell}^\top \rmW_{k,t-1-\ell} \\
\implies & \|\rmX_{k,t}\|^2 \leq \frac{(\lambda_k)^{t-T^0_k}}{\nu^{\min}}\ \|\rmX_{k,T^0_k}\|^2 + \omega_k/C^{\bfP}_k\ \sum_{\ell=T^0_k}^{t-1} (\lambda_k)^{\ell}\ \|\rmW_{t-1-\ell}\|^2,\quad t\geq T^0_k
\end{align*}
Using $\|\rmW_{k,t}\| \leq (C^U_k + \sfC_0 C^\Delta_k)\ ({\txteps_h^{k,0}})^{1/2}\ (1- \eta_k \alpha_k)^{t/2}$ from Lemma \ref{lemma:ConvgRate-W_k,t}
\begin{align*}
\|\rmX_{k,t}\|^2 \leq \frac{(\lambda_k)^{t-T^0_k}}{\nu^{\min}}\ \|\rmX_{k,T^0_k}\|^2 + (C^U_k + \sfC_0 C^\Delta_k)^2\ \omega_k/C^{\bfP}_k\ {\txteps_h^{k,0}}\ \sum_{\ell=T^0_k}^{t-1} (\lambda_k)^{\ell}\ (1- \eta_k \alpha_k)^{t-1-\ell},\quad t\geq T^0_k
\end{align*}
Define $\phi_k={(1- \eta \alpha_k)}/{\lambda_k}$ and simplify $\sum_{\ell=T^0_k}^{t-1} (\lambda_k)^{{\ell}}\ (1- \eta_k \alpha_k)^{t-1-\ell}=(\lambda_k)^{t-1}\sum_{\ell=T^0_k}^{t-1} \phi_k^{t-1-\ell}$ 
to complete the proof.
\item Starting from \eqref{eq:proof2-ineq2}
\begin{align*}
 \rmX_{k,t+1}^\top \bfP_k \rmX_{k,t+1} + s_k\ \Psi_{k,t}^\top (I-\frac{1}{\epsilon_k^2}\Delta_{k,t}^\top \Delta_{k,t}) \Psi_{k,t} \ \preceq \ \lambda_k\ \rmX_{k,t}^\top \bfP_k \rmX_{k,t} + \omega_k\ \rmW_{k,t}^\top \rmW_{k,t}
\end{align*}
Note if $\epsilon_k = \min\left\{C^\Delta_k (\sfD_\KL^{\max})^{1/2},\ \overline \epsilon_k\right\}$, then $\|\Delta_{k,t}\|$ may be larger than $\epsilon_k$ for $0\leq t < T^0_k$. As a result the term $\Psi_{k,t}^\top (I-\frac{1}{\epsilon_k^2}\Delta_{k,t}^\top \Delta_{k,t}) \Psi_{k,t}$ cannot be dropped from the analysis. 

Using $\Psi_{k,t}=\rmX_{k,t}$ from \eqref{eq:dynamic_sys-sep}, then more generally for $0\leq t < T^0_k$
\begin{align*}
 \rmX_{k,t+1}^\top \bfP_k \rmX_{k,t+1} & \ \preceq \ \lambda_k\ \rmX_{k,t}^\top \bfP_k \rmX_{k,t} + \frac{s_k}{\epsilon_k^2}\ \rmX_{k,t}^\top \left(\Delta_{k,t}^\top \Delta_{k,t}-\epsilon_k^2\ I\right) \rmX_{k,t} + \omega_k\ \rmW_{k,t}^\top \rmW_{k,t}\\
 & \ \preceq \ \lambda_k\ \rmX_{k,t}^\top \bfP_k \rmX_{k,t} + \frac{4}{\epsilon_k^2}\ (\|\Delta_{k,t}\|^2-\epsilon_k^2)\ \rmX_{k,t}^\top \rmX_{k,t} + \omega_k\ \rmW_{k,t}^\top \rmW_{k,t} \nonumber\\
 & \ \preceq \begin{aligned}[t]\left(\lambda_k\ + \frac{(C^\Delta_k)^2 \sfD_\KL^{\max}-\epsilon_k^2}{\epsilon_k^2/4\ C^{\bfP}_k}\right)\ \rmX_{k,t}^\top \bfP_k \rmX_{k,t} + \omega_k\ \rmW_{k,t}^\top \rmW_{k,t} \\ \big(\because\ \|\Delta_{k,t}\|\leq C^\Delta_k (\sfD_\KL^{\max})^{1/2},\ C^{\bfP}_k \ I \preceq \bfP_k\big)
 \end{aligned}
\end{align*}
 Define $\tlambda_k=\lambda_k\ + \frac{(C^\Delta_k)^2 \sfD_\KL^{\max}-\epsilon_k^2}{\epsilon_k^2/4\ C^{\bfP}_k}$ then after doing recursion, we have for $0\leq t\leq T^0_k$
\begin{align*}
\rmX_{k,t}^\top \bfP_k \rmX_{k,t} \preceq (\tlambda_k)^{t}\ \rmX_{k,0}^\top \bfP_k \rmX_{k,0} + \omega_k\ \sum_{\ell=0}^{t-1} (\tlambda_k)^{\ell}\ \rmW_{k,t-1-\ell}^\top \rmW_{k,t-1-\ell}
\end{align*}
Using \ $C^{\bfP}_k \ I \preceq \ \bfP_k = C^{\bfP}_k\ \Diag(1/h^{k,*}) \preceq \ C^{\bfP}_k/\nu^{\min}\ I$ \ 
\begin{align*}
\|\rmX_{k,t}\|^2 & \leq \frac{(\tlambda_k)^{\frac{t}{2}}}{\nu^{\min}}\ \|\rmX_{k,0}\|^2 + \omega_k/C^{\bfP}_k\ \sum_{\ell=0}^{t-1} (\tlambda_k)^{\ell}\ \|\rmW_{k,t-1-\ell}\|^2,\quad 0\leq t \leq T^0_k 
\end{align*}
Firstly, use $\rmX_{k,0}=-\sfR_{k,*}$ to bound $\|\rmX_{k,0}\|=\|\sfR_{k,*}\|\leq \sfC_0$ from Lemma \ref{lemma:Bnd-R*_Lip-h*}\ref{lemma:Bounded_R*}. Secondly, use $\|\rmW_{k,t}\| \leq (C^U_k + \sfC_0 C^\Delta_k) ({\txteps_h^{k,0}})^{1/2} (1- \eta_k \alpha_k)^{t/2}$ from Lemma \ref{lemma:ConvgRate-W_k,t} and define $\widetilde\phi_k={(1- \eta \alpha_k)}/{\widetilde\lambda_k}$ to simplify $\sum_{\ell=T^0_k}^{t-1} (\widetilde\lambda_k)^{{\ell}}\ (1- \eta_k \alpha_k)^{t-1-\ell}=(\widetilde\lambda_k)^{t-1}\sum_{\ell=T^0_k}^{t-1} {\widetilde\phi_k}^{t-1-\ell}$ which completes the proof.
\Halmos
\end{enumerate}
\endproof

\section{Smoothness of $F(\xUL)$}\label{appendix:Lip_F}

\proof{Proof of Lemma \ref{lemma:Lipschtz_BilevelGrad}.}
\begin{align*}
    & \nabla F(\xUL) = \nabla_\xUL f\big(h^*(\xUL),\xUL\big) + \sfR_*(\xUL)^\top \nabla_h f\big(h^*(\xUL),\xUL\big)\\
    & \nabla F(\xUL') = \nabla_\xUL f\big(h^*(\xUL'),\xUL'\big) + \sfR_*(\xUL')^\top \nabla_h f\big(h^*(\xUL'),\xUL'\big)
\end{align*}
\begin{align*}
    & \|\nabla F(\xUL)-\nabla F(\xUL')\| \\
    & \leq \big\|\nabla_\xUL f\big(h^*(\xUL),\xUL\big) - \nabla_\xUL f\big(h^*(\xUL'),\xUL'\big)\big\| + \big\|\sfR_*(\xUL)^\top \nabla_h f\big(h^*(\xUL),\xUL\big) - \sfR_*(\xUL')^\top \nabla_h f\big(h^*(\xUL'),\xUL'\big) \big\| \\ \\
    & \|\nabla_\xUL f\big(h^*(\xUL),\xUL\big) - \nabla_\xUL f\big(h^*(\xUL'),\xUL'\big)\| \\
    & = \big\|\nabla_\xUL f\big(h^*(\xUL),\xUL\big) - \nabla_\xUL f\big(h^*(\xUL),\xUL'\big) + \nabla_\xUL f\big(h^*(\xUL),\xUL'\big) - \nabla_\xUL f\big(h^*(\xUL'),\xUL'\big)\big\|\\
    & \leq L_f \|\xUL-\xUL'\| + L_f \|h^*(\xUL)-h^*(\xUL')\| \\
    & \leq L_f (1+\sfC_0) \|\xUL-\xUL'\| \quad \text{(using Lemma \ref{lemma:Bnd-R*_Lip-h*}\ref{lemma:Lipschtz_h^*})} \\ \\
    & \big\|\sfR_*(\xUL)^\top \nabla_h f\big(h^*(\xUL),\xUL\big) - \sfR_*(\xUL')^\top \nabla_h f\big(h^*(\xUL'),\xUL'\big) \big\| \\
    & = \big\|\left(\sfR_*(\xUL)-\sfR_*(\xUL')+\sfR_*(\xUL')\right)^\top \nabla_h f\big(h^*(\xUL),\xUL\big) - \sfR_*(\xUL')^\top \nabla_h f\big(h^*(\xUL'),\xUL'\big) \big\| \\
    & = \big\|\left(\sfR_*(\xUL)-\sfR_*(\xUL')\right)^\top\nabla_h f\big(h^*(\xUL),\xUL\big) \big\| + \big\|\sfR_*(\xUL')^\top \big(\nabla_h f\big(h^*(\xUL),\xUL\big)-\nabla_h f\big(h^*(\xUL'),\xUL'\big) \big\| \\
    & \leq L_{\sfR_*} \Omega_f \|\xUL-\xUL'\| + \sfC_0 L_f (1+\sfC_0) \|\xUL-\xUL'\| \quad \text{(using Lemma \ref{lemma:Lipschtz_matrixR^*})}
\end{align*}
As a result, $L_F=L_{\sfR_*} \Omega_f+  L_f (1+\sfC_0)^2$. \Halmos
\endproof

\proof{Proof of Lemma \ref{lemma:Bnd-R*_Lip-h*}.} \ \\
\begin{enumerate}[label=\normalfont(\alph*)]
\item Consider general form of \eqref{eq:Jacb_Dyn-FixedPoint_v3} below for any given $\xUL$ and $0\leq\alpha\leq\overline\alpha$
\begin{align*}
     \sfR_{*}(\xUL) = (I - M_*(\xUL))^{-1} U_*(\xUL) \implies \|\sfR_{*}(\xUL)\| \leq \|\left(I - M_*(\xUL)\right)^{-1}\|\cdot \|U_*(\xUL)\|
\end{align*}
Use $\|\nabla_\xUL \nabla_{h} g(h,\xUL)\|\leq\Lambda_g^{\xUL h}$ (from Assumption \ref{assump:Bounded_LLCrossHessian}) and $\|B(h)\|\leq 1$ (from Lemma \ref{lemma:B_k,t factor}) to get $\|U_*(\xUL)\|\leq\alpha\Lambda_g^{\xUL h}$. Now, $(I - M_*(\xUL))^{-1} = \sum_{i=0}^\infty {M^i_*(\xUL)}$ implies $\|(I - M_*(\xUL))^{-1}\|\leq \sum_{i=0}^\infty \|M^i_*(\xUL)\|$. Note for a given $\xUL$, the general form of LMI \eqref{eq:LMI2_stability} (in the proof of Lemma \ref{lemma:Robust_Q-stability}) implies
\begin{align*}
\lambda\ \bfP - {M_{*}(\xUL)}^\top \bfP M_{*}(\xUL) \succeq 0 
\end{align*}
were $\bfP=C^{\bfP}\ \Diag(1/h^{*}(\xUL))$, $\lambda = {1-\eta\alpha/2}$. 
Let $y_{i+1}=M_*(\xUL) y_i$ then
\begin{align*}
   & y_i = M_*^i(\xUL) y_0,\quad  y_{i}^\top \bfP y_{i} \leq \lambda\ y_{i-1}^\top\bfP y_{i-1}\leq \hdots \leq \lambda^i y_{0}^\top\bfP y_{0} \\
   &\implies \min\left(\eig(\bfP)\right) \|M_*^{i}(\xUL) y_0\|^2 \leq \lambda^i \max\left(\eig(\bfP)\right) \|y_0\|^2 \\
   & \begin{aligned}[t] \implies &  \frac{\|M_*^{i}(\xUL) y_0\|}{\|y_0\|} \leq \lambda^{i/2} \sqrt{\frac{\max\left(\eig(\bfP)\right)}{\min\left(\eig(\bfP)\right)}}\quad \forall y_0\neq 0 \\
  \implies & \|M_*^{i}(\xUL)\| \leq \lambda^{i/2}/\sqrt{\nu^{\min}}
   \end{aligned}
\end{align*}
Therefore, 
\begin{align*}
    \|(I - M_*(\xUL))^{-1}\|\leq \frac{\sum_{i=0}^\infty \lambda^{i/2}}{\sqrt{\nu^{\min}}} = \frac{1/\sqrt{\nu^{\min}}}{1-\sqrt{\lambda}}
    \implies \|\sfR_*(\xUL)\|\leq \frac{\alpha\Lambda_g^{\xUL h}/\sqrt{\nu^{\min}}}{1-\sqrt{{1-\eta\alpha/2}}} 
\end{align*}
From Remark \ref{remark:JacbFP+B_k,*}\ref{remark:JacbFP-indp_alpha}, $\sfR_*(\xUL)$ does not depend on $\alpha$ then by taking $\sfC_0=\displaystyle\lim_{\alpha\to0^{+}} \frac{\alpha\Lambda_g^{\xUL h}/\sqrt{\nu^{\min}}}{1-\sqrt{{1-\eta\alpha/2}}} = \frac{4\Lambda_g^{\xUL h}}{\eta(\nu^{\min})^{1/2}}$ we have $\|\sfR_*(\xUL)\|\leq \sfC_0$.
\item For any $\xUL, \xUL'\in\calC$
\begin{align*}
    h^*(\xUL)-h^*(\xUL') = \int_{t=0}^1 \sfR_*(t\xUL+(1-t)\xUL')\cdot (\xUL-\xUL') dt \;\; \implies\;\; \|h^*(\xUL)-h^*(\xUL')\| \leq \sfC_0 \|\xUL-\xUL'\| 
\end{align*}
\end{enumerate}    
{This completes the proof.} \Halmos
\endproof

\proof{Proof of Lemma \ref{lemma:Lipschtz_matrixR^*}.}
For $\xUL,\xUL'\in\calC$
\begin{align*}
    & \|\sfR_{*}(\xUL)-\sfR_{*}(\xUL')\| \\
    & = \|\left(I-M_*(\xUL)\right)^{-1} U_*(\xUL)-\left(I-M_*(\xUL')\right)^{-1} U_*(\xUL')\| \\
    & \leq \|\left(I-M_*(\xUL)\right)^{-1}\|\cdot \|U_*(\xUL)-U_*(\xUL')\|+\|\left(I-M_*(\xUL)\right)^{-1}-\left(I-M_*(\xUL')\right)^{-1}\|\cdot \|U_*(\xUL')\|
\end{align*}
From Remark \ref{remark:JacbFP+B_k,*}\ref{remark:JacbFP-indp_alpha}, $\sfR_{*}(\xUL)$ does not depend on choice of $\alpha$. For ease of analysis, the lower-level step size is taken to be $\overline\alpha$ for both $\xUL,\xUL'$
\begin{align*}
    & \|U_*(\xUL)-U_*(\xUL')\| \\
    & = \overline\alpha\|\nabla_\xUL \nabla_{h} g(h^*(\xUL),\xUL)-\nabla_\xUL \nabla_{h} g(h^*(\xUL'),\xUL')\| \\
    & \leq \overline\alpha\|\nabla_\xUL \nabla_{h} g(h^*(\xUL),\xUL)-\nabla_\xUL \nabla_{h} g(h^*(\xUL),\xUL')\|+\overline\alpha\|\nabla_\xUL \nabla_{h} g(h^*(\xUL),\xUL')-\nabla_\xUL \nabla_{h} g(h^*(\xUL'),\xUL')\| \\
    & \leq \overline\alpha (\tau^\xUL_g+\tau_g^h \sfC_0) \|\xUL-\xUL'\|\quad \left(\text{using Lemma \ref{lemma:Bnd-R*_Lip-h*}\ref{lemma:Lipschtz_h^*}, Assumptions \ref{assump:Lipschtz_LL2ndOrder}, \ref{assump:Lipschtz_LL2ndOrder-x}}\right)
\end{align*}
Using bound on $\|\left(I-M_*(\xUL)\right)^{-1}\|$ from proof of Lemma \ref{lemma:Bnd-R*_Lip-h*}\ref{lemma:Bounded_R*} and $\|U_*(\xUL)\|\leq\overline\alpha\Lambda_g^{\xUL h}$ to get
\begin{align*}
    \|\sfR_{*}(\xUL)-\sfR_{*}(\xUL')\| 
    & \leq \frac{\overline\alpha (\tau^\xUL_g+\sfC_0 \tau_g^h)/\sqrt{\nu^{\min}}}{1-{\overline\lambda}^{0.5}} \|\xUL-\xUL'\|+\overline\alpha\Lambda_g^{\xUL h}\|\left(I-M_*(\xUL)\right)^{-1}-\left(I-M_*(\xUL')\right)^{-1}\|
\end{align*}
where $\overline\lambda={1-\eta\overline\alpha/2}$.
\newline Note that \ 
\newline $\begin{aligned}[t]
    \left(I-M_*(\xUL)\right)^{-1}-\left(I-M_*(\xUL')\right)^{-1} & = \left(I-M_*(\xUL)\right)^{-1}(I-M_*(\xUL')-I+M_*(\xUL))\left(I-M_*(\xUL')\right)^{-1} \\
    & = \left(I-M_*(\xUL)\right)^{-1}(M_*(\xUL)-M_*(\xUL'))\left(I-M_*(\xUL')\right)^{-1}
\end{aligned}$
\newline Therefore,
\begin{align*}
    \|\left(I-M_*(\xUL)\right)^{-1}-\left(I-M_*(\xUL')\right)^{-1}\| & \leq \|\left(I-M_*(\xUL)\right)^{-1}\|\cdot \|M_*(\xUL)-M_*(\xUL')\|\cdot \|\left(I-M_*(\xUL')\right)^{-1}\| \\
    & \leq \frac{1/\nu^{\min}}{\big(1-{\overline\lambda}^{0.5}\big)^2} \|M_*(\xUL)-M_*(\xUL')\|
\end{align*}
Recall $M_*(\xUL)-M_*(\xUL')=B(h^*(\xUL))\big[(1-\eta\overline\alpha)\Diag\left(1/h^*(\xUL)\right)-\overline\alpha \nabla^2_{h} g(h^*(\xUL),\xUL)\big] - B(h^*(\xUL'))\big[(1-\eta\overline\alpha)\Diag\left(1/h^*(\xUL')\right)-\overline\alpha \nabla^2_{h} g(h^*(\xUL'),\xUL')\big]$. Therefore,
\begin{align*}
    & \|M_*(\xUL)-M_*(\xUL')\| \\ 
    &\leq (1-\eta\overline\alpha) \|B(h^*(\xUL))-B(h^*(\xUL'))\|\cdot\|\Diag\left(1/h^*(\xUL)\right)\| \\
    & \quad + (1-\eta\overline\alpha) \|B(h^*(\xUL'))\|\cdot \|\Diag\left(1/h^*(\xUL)\right)  - \Diag\left(1/h^*(\xUL')\right)\| \\
    & \quad +\overline\alpha \| B(h^*(\xUL))\|\cdot \|\nabla^2_{h} g(h^*(\xUL),\xUL) - \nabla^2_{h} g(h^*(\xUL'),\xUL')\| + \overline\alpha\|B(h^*(\xUL)) - B(h^*(\xUL'))\|\cdot  \|\nabla^2_{h} g(h^*(\xUL'),\xUL')\|\\
    & \leq \left((1-\eta\overline\alpha)/\nu^{\min}+\overline\alpha L_g\right) L_B \|h^*(\xUL)-h^*(\xUL')\| + (1-\eta\overline\alpha) L_{1/h} \|h^*(\xUL)-h^*(\xUL')\|  \\ 
    & \quad +\overline\alpha \|\nabla^2_{h} g(h^*(\xUL),\xUL) - \nabla^2_{h} g(h^*(\xUL),\xUL')\| + \overline\alpha \|\nabla^2_{h} g(h^*(\xUL),\xUL') - \nabla^2_{h} g(h^*(\xUL'),\xUL')\| \\
    & \qquad \text{(using Lemma \ref{lemma:Lipschtz_matrixB}, Lemma \ref{lemma:Lipschtz_1/h}, $\|B(h^*(\xUL))\|\leq 1$ from Lemma \ref{lemma:B_k,t factor})} \\
    & \leq \left(1/\nu^{\min}+\overline\alpha L_g\right) L_B \sfC_0 \|\xUL-\xUL'\| + L_{1/h}\sfC_0 \|\xUL-\xUL'\| + \overline\alpha(\rho_g^\xUL+\rho_g^h \sfC_0) \|\xUL-\xUL'\|
    \\ & \qquad \qquad \text{(using Lemma \ref{lemma:Bnd-R*_Lip-h*}\ref{lemma:Lipschtz_h^*}, Assumptions \ref{assump:Lipschtz_LL2ndOrder}, \ref{assump:Lipschtz_LL2ndOrder-x})}
\end{align*}
Finally
\begin{align*}
    \|\sfR_{*}(\xUL)-\sfR_{*}(\xUL')\| & \leq \frac{\overline\alpha (\tau^\xUL_g+\sfC_0 \tau_g^h)/\sqrt{\nu^{\min}}}{1-{\overline\lambda}^{0.5}} \|\xUL-\xUL'\| \\
    & \quad +\frac{\overline\alpha\Lambda_g^{\xUL h}/\nu^{\min}}{\big(1-{\overline\lambda}^{0.5}\big)^2} \left[\left(1/\nu^{\min}+\overline\alpha L_g\right) L_B \sfC_0 + L_{1/h}\sfC_0 + \overline\alpha(\rho_g^\xUL+\rho_g^h \sfC_0)\right] \|\xUL-\xUL'\|\\
     & = L_{\sfR_*} \|\xUL-\xUL'\|
\end{align*}
The proof is complete.\Halmos
\endproof

\section{Auxiliary lemmas and proofs.}\label{appendix:AuxProofs}
\proof{Proof of Lemma \ref{lemma:GlobBndErr}.}
Since for each $i=1,\hdots,\sfW$, we have $h^{k,*}_i\in\widetilde \calH_i$ and $h^{k,0}\in\widetilde \calH_i$ from Lemma \ref{lemma:Iters_SmplxIntr}.
\begin{align*}
    \txteps_h^{k,0} & = \sum_{i=1}^\sfW \sfD_\KL(h^{k,*}_i,h^{k,0}_i) \\
    & = \sum_{i=1}^\sfW \sum_{j=1}^{\sfq_\ell^i} h^{k,*}_{i,j} \ln\frac{h^{k,*}_{i,j}}{h^{k,0}_{i,j}} - \left(h^{k,*}_{i,j}-h^{k,0}_{i,j}\right)
    \begin{aligned}[t]
    & = \sum_{i=1}^\sfW \sum_{j=1}^{\sfq_\ell^i} h^{k,*}_{i,j} \ln(h^{k,*}_{i,j}) - h^{k,*}_{i,j} \ln(h^{k,0}_{i,j})\\
    & \leq \nu^{\min} \sum_{i=1}^\sfW \sum_{j=1}^{\sfq_\ell^i} \ln(h^{k,*}_{i,j}) +\sfW \ln(1/\nu^{\min}) \\ 
    & \leq \nu^{\min} \sum_{i=1}^\sfW \sfq_\ell^i \ln(1/\sfq_\ell^i) + \sfW \ln(1/\nu^{\min})\\
    & \left\{\because\ \max_{\substack{\mathbf{1}^\top\rmx=1,\\ \rmx\geq0}} \sum_{j=1}^{\sfq_\rmx} \ln(\rmx_{j}) = \sfq_\rmx \ln(1/\sfq_\rmx) \ \text{ at } \ \rmx^*=1/\sfq_\rmx \right\}
    \end{aligned}
\end{align*}
The proof is complete. \Halmos
\endproof

\begin{lemma}\label{lemma:Lipschtz_matrixB}
    The matrix function $B(h)=\blk\big(B^1,\hdots,B^\sfW\big)$ (where $B^i = \Diag\big(h_i\big) - h_i {h_i}^\top$) is $L_B$--Lipschitz continuous for $h\in\calH$.
\end{lemma}
\proof{Proof of Lemma \ref{lemma:Lipschtz_matrixB}.}
For each $i\in[\sfW]$, $${\partial B(h)}/{\partial h_{i,j}} = \blk\big(0,\hdots,{\partial B^i}/{\partial h_{i,j}},\hdots,0\big),\quad j\in[\sfq^i_\ell]$$ where ${\partial B^i}/{\partial h_{i,j}}=\delta^j_i{\delta^j_i}^\top-\delta^j_i {h_i}^\top - {h_i}{\delta^j_i}^\top$ and $\delta^j_i$ is an indicator vector defined as if $j=j'$ then $[\delta^j_i]_{j'} = 1$, otheriwse $[\delta^j_i]_{j'} = 0$. For $h,h'\in \calH$
\begin{align*}
    B(h)-B(h') & = \sum_{i,j}  \int_{t=0}^1 {\partial B(h)}/{\partial h_{i,j}}\Big|_{th_i+(1-t)h'_i} (h_{i,j}-h'_{i,j}) dt = \blk\big(Q^1,\hdots,Q^\sfW\big)
\end{align*}
where 
\begin{align*}
    Q^i & = \int_{t=0}^1 \sum_{j\in [\sfq^i_\ell]} \frac{\partial B(h)}{\partial h_{i,j}}\Big|_{th_i+(1-t)h'_i} (h_{i,j}-h'_{i,j}) dt\\
    & = \Diag(h_i-h'_i) - \frac{1}{2}\big[(h_i-h'_i)(h_i+h'_i)^\top+(h_i+h'_i)(h_i-h'_i)^\top\big]
\end{align*}
Since the spectral norm of a rank-1 matrix $uv^\top$ is $\|uv^\top\|=\|u\|\cdot\|v\|$.
As a result,
\begin{align*}
    \|Q^i\| \leq \|h_i-h'_i\|_\infty + \|h_i-h'_i\|\cdot\|h_i+h'_i\| & \leq 3\ \|h_i-h'_i\| \quad (\text{using } \|\cdot\|_\infty\leq\|\cdot\|_2 \text{ and } h,h'\in\calH)\\
    & \leq 3\ \|h-h'\|
\end{align*}
Therefore, $\|B(h)-B(h')\big\| \leq L_B\ \|h-h'\|$ where $L_B=3$. \Halmos
\endproof

\begin{lemma}\label{lemma:Lipschtz_1/h}
    The matrix function $\Diag(1/h)$ is $L_{1/h}$--Lipschitz continuous for $h\in\widetilde\calH$.
\end{lemma}
\proof{Proof of Lemma \ref{lemma:Lipschtz_1/h}.}
For each $i\in[\sfW]$, ${\partial \Diag(1/h)}/{\partial h_{i,j}} = \Diag\big(0,\hdots,-1/(h_{i,j})^2,\hdots,0\big),\ j\in[\sfq^i_\ell]$. 
For any $h,h'\in\widetilde \calH$
\begin{align*}
    & \begin{aligned}
    \Diag\left(1/h\right)-\Diag\left(1/h'\right) & = \sum_{i,j}  \int_{t=0}^1 \frac{\partial \Diag(1/h)}{\partial h_{i,j}}\Big|_{th_{i,j}+(1-t)h'_{i,j}} (h_{i,j}-h'_{i,j}) dt \\
    & = \int_{t=0}^1 \Diag(-1/(th+(1-t)h')^2) \cdot (h-h') dt
    \end{aligned}\\ \\
    \implies & \big\|\Diag\left(1/h\right)-\Diag\left(1/h'\right)\big\| \leq L_{1/h}\ \|h-h'\| 
\end{align*}
where $L_{1/h}=1/(\nu^{\min})^2$ since $th+(1-t)h'\in\widetilde\calH$ for $0\leq t\leq 1$ and $\|\Diag\left(-1/h^2\right)\|\leq1/(\nu^{\min})^2\ \forall h\in\widetilde\calH$. 
\Halmos
\endproof

\end{APPENDICES}

\end{document}